\documentclass[sn-mathphys-num]{sn-jnl}% Math and Physical Sciences Numbered Reference Style 
\usepackage{anyfontsize}
\usepackage{amsmath}
\usepackage{amsfonts}
\usepackage{amssymb}
\usepackage{epsfig}
\usepackage{dsfont}
\usepackage{mathrsfs}
\usepackage{bbm}
\usepackage{bm}

\usepackage[linesnumbered,titlenotnumbered,ruled,vlined]{algorithm2e}
\usepackage{float}
\usepackage{makecell} 
\usepackage{boldline}
\usepackage{diagbox}
\usepackage{threeparttable}
\usepackage[dvipsnames]{xcolor}
\usepackage{adjustbox}
\usepackage{booktabs}
\usepackage{hyperref}

\usepackage{enumerate}
\allowdisplaybreaks[4]
\usepackage[title]{appendix}
\usepackage{mathabx} %wide
\usepackage{subfigure}
\usepackage{caption}
\usepackage[normalem]{ulem}

\newtheorem{theorem}{Theorem}[section]
\newtheorem{lemma}{Lemma}[section]
\newtheorem{proposition}{Proposition}[section]
\newtheorem{example}{Example}[section]
\newtheorem{remark}{Remark}[section]
\newtheorem{definition}{Definition}[section]
\newtheorem{assumption}{Assumption}[section]
\newtheorem{corollary}{Corollary}[section]
\newcommand{\dd}{\mathsf {d\kern -0.07em l}} %	nested

\newcommand{\bgeqn}{\begin{eqnarray}}
\newcommand{\edeqn}{\end{eqnarray}}
\newcommand{\bgeq}{\begin{eqnarray*}}
\newcommand{\edeq}{\end{eqnarray*}}
\newcommand{\bec}{\begin{center}}
\newcommand{\enc}{\end{center}}
\newcommand{\R}{{\rm I\!R}}

\newcommand{\inmat}[1]{\mbox{\rm {#1}}}

\renewcommand{\Box}{\framebox{\rule{0.3em}{0.0em}}}

\newcommand{\be}{\begin{equation}}
\newcommand{\ee}{\end{equation}}

\numberwithin{equation}{section}
\numberwithin{definition}{section}

\newcommand{\st}{{\rm s.t.}}
\newcommand{\CVaR}{{\rm CVaR}}

\newcommand{\proof}{\noindent\textbf{Proof.}\;\;}

\newcommand{\ind}{\mathds{1}}
\newcommand{\1}{\textbf{1}}

\newcommand{\essinf}{\text{essinf}\,}
\newcommand{\esssup}{\text{esssup}\,}
\def\bbe{{\mathbb{E}}}

\renewcommand{\theequation}{\thesection.\arabic{equation}}

\begin{document}

\title[Multistage Robust Average Randomized Spectral Risk Optimization]{Multistage Robust Average Randomized Spectral Risk Optimization\footnote{This work is supported by CUHK start-up grant and EPSRC (UK) Grant EP/V008331/1.}}

%%=============================================================%%
%% GivenName	-> \fnm{Joergen W.}
%% Particle	-> \spfx{van der} -> surname prefix
%% FamilyName	-> \sur{Ploeg}
%% Suffix	-> \sfx{IV}
%% \author*[1,2]{\fnm{Joergen W.} \spfx{van der} \sur{Ploeg} 
%%  \sfx{IV}}\email{iauthor@gmail.com}
%%=============================================================%%

\author*[1]{Qiong Wu}\email{qiongwu001@outlook.com}

\author[2]{Huifu Xu}\email{h.xu@cuhk.edu.hk}
% \equalcont{These authors contributed equally to this work.}

\author[2]{Harry Zheng}\email{h.zheng@imperial.ac.uk}
% \equalcont{These authors contributed equally to this work.}

% \cortext[cor1]{Corresponding author}

\affil*[1]{\orgdiv{Department of Systems Engineering and Engineering Management}, \orgname{The Chinese University of Hong Kong}, \orgaddress{\street{Shatin, N.T.}, \city{Hong Kong}, \postcode{999077}, \country{China}}}

\affil[2]{\orgdiv{Department of Mathematics}, \orgname{Imperial College}, \orgaddress{\city{London}, \postcode{SW7 2AZ}, \country{United Kingdom}}}

%%==================================%%
%% Sample for unstructured abstract %%
%%==================================%%

\abstract{In this paper, we revisit the multistage spectral risk minimization models %investigated
proposed 
by Philpott et al.~\cite{PdF13} and  Guigues and R\"omisch \cite{GuR12}  
but with some new focuses. We consider a situation where the decision maker's (DM's) risk preferences may be state-dependent or even inconsistent at some states, and consequently 
there is not a single deterministic spectral risk measure (SRM) which can be used to represent the DM's preferences at each stage. 
We adopt 
the recently introduced 
%a 
average randomized SRM (ARSRM) 
% which is recently introduced 
(in \cite{li2022randomization})
to describe the DM's overall risk preference
%risk preferences
%and then use the average of the RSRM (ARSRM) to 
%model the DM's preference 
at each stage. 
To solve the resulting multistage ARSRM (MARSRM) problem,
%program,
we apply the well-known 
stochastic dual dynamic programming (SDDP) 
method
%approach
%to solve the resulting multistage risk minimization problem
%which allow us to develop  
which 
%allows us to construct 
%upper and lower approximations 
%at each stage
%to produce 
generates
a sequence of 
lower and upper bounds in an iterative manner. Under some moderate conditions, 
we prove that the optimal solution can be
% we 
% demonstrate the optimal solution to be 
found in a finite number of iterations. 
%respectively.
The 
%approach 
%multistage 
MARSRM model 
% provides 
%simplifies 
%a generalization 
%of 
generalizes the one-stage ARSRM
%as well as
and 
%a simplification of
simplifies the existing 
multistage 
%risk minimization models where the DM's risk preference is 
state-dependent
preference robust model
\cite{liu2021multistage}, %{\color{red}it
%we describe 
%covers
%and 
% it also subsumes 
while also encompassing
the mainstream 
multistage risk-neutral and risk-averse  optimization models \cite{GuR12,PdF13}. 
In the absence of complete information on the probability distribution of the 
DM's random preferences, we propose to use distributionally robust ARSRM (DR-ARSRM) to describe the DM's preferences at each stage.
% in the risk spectrum.
%, we propose a distributionally robust model.
%Consequently the multistage problem is solved based on the average random spectral risk measure (ARSRM). 
%???}. and 
We 
% discuss in detail 
detail
computational schemes
for solving both MARSRM and DR-MARSRM. 
%the multistage optimization problems with ARSRM based on the moment-type ambiguity set of the probability distribution of random preference parameter.
Finally, 
%a  numerical example is provided to 
%show 
we examine the 
%efficiency
performance of 
%the multistage 
MARSRM and DR-MARSRM 
by applying them to 
%a stage-wise independent 
%case with application to 
an asset allocation problem with transaction costs and 
%comparison 
compare them with 
%the common 
standard risk neutral and risk averse multistage linear stochastic programming (MLSP) models. }


%%\pacs[JEL Classification]{D8, H51}

%%\pacs[MSC Classification]{35A01, 65L10, 65L12, 65L20, 65L70}

\maketitle

\section{Introduction}\label{sec1}

Spectral risk measure (SRM) introduced 
by Acerbi \cite{acerbi2002spectral} is defined as the weighted average of the quantile function of 
a random variable over $[0,1]$.
It subsumes a number of important risk measures
including conditional value-at-risk (CVaR),
essential supremum, Wang's proportional Hazard transform,
Gini's risk measure and convex combination of expected value and CVaR.
% risk measure.
SRM is also known as a distortion risk measure when the weighting function (also known as risk spectrum) is the derivative of a distortion 
function. 
% The latter 
This concept
is rooted in Yaari's dual theory of choice \cite{yaari1987dual} and 
is widely used in actuarial science \cite{denneberg1990distorted,wang1995insurance,tsanakas2003risk}. SRM has several advantages: it is a coherent risk measure \cite{acerbi2002spectral} and the set of SRMs is dense in the cone of the law-invariant coherent risk measures \cite{li2022randomization}; it can be represented in terms of CVaR \cite{acerbi2002portfolio} and its risk spectrum can be effectively used to describe a decision maker's (DM's) risk preference \cite{WaX20,guo2022robust}. 
Guigues and R\"omisch \cite{GuR12} 
% seem to be 
were among
the first to
apply SRM in multistage 
decision making. By restricting the DM's decision-making problem at each stage to a stochastic linear program with recourse, they use the well-known stochastic dual dynamic programming method (SDDP) to solve the multi-stage SRM-based risk minimization problem. Philpott et al.~\cite{PdF13}  apply a specific SRM which is a convex combination of the expected value and CVaR to 
a class of multistage stochastic linear programs in which at each stage the spectral risk of future costs
is to be minimized. Moreover, they apply the SDDP method to solve the multistage risk minimization problem and demonstrate 
efficient methods to obtain upper and lower bounds of the optimal value.
% how an upper and a lower bound
% of the optimal value may be obtained efficiently.

% A general computational approach based on dynamic programming is derived that can be shown to
% converge to an optimal policy. By computing an upper approximation to future cost functions, we can evaluate an upper
% bound on the cost of an optimal policy, and an lower approximation delivers a lower bound. The approach we describe
% is particularly useful in sampling-based algorithms, and a numerical example is provided to show the efficacy of the
% methodology when used in conjunction with stochastic dual dynamic programming

In this paper, we 
extend this line of research
% take on the stream of research 
but with slightly different focuses.
% {\color{red} First}, 
In  both \cite{GuR12}
%, the authors use a single SRM to measure the  
and \cite{PdF13}, the DM's risk preference at each stage is described by a single deterministic 
SRM which is independent of state and stage. 
While this simplifies the derivation of efficient 
numerical methods for solving the multistage program,
% the SRM models 
it
may not 
% cover
account for
the cases where the decision maker's risk preference is state or stage dependent. One way to 
%address 
%{\color{red} this} 
% {\color{blue} 
fill the gap is to introduce 
a state-dependent SRM at each stage and then apply SDDP to solve the resulting program. 
For instance, Liu et al.~\cite{liu2021multistage}
%, the authors 
%propose 
consider a multistage 
expected utility maximization 
%utility preference robust
%model
problem where the DM's utility at each stage is state-dependent 
%{\color{blue} 
and 
at each state
the DM's 
%risk 
preference over gains and losses can be described by 
Von Neumann-Morgenstern expected utility.
% is consistent, 
 % and
In the case that information on the state-dependent 
utility of the DM is incomplete, they
propose 
 %develop 
a multistage utility 
preference robust model where 
 % }
the DM's decision is based on the worst-case utility function from a set of plausible utility functions. 
% {\color{blue} 
To solve the multistage maximin problem, they require the ambiguity set
of state-dependent utility functions to have some rectangularity structure
which ensures time consistency of optimal decisions and dynamic formulations of the problem, enabling 
the application of SDDP.
% one to apply SDDP for solving the problem.
% Under some moderate conditions they demonstrate that the optimal decisions are time consistent and the multi-stage PRO model can be solved by SDDP method fairly efficiently. 
% }
While their 
%work
model effectively addresses the issues of 
preference ambiguity and state-dependence, 
%successfully, 
it does not cover the case that the DM's preferences at 
each fixed 
%a 
state may be inconsistent. 
Inconsistency occurs for various reasons, 
see e.g. \cite{allais1953comportement, armbruster2015decision, bertsimas2013learning, delage2015robust, hu2022distributionally}. 
% {\color{red} HX: add Hu,Zhang,Xu and Zhang's paper and Alais paradox paper}
% {\color{blue} 
Moreover, like 
%all
other robust optimization models, the worst-case argument might be overly conservative particularly in finance. Furthermore, the mathematical structure of the multistage maximin problem is intrinsically 
%complicated 
complex, making 
% it difficult for 
% %in terms of
% for 
numerical resolution challenging.
%computation, 
Thus some kind of simplification 
of the model might be desirable.
% The disadvantage is that the SRM 
% will be historic path dependent and the
% resulting mathematical program is very complex and difficult to solve.  
%With all these in mind, here we adopt a slightly simpler approach where

In this paper, we aim to address these issues in the context of
multistage spectral risk minimization.
Specifically, 
we adopt  a stage-dependent 
average randomized SRM (ARSRM)
%to 
%which captures 
to describe the DM's 
%varying risk 
state-dependent preferences at each stage. 
%at different states.
ARSRM is recently proposed by 
Li et al.~\cite{li2022randomization} in one-stage decision-making. The key idea
is to consider a randomly parameterized SRM (RSRM) to characterize 
a DM's ``random'' preferences and then use the average (expected value)
of the RSRM to represent the DM's risk preference.
% Specifically, we introduce a parameter in the definition of risk spectrum (the weighting function) and then randomize the parameter. 
% %The randomization is 
% {\color{red}
% Specifically, 
% the random parameter 
% may follow the probability distribution of the preference states at a stage. Then we can use the mean value of the randomized SRM to describe the DM's risk preference at each stage.
% }
% Of course, to further simplify the model and the follow-up computational scheme, we may 
% assume that the distribution of the random parameter is dependent on stage but is independent on state. 
% The randomization approach is recently introduced by Li et al.~\cite{li2022randomization} in static stochastic decision-making and this is the first attempt to apply 
% it to multistage decision-making.
%{\color{red} 
One of the main challenges is to define the random parameter:
is it stage-dependent (vs independent).
This distinction determines whether
% meaning that 
the DM's 
random preferences depends (vs does not depend) on the stage variable and/or previous stages as well? 
Additionally, Is information on the probability distribution of the random parameter complete (vs incomplete)? What can we do if the information is incomplete?
Another main challenge is to develop efficient computational methods for solving the multistage decision-making problem with
stage-dependent (or independent) ARSRM and its distributionally robust counterpart in the absence of complete information of the distribution of the random parameter. 
% Like \cite{li2022randomization}, we 
We
propose to use SDDP but with significant simplifications. 
The main contributions 
%of  are 
of this paper can be summarized 
as follows.
% }

% In the literature on the multistage framework, it is commonly assumed that the underlying distribution of the uncertain parameter is known, while these assumptions are not entirely realistic.
% In the proposed multistage decision making with ARSRM, if there is absence of complete information of the probability distribution of the preference parameter, the modeller may consider a robust formulation of ARSRM.
% In this paper, we will also focus on the a situation where the probability distribution of DM's random preference parameter is ambiguous but it is possible to use partially available information to construct an ambiguity set which contains the true probability distribution with high confidence.
% The main challenge to be tackled is to develop efficient computational schemes for constructing the ambiguity set and solving the resulting multistage DRO model.

% Second, we apply the SDDP method to solve the multistage risk minimization problem based on stage-dependent randomized SRM.
% %Stochastic dual dynamic programming method 
% SDDP is a Monte Carlo sampling-based variant of Benders decomposition method which is first proposed by Pereira and Pinto \cite{pereira1991multi} in the risk-neutral setting, and has been employed successfully in a wide range of applications with large sample space, see e.g. \cite{homem2011sampling,rebennack2011stochastic}.
% {\color{red} Please add the unique challenges when SDDP is applied to our model in this paper. }
%}

\begin{itemize}
    \item %First,
    % {\color{blue}
    \textbf{Modelling}. 
    We propose a multistage risk minimization model
    where the DM's risk preferences at each stage 
    %is
    are
    described 
    by a stage-dependent ARSRM. 
    In comparison with previous models
    \cite{GuR12,PdF13}, the new model accommodates
    both stage dependence of the DM's risk preferences and randomness of such preferences. When the distribution of the random parameter in ARSRM is ambiguous, we propose a distributionally robust model 
    %where
    with the ambiguity set 
    being constructed by moment-type conditions. 
%     {\color{red} The ambiguity
% set at each stage is updated from the previous stage by
% the change of the mean value and variance, which means that
% the randomness of the DM's preference at the current stage depends on the average and standard deviation of the preferences at the previous stage HX: Apparently this is not the case}. 
The new model extends the work of 
%    In comparison with 
\cite{liu2021multistage}
%, our new model
by incorporating
% allowing 
random preferences and simplifies
the model by removing 
%state
historical path-dependence 
in describing the DM's preferences.
  % }  
        
%     PRO model where the DM's risk preferences at different stages are characterized by stage-dependent average randomized spectral risk measure.
% In the case that the probability distribution of the preference parameter is ambiguous, we propose the preference distributionally robust multistage decision making model.
% We define the moment-type ambiguity set of preference parameter from the past and known mean and variance and the optimal policy is obtained based on the worst-case ARSRM of the cost-to-go functions at different stages in the ambiguity set.

\item \textbf{Computational method}. 
We propose to solve the
%proposed https://cuhk.zoom.us/saml/mobile_success?token=ec003cf#success
multistage 
%decision-making 
%ARSRM based 
decision-making problem with ARSRM as an objective at each stage
by  using the well-known SDDP method. 
In the case that 
%{\color{blue}
the constructed scenario tree is stagewise independent, and
the probability distribution 
of the random preference parameter is known,
we propose an iterative scheme where
a single-cut lower approximation of the stagewise subproblem is constructed 
% to 
% form a lower approximation 
% of the Bellman's function 
% at each stage 
from the last stage, 
% and subsequently a lower bound for the original minimization 
% problem. 
thereby providing a lower bound for the original minimization problem.
Likewise, in each iteration,
%the outer 
%a lower approximation with a single cut 
%and 
%inner
% {\color{red} HX: Please rewrite the rest of the paragraph:
% an upper approximation at each stage 
% in each iteration to obtain 
% %the 
% a lower bound and 
% {\color{red}accurate??} 
% upper bound through  cutting plane and dual representation of CVaR when the random prospect is stage-independent.
% When the probability distribution of the preference parameter is ambiguous, 
% we derive a multi-cut version 
% %outer 
% of lower approximation and 
% %inner
% upper approximation to the equivalent formulation of $\CVaR$ with the $\VaR$ representative variable and the Lagrange dual.
% }
% {\color{blue}
we construct %iteratively
%the 
an upper approximation at each stage 
and obtain 
%the 
a deterministic upper bound of the original 
multistage minimization problem
by developing a lower convex envelope of 
the previous upper bounds.
When the probability distribution of the preference parameter is ambiguous, we derive  lower and upper bounds in a per-scenario manner. That is, we build 
%the 
a multi-cut version of  lower approximations 
and construct lower convex envelope of upper bounds
 per scenario.
% 
%\item \textbf{Error bound}. 
We prove the convergence of the upper and lower bounds 
for the true optimal values of the stagewise subproblems
obtained from the 
%inner and outer 
approximations.
%are valid.
%Furthermore, 
%Moreover, 
Specifically, we prove the finite time convergence of the SDDP algorithm to the optimal solution. 
Furthermore, we derive 
%the 
error bound of the optimal value of the multistage decision-making problem incurred by
the step-like approximation of the risk spectrum.

\item \textbf{Application}. We apply the proposed model and algorithm to a multistage asset allocation problem with transaction costs.
We examine the convergence of the optimal value with the increasing iterations and compare the multistage decision-making problem with ARSRM, expectation (risk neutral case), and the combination of expectation and $\CVaR$ (in mild/strong risk averse case).
We also carry out numerical tests on the preference distributionally robust multistage asset allocation problem.
% and examine the out-of-sample performance of the model.

\end{itemize}

The paper is organized as follows.
In Section 2, we 
%introduce
% recall
revisit
the definition of the randomized spectral risk measure and present tractable formulations for computing the 
%computational scheme of 
general average random spectral risk measures.
In Section 3, we 
%formulate the 
propose a multistage stochastic 
%linear programs 
minimization model 
%with a risk averse 
where the decision maker's 
%whose 
risk preferences can be characterized by 
%an average random spectral risk measures 
stage-dependent ARSRM
and discuss details as to how to use SDDP method for solving 
the problem.
%study the lower approximation and the upper approximation.
% discuss the computation of the general average random spectral risk measures.
In Section 4, we take a step further to
%discuss the formulation in 
consider the case that the probability distribution of the random preference parameter 
in the definition of ARSRM is 
%unknown
ambiguous, 
and propose a distributionally robust model which
mitigates the model risk arising from ambiguity 
of the distribution. We then 
detail the application of SDDP to 
% discuss 
% how
% the SDDP may be applied to 
solve the multistage minimax robust optimization problem.
%f the corresponding distributional preference robust multistage problem, tractable formulations of inner and lower approximations, convergence and error bound.
% study the lower approximation and the upper approximation algorithm.
In Section 5, we report
%presents
numerical test results of the proposed models and computational schemes. Finally, 
%applying the methodology to a asset allocation problem.
we conclude the paper with some closing
%discussion 
remarks in Section 6.

Throughout the paper, we use the following notation. 
By convention, we use $\R^n$ to denote the Euclidean space of $n$-dimensional vectors, 
$\R^n_+$ to denote those with non-negative components and 
$\R^n_{++}$ to denote those with positive components. 
 For any two vectors $a,b\in \R^n$, we write for simplicity of notation 
 $ab$ to denote their scalar product and
for two matrices $A$ and $B$, we write $A\cdot B$ for the Frobenius product which is the trace of $A^{\top}B$ where $A^{\top}$ denotes the transpose of $A$.
We write $[T]$ for a discrete set $T=\{1,\cdots,T\}$ and 
$|T|$ for the cardinality of the set.
For a random variable mapping from $(\Omega, {\cal F}, \mathbb{P})$ to $\R$, we use $\text{essinf} \; X$ an 
$\text{esssup} \; X$ to denote respectively the essential infimum and essential supremum of $X$, i.e, 
$\text{essinf} \; X = \sup\{r: \mathbb{P}(X<r)=0\}$ and 
$\text{esssup} \; X = \inf\{r: \mathbb{P}(X>r)=0\}$.

\section{Spectral Risk Measure and Its Randomization}
{\color{black}
% Risk-averse stochastic optimization problem have been studies for many years.
% The concept of coherent risk measures (CRM), introduced in \cite{ADEH99} is now widely accepted for time-static risk management and finance.
% Many risk measures are known to satisfy the four properties of CRM but there is not obvious advice to which particular risk measure an decision maker (DM) may use.
% The well-known representation of CRM is in term of weighted value-at-risk, that is, a real valued function $\rho$ is a coherent risk measure if and only if 
Let $\xi: (\Omega, {\cal F}, \mathbb{P})\to \R$ be a random variable representing financial losses.
The spectral risk measure is defined as 
\begin{equation}
\label{eq:SRM}
   \inmat{(SRM)} \quad  \rho(\xi) = 
    \int_0^1 F_X^{\leftarrow} (z) \sigma (z) d z,
\end{equation}
where $F_{\xi}^{\leftarrow} (z)$ is the quantile function of $\xi$ and $\sigma:[0,1)\to\R_+$ is a non-decreasing and right continuous function with $\int_0^1 \sigma(z) d z =1$.
%{\color{red} 
In this formulation, 
$\sigma$ plays a role of weighting which
reflects the decision maker's view on the losses at different levels of confidence. It is therefore known as risk spectrum.
By \cite[Theorem 1]{acerbi2002portfolio},
$\rho$ is a coherent risk measure if and only if $\sigma$ is non-decreasing.
% $\sigma$ is known as risk spectrum.
%) and $\mathcal{M}$ is a set of risk spectra.
When $\sigma$ takes some specific forms, the SRM recovers 
a number of important risk measures including CVaR, essential supremum, Wang's proportional hazard transform, Gini's risk measure and convex combination of the mean and CVaR. 
%{\color{blue} 
Since $\sigma$ may be viewed as the derivative of a convex function, SRM is also known as distortion risk measure (DRM), see \cite{acerbi2002spectral}.  
The latter is widely used in actuarial science and finance; see, e.g., \cite{denneberg1990premium,wang1995insurance, acerbi2002spectral,PiS15} and monograph \cite{SDR09}.
%In the case where the true risk spectrum is ambiguous, we have the following assumption:
% \begin{assumption}
% \label{A-RSRM}
%     The DM's risk preference 
%     % at every stage 
%     can be represented by an SRM but there does not exist a deterministic risk spectrum $\sigma$ such that the DM's risk preference can be described by $\rho_{\sigma}$.
% \end{assumption}
% Specifically, Wang and Xu \cite{WaX20} propose a robust SRM
% \begin{equation}
% \label{eq-CRM}
%     \rho(X):=\sup_{\sigma\in\mathcal{M}} \rho_\sigma(X) = \sup_{\sigma\in\mathcal{M}}
%     \int_0^1 F_X^{\leftarrow} (z) \sigma (z) d z,
% \end{equation}
% where $\mathcal{M}$ is an ambiguity set of plausible risk spectra. 

\subsection{Randomization of SRM}

Li et al.~\cite{li2022randomization} 
%take a step further to
consider the case that a DM's spectral risk preference is not only 
ambiguous but also inconsistent,
which means that there is no single $\sigma$ such that $ \rho_{\sigma_s}$ can be used to describe the DM's risk preference. Consequently, they introduce a randomized SRM
\begin{equation}
\label{eq-RSRM}
    \text{(RSRM)\quad} \rho_{\sigma_s}(\xi) = 
    \int_0^1 F_{\xi}^{\leftarrow} (z) \sigma (z,s) dz,
\end{equation}
where $s\in\R^{d_s}$ is a random parameter. 
% \begin{assumption}
% \label{A-RSRM}
In this set-up, the DM's risk preference 
% at every stage 
can be represented by an SRM but there does not exist a deterministic risk spectrum $\sigma$ such that the DM's risk preference can be described by $\rho$. 
%\end{assumption}
In the case that 
the probability distribution of $s$ is known, they consider
%Let $Q$ denote the probability distribution of $s$, 
the average of RSRM: 
%$\rho_{\sigma(\cdot,s)}(X)$
%can be calculated by 
\begin{equation}
\label{eq-ARSRM}
    \text{(ARSRM)\quad} \rho_{Q}(\xi) := \bbe_{Q}\left[\rho_{\sigma_s}(\xi) \right] 
    = \bbe_{Q}\left[ \int_0^1 F^{\leftarrow}_{\xi}(z) \sigma(z,s) dz \right],
\end{equation}
where $Q$ denotes the probability distribution of $s$. Note that since $\rho_{Q}(\xi)$ can also be written as
% equals to
$\int_0^1 F^{\leftarrow}_{\xi}(z) \bbe_{Q}\left[ \sigma(z,s)\right] dz 
$,
and $\bbe_{Q}\left[ \sigma( \cdot, s)\right] $ is a deterministic risk spectrum,  $\rho_{Q}(\xi)$
is also an SRM. The key point here is that
$\bbe_{Q}\left[ \sigma(\cdot,s)\right]$ does not really 
capture the DM's true risk preference (which is random), rather it is the ``average'' from modeller's perspective. In other words, given uncertainty of the DM's risk preference, it might be reasonable to suggest that 
$\bbe_{Q}\left[ \sigma(\cdot, s)\right]$ represents 
the DM's measure of risk in the next decision making stage.

Since $\mathbb{E}_Q\left[\sigma(\cdot,s)\right]: [0,1)\to \R_+$ and  
$\int_0^1\mathbb{E}_Q\left[\sigma(t,s)\right]dt=1$,
then $\mathbb{E}_Q\left[\sigma(\cdot,s)\right]$
is a risk spectrum, and  $\rho_{Q}(X)$
is a deterministic SRM. 
%{\color{blue}
Thus, any DM whose  
risk preferences 
can be represented by ARSRM 
%{\em virtually}
%as if
%(s)he  was a
is ``equivalent'' to a deterministic
SRM decision maker from mathematical modelling perspective.
%}
Unfortunately, in practice, we (as a modeller) 
are unable to obtain $\mathbb{E}_Q\left[\sigma(\cdot,s)\right]$ via preference elicitation when the DM's preferences are varying/inconsistent.
%At this point,
Moreover, it might also be helpful to clarify 
the motivations behind ARSRM in (\ref{eq-ARSRM}). 
% There are at least two reasons why we consider average rather than adopt a risk measure $\phi(\cdot)$ to measure
% the risk of random variable $\rho_{\sigma(\cdot,s)}(X)$.
First, it is the modeller  that determines whether to choose $\rho_{Q}(X)$ or 
$\phi\left(\rho_{\sigma(\cdot,s)}(X)\right)$ (where $\phi$ is a risk measure)
to measure the risk arising from  uncertainty of
$\rho_{\sigma(\cdot,s)}(X)$. The ``modelling selection''  is not necessarily related to the DM's risk preference, rather \underline{it is the modeller's ``taste''} as to which risk measure might be most appropriate to 
capture the DM's risk preference in future, see \cite[Example 1.1]{li2022randomization}
%The next example explains the 
%To 
for an explanation of  
the argument.
This kind of formulation should be distinguished 
from the so-called %randomized
state-dependent distortion risk measure (SDRM) \cite{wang2023preference} in act space where the consequence of a random act at a particular state is a random variable and the DM's risk preference at different states may differ. In that set-up, the DM's 
risk preference at each state can be described by a deterministic distortion risk measure
which is equivalent to an SRM, and 
SDRM of an act is the statistical average of the DRM at different states.
% }

% We begin by considering some specific parametric SRMs where the DM's risk preferences may be described by one of them.

%Next we provide some specific parametric SRMs where the DM's risk preferences may be described by one of them.
The next example gives rise to a specific parametric 
SRM which is a combination of mathematical expectation and CVaR and explain how it may be randomized.
%and averaged after the randomization.

\begin{example}
\label{exm-cvar-mean-cmb}
(i) Let $\alpha\in(0, 1)$ and $\sigma(z, \alpha) = \frac{1}{1-\alpha} \ind_{\alpha, 1}(z)$, where $\ind_{[a,b]}$ denotes the indicator function over the interval $[a,b]$. 
Then $\rho_s(\xi)=\CVaR_\alpha(\xi)$ with $s=\alpha$ is the most widely used spectral risk measure known as \textit{Expected Shortfall} or Average/Conditional Value at Risk (AVaR/CVaR).

{\color{black}
(ii) Let $s= (\lambda, \alpha)$ with $\alpha$ and $\lambda$ being 
parameters which take values over $[0,1]$, and let  
\begin{equation}
\label{eq:sigmaCombination}
    \sigma(z, s):=\lambda\ind_{[0,1]}(z) + \frac{1-\lambda}{1-\alpha}\ind_{[\alpha,1]}(z).
\end{equation}
Then
\begin{equation}
\rho_{s}(X) = \lambda\bbe[X]+(1-\lambda)\CVaR_\alpha(X), 
\label{eq:SRM-comb-Exp-CVar}
\end{equation}
is a spectral risk measure which balances between risk neutral ($\bbe[X]$) and risk averse ($\CVaR_\alpha(X)$). The risk measure is widely used in operation research and management sciences, see e.g.~\cite{PdF13, KoM15}.
Our interest here is the case when a DM's risk preference can be described by $\rho_s$
but there is no deterministic $(\lambda,\alpha)$ such that the risk measure can be used to describe the DM's preference consistently primarily because the DM's preference is random. In that case, we may randomize the parameters $(\lambda,\alpha)$  to obtain a random spectral risk measure $\rho_s(X)$.
% {\em randomized} convex combination of the expected value and {\em randomized} conditional-value-at-risk (CVaR).
Randomization of $\lambda$ means that the investor's risk preference
switches between risk neutral and risk averse under the framework of SRM 
whereas randomization of $\alpha$ 
%indicates 
means that there is a  variation of investor's  risk averse level on the tail losses.
In 
%this 
the case that $\lambda$ and $\alpha$ are deterministic, the $\rho_s$ 
%becomes
reduces to standard spectral risk measure.}

% {\color{red} We may continue to polish the comments, including why randomized ARM is necessary?
% }

\end{example}

\subsection{Computation of ARSRM}

Let $\xi$ be a random variable with finite
%ly distributed with 
distribution, i.e.,
$\mathbb{P}(\xi=\xi_k)=p_k$ where $-\infty=\xi_0 <\cdots< \xi_{K+1}=+\infty$.
Let $\pi_k:=\sum_{i\in[k-1]} p_i$ for $k\in[K+1]$, where $\pi_1=0$ and $\pi_{K+1}=1$.
The quantile function of $\xi$ is a non-decreasing and left-continuous step-like function with $K$ breakpoints, i.e.,
\begin{equation*}
    F_{\xi}^\leftarrow(z)=\xi_k, \quad\inmat{for\;} z\in(\pi_k,\pi_{k+1}], \;k\in[K].
\end{equation*}
Let $\mathcal{S}$ be the support set of $s$. For $s\in \mathcal{S}$ and the random risk spectrum $\sigma(\cdot,s)$, the discrete distribution of $\xi$ allows us to simplify the RSRM:
\begin{equation}
\label{eq-RSRM-discrete}
    \rho_s(\xi)=\int_0^1 F_{\xi}^\leftarrow(z)\sigma(z,s)dz = \sum_{k\in[K]} \xi_k \int_{\pi_k}^{\pi_{k+1}} \sigma(z,s)dz=\sum_{k\in[K]} \xi_k\psi_k(s),
\end{equation}
where $\psi_k(s):=\int_{\pi_k}^{\pi_{k+1}} \sigma(z,s)dz$ and $\sum_{k\in[K]} \psi_k(s)=1$.
Since $\pi_k$ is the cdf of $\xi$, this formulation relies on the order of the outcomes of $\xi$.

In practice, we often consider the step-like risk spectrum which simplifies the calculation.
%and i
It has been shown that the continuous risk spectrum can be effectively approximated by a step-like function;
%with enough breakpoints; 
see e.g. \cite{WaX20}.
Let $0=z_0<\cdots<z_{J+1}=1$ and
$\sigma(z,s):=\sum_{j\in[J]\cup\{0\}} \sigma_j(s)\ind_{[z_j,z_{j+1})}(z)$ be a random step-like risk spectrum.  Then we can write (\ref{eq-RSRM-discrete}) as
\begin{equation*}
    \rho_{\sigma_s}(\xi) = \sum_{k\in[K]} \xi_k \sum_{j\in[J]\cup\{0\}} \left(\int_{\pi_k}^{\pi_{k+1}}\sigma_j(s)\ind_{[z_j,z_{j+1})}(z) d z \right).
\end{equation*}

Let $Q$ be the probability distribution of $s$. Then the ARSRM with step-like random risk spectrum can be expressed by
\begin{align}
    \rho_{Q}(\xi) & = \int_{\mathcal{S}} \rho_{\sigma_s}(\xi) Q(ds) 
    = \int_{\mathcal{S}} \left[\sum_{k\in[K]} \xi_k \sum_{j\in[J]\cup\{0\}} \left( \int_{\pi_k}^{\pi_{k+1}}\sigma_j(s)
    \ind_{[z_j,z_{j+1})}(z) d z\right) \right]Q(ds) \nonumber \\
    & = \sum_{k\in[K]} \xi_k \sum_{j\in[J]\cup\{0\}} \left(\int_{\mathcal{S}} \sigma_j(s)Q(ds)\right) 
    \int_{\pi_k}^{\pi_{k+1}} \ind_{[z_j,z_{j+1})}(z) dz.
    \label{eq-ARSRM-step}
\end{align}
In the case that $s$ is discretely distributed over $\mathcal{S}$ with $\mathbb{P}(s=s_l)=q_l$, for $l\in[|\mathcal{S}|]$ (where $|\mathcal{S}|$ denotes the cardinality of the set $\mathcal{S}$), we have, %following 
by taking a similar strategy to that of \cite{li2022randomization}, that 
% we have 
\begin{align}
    \rho_{Q}(\xi) & = \sum_{k\in[K]} \xi_k \sum_{j\in[J]\cup\{0\}} \left(\sum_{l\in[|\mathcal{S}|]} q_i\sigma_j(s_l)\right) \int_{\pi_k}^{\pi_{k+1}} \ind_{[z_j,z_{j+1})}(z) dz \nonumber \\
    & = \sum_{l\in[|\mathcal{S}|]} q_l \sum_{k\in[K]} \psi_{l,k}\xi_k
    =\sum_{l\in[|\mathcal{S}|]} \sum_{k\in[K]} q_l(\psi_{l,k}-\psi_{l,k-1}) \left(\sum_{j\in[K]{\color{blue}\setminus}[k-1]} \xi_j \right) \nonumber \\
    & = \sum_{k\in[K]} \sum_{l\in[|\mathcal{S}|]}  q_l\beta_{l,k} \CVaR_{\alpha_k}(\xi), 
    \label{eq-ARSRM-sDiscretized} 
\end{align}
where $\beta_{l,k}:=(\psi_{l,k}-\psi_{l,k-1})(K-k+1)$,
$$\psi_{l,k}:=\int_{\pi_k}^{\pi_{k+1}} \sigma(z,s_l)dz=\sum_{j\in[J]\cup\{0\}} \sigma_j(s_l) \int_{\pi_k}^{\pi_{k+1}}\ind_{[z_j,z_{j+1})}(z) dz, \; \psi_{l,0}=0, $$
and
$\alpha_{k}:=\frac{k-1}{K}$ for $l\in[|\mathcal{S}|], k\in[K]$.
It is easy to verify that $\sum_{l\in[|\mathcal{S}|]} \sum_{k\in[K]} q_l \beta_{l,k}=1$, which implies that $\rho_Q(\xi)$ is a convex combination of $\CVaR_{\alpha_k}(\xi)$, for $k\in[K]$.

\section{Multistage Average Randomized Spectral Risk %optimization
Minimization}

With preparations in the previous section, we are now ready to discuss multistage
decision-making with 
stage-dependent ARSRM. 

\subsection{The MARSRM Model}

We consider a discrete time-horizon problem  
%study a class of multistage stochastic optimization problems 
with a finite number of stages $T>1$.
Let $(\Omega,\mathcal{F},\mathbb{P})$ be a probability space, and $\bm{\xi}:=\{\xi_t\}_{t=1}^T$ be a stochastic process on $(\Omega,\mathcal{F},\mathbb{P})$, where $\xi_t:\Omega\to\Xi_t\subset \R^{d_t}$ represents the random parameters in stage $t$ for $t\in[T]$.
% We assume that $\xi_t, t\in[T]$ has finite possible realizations, where $\xi_1$ is a degenerate random vector.
Let $\bm{\xi}_{[t]}:=(\xi_1,\ldots,\xi_t)$ denote the history of the stochastic process up to time $t$, and $\mathcal{F}_t:=\sigma(\bm{\xi}_{[t]})$ be the filtration associated with $\bm{\xi}_{[t]}$, where $\sigma(\bm{\xi}_{[t]})$ is the $\sigma$-algebra generated by $\bm{\xi}_{[t]}$ and $\{\emptyset, \Omega\}=\mathcal{F}_1\subset\mathcal{F}_2 \subset\cdots \subset\mathcal{F}_T=\mathcal{F}$.
As the 
%information
uncertainty 
%on 
of the stochastic process $\bm{\xi}$ is gradually 
%realized
observed, decisions can be made based on the available information.
%so far, 
Overall, these decisions %define 
form a decision rule or a policy $x=[x_1,\ldots, x_T]$.
Similar to $\xi_{[t]}$, we %use 
write 
$x_{[t]}$ 
%to denote a 
for the policy up to time $t$, i.e., $x_{[t]}:=(x_0,\ldots,x_t)$. 
We begin with a formal definition of 
multistage ARSRM.

\begin{definition}[Multistage ARSRM]
\label{def:MASRM}
Let $\mathcal{L}_1(\Omega, \mathcal F,\mathbb{P})$ denote
the set of random variables $Z:
(\Omega, \mathcal F,\mathbb{P})\to \mathbb{R}$
with finite first order moment, i.e.,
$\int_{\Omega} |Z(\omega)| d\mathbb{P}(\omega)<\infty$.
For $t=2,\cdots,T$, let  $\rho_{Q_t|\mathcal{F}_{t-1}}:
\mathcal{L}_1(\Omega, \mathcal F,\mathbb{P}) \to \R
$ 
be ARSRMs such that 
$$
\rho_{Q_t|\mathcal{F}_{t-1}}(Z_t) = \int_{\mathcal{S}_t} \rho_{\sigma_s}(Z_t|\mathcal{F}_{t-1}) Q_t (ds), \quad \text{for}\; Z_t\in 
\mathcal{L}^1(\Omega, \mathcal F,\mathbb{P}),
$$
where $Q_t$ is the probability distribution of $s$ at stage $t$ with support set ${\cal S}_t$. 
\end{definition}
With the definition of  the multistage ARSRM,  
we are able to introduce our
multistage decision making problem with 
MSRSRM:
\begin{equation}
\label{eq-MARSRM}
\begin{split}
    & \min_{x_1\in \mathcal{X}_1 (x_0, \xi_1)} \;\; c_1x_1  + 
    % \sup_{Q_1\in\mathcal{Q}_1(\xi_1)} 
    \rho_{Q_2|\mathcal{F}_1}  \Bigg[ \inf_{x_2\in\mathcal{X}_2(x_1,\xi_2)} c_2x_2 +\cdots \\
     & \qquad + \rho_{Q_{T-1}|\mathcal{F}_{T-2}} 
    \bigg[ \inf_{x_{T-1}\in\mathcal{X}_{T-1}(x_{T-2},\xi_{T-1})} c_{T-1}x_{T-1} 
      + \rho_{Q_{T}|\mathcal{F}_{T-1}}
    \Big[ \inf_{x_T\in\mathcal{X}_T(x_{T-1},\xi_T)} c_Tx_T \Big]\bigg]\Bigg],
\end{split}
\end{equation}
where 
$\mathcal{X}_1 (x_{0}, \xi_{1} ):=\left\{x_1\in \R^{n_1}_+: A_1x_1= b_1\right\}$,
% {\color{red}Why do we restrict $x$ to taking non-negative values?}
and
\begin{equation}
\label{eq-mathcalX_t}
    \mathcal{X}_t(x_{t-1},\xi_{t}):=\{x_t\in\R^{n_t}_+: A_t x_t= b_t-E_t x_{t-1}\}, \; \text{for} \; t\in[T]\setminus\{1\}, 
\end{equation}
where
$\xi_t :=(c_t, b_t, A_t,E_t)$ is
a vector of random variables
with 
%the degenerated (deterministic) first stage parameters 
$\xi_1 :=(c_1, b_1, A_1,E_1)$ 
being deterministic 
and $x_0$ is a given input,
%a dummy variable, 
$A_t, E_t\in\R^{m_t}\times\R^{n_t}$ are random matrices, for $t\in[T]$.
Let
$$
f_t(x_t,\xi_t,Q_{t+1}) := c_t x_t+\rho_{Q_{t+1}|\mathcal{F}_{t}}
\left( V_{t+1}(x_{t}, \xi_{t+1})
%\min_{x_{t+1}\in \mathcal{X}_t(x_{t-1},\xi_{t})}f_{t+1}(x_{t+1},\xi_{t+1},Q_{t+2}) 
\right). 
$$
Then problem (\ref{eq-MARSRM}) can be written recursively as
\bgeqn 
\label{eq:recourseobjective}
V_t(x_{t-1}, \xi_t) :=
\min_{x_t\in \mathcal{X}_t(x_{t-1},\xi_{t})} f_t(x_t,\xi_t,Q_{t+1}) 
%:= c_t x_t+\rho_{Q_{t+1}|\mathcal{F}_{t}}
%\left( V_{t+1}(x_{t}, \xi_{t+1})
%\min_{x_{t+1}\in \mathcal{X}_t(x_{t-1},\xi_{t})}f_{t+1}(x_{t+1},\xi_{t+1},Q_{t+2}) 
%\right)
% \inf_{x_{t+1}\in\mathcal{X}_{t+1}(x_t,\xi_{t+1})} c_{t+1} x_{t+1}
% + \rho_{Q_{t+2}|\mathcal{F}_{t+1}}
% \left(f_{t+1}(x_{t}, \xi_{t+1},Q_{t+2} \right)
% \right)
\edeqn
for $t=T-1, T-2,\ldots, 1$, where
\bgeq 
V_T(x_T,\xi_T) = \inf_{x_T\in\mathcal{X}_T(x_{T-1},\xi_T)} f_T(x_T,\xi_T,Q_{T+1})
\quad \inmat{and} \quad
f_T(x_T,\xi_T,Q_{T+1}) = c_Tx_T. 
\edeq
In this model, the DM chooses
an optimal decision 
$x_t$ at stage $t$, after observing 
realization of $\xi_t$, 
from feasible set 
$\mathcal{X}_t(x_{t-1},\xi_{t})$, i.e., 
which minimizes
$f_t(x_t,\xi_t,Q_{t+1})$.
% \bgeqn 
% \label{eq:recourseobjective}
% V_t(x_{t-1}, \xi_t) :=
% \min_{x_t\in \mathcal{X}_t(x_{t-1},\xi_{t})} f_t(x_t,\xi_t,Q_{t+1}) := c_t x_t+\rho_{Q_{t+1}|\mathcal{F}_{t}}
% \left( V_{t+1}(x_{t}, \xi_{t+1})
% %\min_{x_{t+1}\in \mathcal{X}_t(x_{t-1},\xi_{t})}f_{t+1}(x_{t+1},\xi_{t+1},Q_{t+2}) 
% \right)
% % \inf_{x_{t+1}\in\mathcal{X}_{t+1}(x_t,\xi_{t+1})} c_{t+1} x_{t+1}
% % + \rho_{Q_{t+2}|\mathcal{F}_{t+1}}
% % \left(f_{t+1}(x_{t}, \xi_{t+1},Q_{t+2} \right)
% % \right)
% \edeqn
% for $t=T-1, T-2,\ldots, 1$, where
% \bgeq 
% f_T(x_T,\xi_T,Q_{T+1}) = c_Tx_T \quad \inmat{and} \quad
% V_T(x_T,\xi_T) = \inf_{x_T\in\mathcal{X}_T(x_{T-1},\xi_T)} f_T(x_T,\xi_T,Q_{T+1}).
% \edeq
In other words, at each stage, we solve a recourse problem. The term
$\rho_{Q_{t+1}|\mathcal{F}_{t}}
\left( V_{t+1}(x_{t}, \xi_{t+1})\right)$
% The second term in (\ref{eq:recourseobjective}) 
%$$
% \rho_{Q_{t+1}|\mathcal{F}_{t}}
% \left(\inf_{x_{t+1}\in\mathcal{X}_{t+1}(x_t,\xi_{t+1})} c_{t+1} x_{t+1} +
%\right)
%$$
quantifies the DM's average randomized SRM
of $
f_{t+1}(x_{t+1},\xi_{t+1},Q_{t+2})$.
This is because from modeller's perspective,
the DM's valuation of the risk on the remaining losses is random.
In a particular case when
%If 
$\sigma$ is defined as in Example \ref{exm-cvar-mean-cmb} (ii) with $\lambda=1$, 
%then 
$\rho_{\sigma(\cdot,1,\alpha)}(X)=\bbe[X]$ and MARSRM
collapses to a classical risk-neutral multistage problem 
% {\color{blue} 
see e.g., \cite{pereira1991multi, ruszczynski1997decomposition}:
\begin{equation}
\label{eq-ms-riskneutral}
\begin{split}
    & \min_{x_1\in \mathcal{X}_1 (x_0, \xi_1)} \;\; c_1x_1  + 
    % \sup_{Q_1\in\mathcal{Q}_1(\xi_1)} 
    \bbe_{P_2|\mathcal{F}_1}  \Bigg[ \inf_{x_2\in\mathcal{X}_2(x_1,\xi_2)} c_2x_2 +\cdots \\
     & \qquad + \bbe_{P_{T-1}|\mathcal{F}_{T-2}} 
    \bigg[ \inf_{x_{T-1}\in\mathcal{X}_{T-1}(x_{T-2},\xi_{T-1})} c_{T-1}x_{T-1} 
      + \bbe_{P_{T}|\mathcal{F}_{T-1}}
    \Big[ \inf_{x_T\in\mathcal{X}_T(x_{T-1},\xi_T)} c_Tx_T \Big]\bigg]\Bigg].
\end{split}
\end{equation}
% }
In this case, the DM's valuation 
of future loss is a deterministic risk measure 
($\bbe[\cdot]$).
On the other hand, 
%However, if we set 
% {\color{blue}
when
$Q_t$ is the Dirac probability measure at a particular pair $(\lambda, \alpha)$ with $\lambda, \alpha \in (0,1)$,
% = \delta(\alpha)$ 
$\rho_{\sigma(\cdot, \lambda, \alpha)}(X) = \lambda\bbe[X]+ (1-\lambda) \CVaR_{\alpha}[X]$ and consequently 
MARSRM reduces to 
\begin{equation}
\label{eq-ms-riskaverse}
\begin{split}
    & \min_{x_1\in \mathcal{X}_1 (x_0, \xi_1)} \;\; c_1x_1  + 
    % \sup_{Q_1\in\mathcal{Q}_1(\xi_1)} 
    \rho_{2|\mathcal{F}_1 }\Bigg[ \inf_{x_2\in\mathcal{X}_2(x_1,\xi_2)} c_2x_2 +\cdots \\
     & + \rho_{T-1|\mathcal{F}_{T-2}} 
    \bigg[ \inf_{x_{T-1}\in\mathcal{X}_{T-1}(x_{T-2},\xi_{T-1})} c_{T-1}x_{T-1} 
      + \rho_{T|\mathcal{F}_{T-1}}
    \Big[ \inf_{x_T\in\mathcal{X}_T(x_{T-1},\xi_T)} c_Tx_T \Big]  \bigg] \Bigg].
\end{split}
\end{equation}
In that case,
the DM's risk preference is evaluated by a deterministic combination of expectation and CVaR, 
see e.g., \cite{PdF13, KoM15, ZRG16}.
% }
% which is a risk-averse multistage problem with randomized Conditional-Value-at-Risk. We will come back to the specific cases in Section 5.
%{\color{red}HX: please add some references to the above two models.}

From the computational perspective, problem (\ref{eq-MARSRM}) can be solved 
in a dynamic recursive manner:
%by setting
\begin{equation}
\label{eq-costToGoFunction}
    V_t(x_{t-1},\xi_{t}) = \min_{x_t\in\mathcal{X}_t(x_{t-1},\xi_t)} \; c_tx_t+
    % \max_{Q_t\in\mathcal{Q}_t(\xi_t)} 
    \rho_{Q_{t+1}}(V_{t+1}(x_{t},\xi_{t+1})),
\end{equation}
for $t=T,T-1,\ldots,1$, where
%By setting 
$\rho_{Q_{T+1}}(V_{T+1}(x_{T},\xi_{T+1}))\equiv 0$
and $V_1(x_0,\xi_1)$ is the optimal value of problem (\ref{eq-MARSRM}).
By convention, we call 
$V_t(x_{t-1},\xi_{t})$
the {\em cost-to-go function}
at stage $t$.
To %evaluate 
%the cost-to-go function,
calculate $V_t(x_{t-1},\xi_{t})$,
we consider the case that 
both $\xi_t$ and $s_t$ have finite number of scenarios. This can be viewed as sample average approximation if they are 
continuously distributed originally. 
In the latter case, one may be interested 
in the error bounds on the approximations.
We refer readers to \cite{li2022randomization}
for the case when $T=1$. The next assumption 
states the set-up of the discrete distributions.

% where in the last stage we assume for simplicity that $\rho_{Q_{T+1}}(V_{T+1}(x_{T},\xi_{T+1}))=0$.
% We define $\xi_t=(A_t,E_t,b_t,c_t)$ and the feasibility sets as: 
% \begin{equation*}
%     \mathcal{X}_1 (x_{0}, \xi_{1} ):=\{x_1\geq 0: A_1x_1= b_1\}
% \end{equation*}
% and 
% \begin{equation*}
%     \mathcal{X}_t(x_{t-1},\xi_{t}):=\{x_t\geq 0: A_tx_t= b_t-E_t x_{t-1}\}, \; t\in[T]/\{1\}.
% \end{equation*}

% Define $\xi_t=(A_t,E_t,b_t,c_t)$. 
% To evaluate the cost-to-go function, 
%We make the following assumptions:
\begin{assumption}
\label{assumption-multistageproblem}
    \begin{enumerate}[(a)]
        \item For $t\in[T]$, $\xi_t$ has finite possible realizations, 
        % and the support set $\Xi_t$ of $\xi_t$ is bounded, where $\xi_1$ is a degenerate random vector.
        %Moreover, t
         and the random process $\{\xi_t\}_{t=1}^T$ is stagewise independent.
        \item The random preference parameter $s_t$ is discretely distributed over finite support set $\mathcal{S}_t$ with $\mathbb{P}(s_t=s_{t, l})=q_{t,l}$, for $l\in[|\mathcal{S}_t|]$, $t=2,\ldots,T$.
        %at stage $t$.
        
        \item  For $t\in[T]$, problem (\ref{eq-MARSRM}) has  relatively complete recourse, i.e.,
        $\mathcal{X}_t(x_{t-1},\xi_{t})\neq \emptyset$,
        %for $t\in[T]\setminus\{1\}$ 
        %       the first-stage problem is feasible and stage t problem is feasible given any 
        for almost every 
        realization of $\xi_t$ 
                and any feasible solution $x_{t-1}$ 
        from previous stage.
        % for $t\in[T]\setminus\{1\}$. 
         Moreover, the set $\mathcal{X}_t(x_{t-1},\xi_{t})$ is compact.
    \end{enumerate}
\end{assumption}
%{\color{red} 

Assumption \ref{assumption-multistageproblem} (a) and (c) are standard, see e.g., \cite{shapiro2013risk, PdF13, ZRG16,duque2020distributionally}.
Assumption \ref{assumption-multistageproblem} (b) is made 
for tractable reformulation of problem (\ref{eq-MARSRM}). In practice, the DM's preference may change continuously, in which case we may view (\ref{eq-MARSRM}) as an approximation of the real. Under Assumption \ref{assumption-multistageproblem},
% By combining (\ref{eq-ARSRM-sDiscretized}),  (\ref{eq-costToGoFunction}) and 
%, for given $\tilde{\xi}_t$, 
we can use (\ref{eq-ARSRM-sDiscretized})
% and  (\ref{eq-costToGoFunction}) 
to 
%reformulate  problem
%simplify 
derive a discretized 
recursive formulation for 
%from
(\ref{eq-costToGoFunction}):
% \begin{equation}
% \label{eq-costToGoFunction-ARSRM}
%     \begin{split}
%         & V_t (x_{t-1}, \tilde{\xi}_t) \\
%         & :=  \min_{x_t\in\mathcal{X}_t(x_{t-1}, \xi_t)} \; c_t x_t 
%         + \sum_{k\in[K_{t+1}]} \sum_{l\in[I_{t+1}]}  q_{t+1,l}\beta_{t+1,l,k} \CVaR_{\alpha_{t+1,k}} \left(V_{t+1}(x_{t},\xi_{t+1})| \xi_t=\tilde{\xi}_t \right), 
%     \end{split}
% \end{equation}
\begin{equation}
\label{eq-costToGoFunction-ARSRM-0}
\resizebox{0.925\linewidth}{!}{$
    V_t (x_{t-1}, \tilde{\xi}_t) :=  \displaystyle{\min_{x_t\in\mathcal{X}_t(x_{t-1}, \xi_t)} } \; c_t x_t 
    + \sum_{k\in[K_{t+1}]} \sum_{l\in[I_{t+1}]}  q_{t+1,l}\beta_{t+1,l,k} \CVaR_{\alpha_{t+1,k}} \left(V_{t+1}(x_{t},\xi_{t+1})| \xi_t=\tilde{\xi}_t \right), 
    $}
\end{equation}
where $K_{t+1}$ is the number of scenarios of $\xi_{t+1}$,
$I_{t+1}$ is the number of scenarios of random preference parameter at stage $t+1$,
$\beta_{t,l,k}:=(\psi_{t,l,k}-\psi_{t,l,k-1})(K_t-k+1)$, 
\bgeqn 
\psi_{t,l,k}:=\int_{\pi_{t,k}}^{\pi_{t,k+1}} \sigma_t(z,s_{t,l})dz=\sum_{j\in[J_t]\cup\{0\}}\sigma_{t,j}(s_{t,l}) \int_{\pi_{t,k}}^{\pi_{t,k+1}}\1_{[z_{t,j},z_{t,j+1})}(z) dz, \; \psi_{t,l,0}=0, 
\edeqn 
%and
$\alpha_{t, k}:=\frac{k-1}{K_t}$ for $l\in[|\mathcal{S}_t|], k\in[K_t]$, and 
$\tilde{\xi}_t$
represents a particular scenario of $\xi_t$.

\subsection{SDDP for Solving MARSRM}

%In the rest of this section, 
We now move on to discuss how to use the well-known SDDP method 
to solve the MARSRM. As in \cite{dias2010stochastic, PdF13, georghiou2019robust, silva2021data, lan2024numerical}, we develop 
a process which is a combination of 
lower/upper approximations for constructing  
lower/upper bounds of the optimal values at different stages.
To this end, we define
\begin{equation}
\label{eq-secondTerm}
    \mathcal{V}_{t+1}(x_t, \tilde{\xi}_t):=
    \sum_{k\in[K_{t+1}]} \sum_{l\in[I_{t+1}]} q_{t+1,l}\beta_{t+1,l,k} \CVaR_{\alpha_{t+1,k}} \left(V_{t+1}(x_{t},\xi_{t+1})| \xi_t=\tilde{\xi}_t \right)
\end{equation}
and write (\ref{eq-costToGoFunction-ARSRM-0})
\begin{equation}
\label{eq-costToGoFunction-ARSRM}
%\resizebox{0.925\linewidth}{!}{$
    V_t (x_{t-1}, \tilde{\xi}_t) =  \displaystyle{\min_{x_t\in\mathcal{X}_t(x_{t-1}, \xi_t)} } \; c_t x_t 
    + \mathcal{V}_{t+1}(x_t, \tilde{\xi}_t)
    %\sum_{k\in[K_{t+1}]} \sum_{l\in[I_{t+1}]}  q_{t+1,l}\beta_{t+1,l,k} \CVaR_{\alpha_{t+1,k}} \left(V_{t+1}(x_{t},\xi_{t+1})| \xi_t=\tilde{\xi}_t \right). 
%    $}
\end{equation}
Observe that $\mathcal{V}_{t+1}(x_t, \tilde{\xi}_t)$ is a convex function of $x_t$ when $V_{t+1}(x_t,\xi_{t+1})$ is convex 
in $x_t$ for every realization of $\xi_{t+1}$. 
% because $\mathcal{V}_{t+1}(x_t, \tilde{\xi}_t)$ is the convex combination of $\CVaR_{\alpha_{t+1,k}}$ 
% %and $\CVaR_{\alpha_{t+1,k}}$ 
% is convex for $k\in[K_{t+1}]$.
%Furthermore
%Moreover, t
% This implies $V_{t}(x_{t-1}, \tilde{\xi}_{t})$ is convex in $x_{t-1}$ for every fixed $\tilde{\xi}_{t}$.
% %when $V_{t+1}(x_t,\xi_{t+1})$ is convex in $x_t$ for every realization of $\xi_{t+1}$ and $t\in[T]\setminus\{1\}$.
By induction backward from stage $T$, we 
%can have 
conclude that $V_{t}(x_{t-1}, \tilde{\xi}_{t})$ is convex in $x_{t-1}$ for every realization of $\tilde{\xi}_{t}$.
%Thus we can use 
The convexity allows us 
to 
construct 
cutting planes 
which form 
%to construct 
%the
a lower approximation of $V_{t}(x_{t-1}, \xi_{t})$. Here we give a sketch of the process.
% and upper approximations for constructing  
% %lower
% upper bounds of the optimal values.

% {\color{red}
% We now proceed to discuss schemes to obtain optimal policies for MARSRM (\ref{eq-MARSRM}).
% One method for building lower approximation which gives lower bounds on the optimal value for the problem is the algorithm based on \cite{kelley1960cutting}.
% % , and has been widely used in e.g. \cite{}.
% The method for constructing the upper approximation which gives upper bound on the optimal value for the problem is the algorithm developed by \cite{dias2010stochastic, PdF13, georghiou2019robust, silva2021data}.
% }
% % and used in \cite{}.
% These algorithms are simplest to describe when the risk measure is expectation, so we first consider this case and extend them to the ARSRM case.

% From the above discussions, we have the ARSRM at stage $t$:
% % At each stage $t-1$, we have the following subproblem:
% \begin{align}
% \label{eq-rho-t}
%     \rho_{Q_t}(V_t(x_{t-1},\xi_{t})) & = \sum_{i\in[I_t]} \sum_{k\in[K_t]} q_{t,i}\beta_{t,i,k} \CVaR_{\alpha_k}(V_t(x_{t-1},\xi_{t})|\mathcal{F}_{t-1}), 
%     \\
%     & =\min_{{\bm \eta}_{t-1}} \;\; \sum_{i\in[I_t]} \sum_{k\in[K_t]} q_{t,i}\beta_{t,i,k} 
%     \left\{\eta_{t-1,k}+\frac{1}{1-\alpha_{t,k}}\bbe[(V_t(x_{t-1},\xi_{t})-\eta_{t-1,k})_+|\mathcal{F}_{t-1}]\right\} \label{eq-rho-eta}.
% \end{align}
% Under Assumption~\ref{A-RSRM} (a), $\rho_{Q_{t+1}|\mathcal{F}_t}$ coincides with its unconditional counterpart $\rho_{Q_{t+1}}$.

\subsubsection{Lower Bound of \texorpdfstring{$V_t (x_{t-1}, \tilde{\xi}_t)$}{} } 
%Lower approximation of MARSRM}
 
We 
%start from 
begin by constructing 
%the
lower approximations of 
%the cost-to-go function
$\mathcal{V}_{t+1}(x_t, \tilde{\xi}_t)$ defined in (\ref{eq-costToGoFunction-ARSRM}) 
based on 
%cutting planes of
lower linear envelope of the 
convex function to obtain 
%its
a sequence of lower bound of $ V_t (x_{t-1}, \tilde{\xi}_t)$.
% {\color{red}HX: the terminologies of lower approximation and lower approximation are confusing.}
%{\color{blue}
In the case when the risk measure reduces to the mathematical expectation, this kind of scheme is widely used in the literature of Bender's decomposition method and 
%stochastic dual decomposition method 
SDDP for solving mutistage stochastic programming problems, see e.g.~\cite{PdF13, ZRG16}.
%The scheme for building the lower approximation is simplest to describe when the risk measure is expectation.
% , where the approximation is constructed using cutting planes.
% The details are described 
% %at length, e.g.
% in \cite{PdF13, ZRG16},
%so
Here we use the same
%directly 
%extend 
%them
%the approximation scheme 
approach to approximate 
%the
MARSRM problem.
%The key point 
%here
A key step of the approach 
is to construct, at each iteration,
a 
cutting plane 
%approximation 
of 
$\mathcal{V}_{t+1}(x_t, \tilde{\xi}_t)$
%the second term of 
in (\ref{eq-costToGoFunction-ARSRM}).
%, i.e., 
% \begin{equation}
% \label{eq-secondTerm}
%     \mathcal{V}_{t+1}(x_t, \tilde{\xi}_t):=
%     \sum_{k\in[K_{t+1}]} \sum_{l\in[I_{t+1}]} q_{t+1,l}\beta_{t+1,l,k} \CVaR_{\alpha_{t+1,k}} \left(V_{t+1}(x_{t},\xi_{t+1})| \xi_t=\tilde{\xi}_t \right).
% \end{equation}
% Observe that $\mathcal{V}_{t+1}(x_t, \tilde{\xi}_t)$ is a convex function of $x_t$ when $V_{t+1}(x_t,\xi_{t+1})$ is convex 
% in $x_t$ for every realization of $\xi_{t+1}$. 
% % because $\mathcal{V}_{t+1}(x_t, \tilde{\xi}_t)$ is the convex combination of $\CVaR_{\alpha_{t+1,k}}$ 
% % %and $\CVaR_{\alpha_{t+1,k}}$ 
% % is convex for $k\in[K_{t+1}]$.
% %Furthermore
% %Moreover, t
% % This implies $V_{t}(x_{t-1}, \tilde{\xi}_{t})$ is convex in $x_{t-1}$ for every fixed $\tilde{\xi}_{t}$.
% % %when $V_{t+1}(x_t,\xi_{t+1})$ is convex in $x_t$ for every realization of $\xi_{t+1}$ and $t\in[T]\setminus\{1\}$.
% By induction backward from stage $T$, we 
% %can have 
% conclude that $V_{t}(x_{t-1}, \tilde{\xi}_{t})$ is convex in $x_{t-1}$ for every realization of $\tilde{\xi}_{t}$.
% %Thus we can use 
% The convexity allows us 
% to 
% construct 
% cutting planes 
% which form 
% %to construct 
% %the
% a lower approximation of $V_{t}(x_{t-1}, \xi_{t})$. Here we give a sketch of the process.
% %at each stage.

Let $i$ denote the number of iterations starting with $i=0$
and $\underline{V}_{t, i} (\cdot, \tilde{\xi}_t)$ the lower approximation 
of $V_{t, i} (\cdot, \tilde{\xi}_t)$ in $x_{t-1}$ for each observed $\tilde{\xi}_t$. 
We consider the $i$-th iteration.
At 
%each
stage $t$ ($t<T$), we are given
$\underline{V}_{t+1, i} (\cdot, \xi_{t+1})$ and  we 
solve 
% we have 
the following subproblem:
\begin{subequations}
\label{eq-sub-re2}
\begin{align}
    \underline{V}_{t, i} (x_{t-1}, \tilde{\xi}_t) = \min_{x_t\geq0,\theta_{t+1}} \;\;& \tilde{c}_t x_t +\theta_{t+1} \label{eq-sub-obj-re2} \\
    \st \;\;\; & \tilde{A_t} x_t= \tilde{b}_t-\tilde{E}_t x_{t-1}, \label{eq-sub-cons-o-re2} \\
    & g_{t+1,n}+G_{t+1,n} x_t \leq \theta_{t+1}, \; \text{for}\; n\in [\mathcal{N}_{t+1, i}],
    \label{eq-sub-cons-theta-re2}
\end{align}
\end{subequations}
%{\color{blue} 
where
$[\mathcal{N}_{t+1, i}]$ denotes
the index set of linear pieces of the 
piecewise linear convex function $\underline{V}_{t+1, i} (x_t, \xi_{t+1})$
and 
%This formulation is based on the fact that 
%$\underline{V}_{t+1, i} (x_t, \xi_{t+1})$
%is a piecewise linear convex function 
%with linear pieces 
$g_{t+1,n}+G_{t+1,n} x_t=\theta_{t+1}$ 
denotes the $n$-th linear piece.
%for $n\in [\mathcal{N}_{t+1, i}]$.
Each linear piece represents a cut constructed 
from 
%each 
%{\color{blue}
one of the previous iterations,
which means 
$[\mathcal{N}_{t+1, i}] =\{1,\cdots, i*N\}$, 
where $N$ denotes the number of sampled 
paths of $\xi_{[T]}$. 
Inequality (\ref{eq-sub-cons-theta-re2}) is no more than epigraphical formulation of 
$\underline{V}_{t+1, i} (x_t, \xi_{t+1})$
with  intermediate variable $\theta_{t+1}$.
% is an intermediate variable for epigraphical formulation of 
% the piecewise linear function
% , 
% in the objective function (\ref{eq-sub-obj-re2}),
%cut constraints (\ref{eq-sub-cons-theta-re2}), forms the current approximation, denoted by $\underline{V}_{t+1, i} (x_t, \cdot)$,
% using $\mathcal{N}_{t+1, i}$ cuts,
% where $\mathcal{N}_{t+1, i}$  
% denotes the index set of previous cuts at stage $t+1$,
% %which depends on the number $i$ of iteration, 
% $G_{t+1, n}$ and $g_{t+1, n}$ denote the gradients and intercepts of the optimality cuts. 
% {\color{red} I think we need to define iteration index $i$
% and lower approximation of the future risk after $t$ $\underline{V}_{t+1, i} (x_t, \cdot)$ first.
% Since we assume the relative complete recourse of the problem, we do not need any feasibility cuts. ?? I don't understand.
% }
At the last stage $T$,
$\rho_{Q_{T+1}}(V_{T+1}(x_{T},\xi_{T+1}))=0$ and consequently 
% there is no recourse function to be
% %approximated;
% taken, therefore, 
% %at this stage, 
$\theta_{T+1}$ and the cut
constraints (\ref{eq-sub-cons-theta-re2}) are absent.
% For each realization $\xi_{t,j}$ of $\xi_t$, let $\tau_{t, j}$ be the dual variables associate with constraint (\ref{eq-sub-cons-o-re2}.

We need to update $\underline{V}_{t, i} (x_{t-1}, \cdot)$ by 
constructing a new cut. To this end,
%of $\underline{V}_{t, i} (x_{t-1}, \cdot)$, 
we  reformulate 
(\ref{eq-secondTerm}) by exploiting 
 a dual formulation of CVaR.
% rewrite (\ref{eq-secondTerm}) by using the expectation with a ``changed'' probability measure to replace the $\CVaR$ calculation, and then we can
% % $\mathcal{V}_{t+1}(x_t, \tilde{\xi}_t)$ 
% obtain its
% % its subgradient 
% % and 
% cutting plane approximation.
% According to \cite{ADEH99,SDR09}, any coherent risk measure has a dual representation. 
Specifically, let $Z$ be a random variable with finite realizations $Z_1\leq \cdots \leq Z_K$ and corresponding probabilities $p_1, \ldots, p_K$. 
Then
% \begin{equation}
% \label{eq-CRM-dual}
%     \rho(\xi) = \sup_{\lambda\in\Lambda} \sum_{k\in[K]} p_k\lambda_k \xi_k, 
% \end{equation}
\begin{equation}
\label{eq-CRM-dual}
    \CVaR_{\alpha}(Z) = \sup_{\lambda\in\Lambda_{\alpha}} \sum_{k\in[K]} p_k\lambda_k Z_k, 
\end{equation}
where 
% $\Lambda$ is a convex subset of 
% % \begin{equation*}
% $
% \mathscr{B} = \left\{\lambda\in\R^K:\sum_{k\in[K]} p_k\lambda_k=1, \lambda\geq0 \right\}.
% $
% % \end{equation*}
% In the special case where the risk measure is $\CVaR_{\alpha}(\xi)$, we have 
\begin{equation*}
    \Lambda_{\alpha} = \left\{\lambda\in \R^K_+ : \sum_{k\in[K]} p_k\lambda_k=1, \lambda_k\leq \frac{1}{1-\alpha}, \;k\in[K] \right\},
\end{equation*}
see e.g.~\cite{ADEH99,SDR09}.
%which 
Since $ \Lambda_{\alpha}$ is a compact set,
%and 
the supremum in (\ref{eq-CRM-dual}) is attainable.
%can be attained.
% Without loss of generality, suppose $X(\omega_1)\leq\cdots\leq X(\omega_K)$.
% Since we assume that $\xi_1 \leq \cdots \leq \xi_k$.
Define index $\hat{k}$ such that $\sum_{k\in[\hat{k}-1]} p_k\leq\alpha<\sum_{k\in[\hat{k}]} p_k$,
where $\hat{k}=K$ if $1-p_K\leq\alpha$. 
Then $\CVaR_{\alpha}(Z)=\sum_{k\in[K]} p_k \hat{\lambda}_k Z_k$,
where 
\begin{equation}
\label{eq-obtain-lambda}
    \hat{\lambda}_k=
    \begin{cases}
        0, & \text{for}\; k<\hat{k}, \\
        \left(1-\frac{1}{1-\alpha}\sum_{i\in[K]\setminus[\hat{k}]} p_i\right)/p_{\hat{k}}, & \text{for}\; k=\hat{k}, \\
        \frac{1}{1-\alpha}, & \text{for}\; k>\hat{k}.
    \end{cases}
\end{equation}
%Then t
Based on (\ref{eq-CRM-dual}) and (\ref{eq-obtain-lambda}), 
we can rewrite the lower approximation of $V_t (x_{t-1}, \tilde{\xi}_t )$ as 
%the change-of-measure-based subproblem 
%can be written 
%as 
\begin{equation}
\label{eq:change-of-measure-based subproblem}
\resizebox{0.9\linewidth}{!}{$
    \underline{V}_t (x_{t-1}, \tilde{\xi}_t ) :=  \displaystyle \min_{x_t\in\mathcal{X}_t(x_{t-1},\tilde{\xi}_t)} \left\{ c_t x_t 
    + \sum_{k\in[K_{t+1}]} \sum_{l\in[\mathcal{S}_{t+1}]} q_{t+1,l}\beta_{t+1,l,k} 
    % \hat{\bbe}_{t+1|\tilde{\xi}_t} [V_{t+1}(x_{t},\xi_{t+1})],
    \sum_{j\in[K_{t+1}]} p_{t+1,j} \hat{\lambda}_{t+1, k, j} \underline{V}_{t+1}(x_{t},\xi_{t+1, j})\right\}.
    $}
\end{equation}
% {\color{red} Is it from (\ref{eq-costToGoFunction-ARSRM})?
% The set of samples at stage $t+1$
% ${\cal S}_{t+1}$ is not defined.
% }
% where $\hat{\bbe}_{t+1|\tilde{\xi}_t}$ denotes the expectation of the ``changed'' probability measure $\hat{P}$ with $\hat{P}(Z_t=Z_{t,k}|\xi_t = \tilde{\xi}_t) = p_{t,k} \hat{\lambda}_{t,k}$. 
%{\color{red}
It is important to note that the objective function in the rhs of minimization problem  (\ref{eq:change-of-measure-based subproblem}) is independent of $x_{t-1}$.
However, we write down the Lagrange dual of the problem, its objective function depends on $x_{t-1}$. 
Let $\tau_{t, j}$ be the dual variables associated with constraint (\ref{eq-sub-cons-o-re2}) for
% $x_{t-1}= \hat{x}_{t-1}$ and 
$\tilde{\xi}_t = \xi_{t, j}$. 
Then the new cut 
associated with $x_{t-1}$ 
can be constructed 
%is defined 
by
\begin{subequations}
\begin{align}
    & G_{t, i} = -\sum_{l\in[|\mathcal{S}_t|]} \sum_{k\in[K_t]} q_{t,l} \beta_{t,l,k} \sum_{j\in[K_t]} p_{t,j}\hat{\lambda}_{t,k,j} \tau_{t, j} E_{t, j}, \label{eq-cut-re2-a} \\
    & g_{t, i} = \sum_{l\in[|\mathcal{S}_t|]} \sum_{k\in[K_t]} q_{t,l} \beta_{t,l,k} \sum_{j\in[K_t]} p_{t,j}\hat{\lambda}_{t,k,j} \underline{V}_{t, i} (\hat{x}_{t-1}, \xi_{t, j})- G_{t,i} \hat{x}_{t-1}, \label{eq-cut-re2-b}
\end{align}
\end{subequations}
where $\bm{\hat{\lambda}}_{t, k}$ is chosen (as in (\ref{eq-obtain-lambda})) to maximize $\sum_{j\in[K_t]} p_{t,j} \lambda_{t, k, j} \underline{V}_{t, i} (\hat{x}_{t-1},\xi_{t, j})$ over $\Lambda_{\alpha_{t,k}}$ which is defined as:
$$
\Lambda_{\alpha_{t,k}}:= \left\{\bm{\lambda}\in \R^{K_t}_+ : \sum_{j\in[K_t]} p_{t, j} \lambda_{t, k, j}=1,
\lambda_{t, k, j}\leq\frac{1}{1-\alpha_{t,k}},\; \text{for}\;j\in[K_t] \right\}, \; \text{for } k\in[K_t].
$$
% {\color{red}Remind the meaning 
% of this definition}\\
% which depends on the current variable values.
Note that $G_{t,i}$ is a subgradient 
of $\underline{V}_{t}(x_{t-1}, \tilde{\xi}_t)$ 
in $x_{t-1}$ at $x_{t-1}=\hat{x}_{t-1}$
and its formulation
in (\ref{eq-cut-re2-a})
%it is calculated 
is %essentially derived by virtue of 
established by the application of the
 well-known Danskin's theorem \cite{bertsekas2016nonlinear}
 %, i.e.,
%a subgradient of the objective function of 
%applied 
to the Lagrange dual of problem
(\ref{eq:change-of-measure-based subproblem}) 
%with respect to
parameterized by $x_{t-1}$ (evaluated as $x_{t-1}=\hat{x}_{t-1}$).
% is evaluated at the optimal solution of the minimization problem (with respect to $x_t$). 
%This subgradient 
% Since the objective function of problem (\ref{eq:change-of-measure-based subproblem}) is a piecewise linear convex function, we can use \cite[Theorem 23.8]{Roc70} to obtain
% \bgeq
% && \partial_{x_{t-1}} \Bigg(
% \sum_{k\in[K_{t+1}]} \sum_{l\in[\mathcal{S}_{t+1}]} q_{t+1,l}\beta_{t+1,l,k} 
%     % \hat{\bbe}_{t+1|\tilde{\xi}_t} [V_{t+1}(x_{t},\xi_{t+1})],
%     \sum_{j\in[K_{t+1}]} p_{t+1,j} \hat{\lambda}_{t+1, k, j} \underline{V}_{t}(x_{t-1},\xi_{t, j})\Bigg) \Bigg|_{x_{t-1}=\hat{x}_{t-1}} \\
% && = \sum_{l\in[|\mathcal{S}_t|]} \sum_{k\in[K_t]} q_{t,l} \beta_{t,l,k} \partial \CVaR_{\alpha_k} (\underline{V}_{t, i}(\hat{x}_{t-1},\xi_t)|\xi_{t-1}=\tilde{\xi}_{t-1}).
% \edeq
The intercept $g_{t, i}$ can be interpreted 
geometrically 
as a convex combination of the intercepts of the cutting planes
of ``active pieces'' in the objective function of the Lagrange dual problem (we can imagine that the Lagrange dual is function of $x_{t-1}$ parameterized by $x_t$ and the active pieces means that
%when 
$x_t$  corresponds to the optimal solution of problem (\ref{eq:change-of-measure-based subproblem})).

\begin{remark}
There is another way to construct 
%the 
lower approximations of cost-to-go function $V_t (x_{t-1}, \tilde{\xi}_t)$ in  (\ref{eq-costToGoFunction-ARSRM})
by using the equivalent formulation of $\CVaR$:
\begin{equation}
\label{eq-CVaR-optimizationFormulation}
    \CVaR_{\alpha}(Z) 
    = \min_{\eta}
    \left\{\eta+\frac{1}{1-\alpha} \bbe [Z-\eta]_+
    % \sum_{j\in[K]} p_j (\xi_j-\eta)_+
    \right\}. 
\end{equation}
%then 
% {\color{blue}
It is well-known that the minimum is attained at 
$\inmat{VaR}_\alpha(Z)$.
% }.
Consequently, we can write $\mathcal{V}_{t+1}(x_t, \tilde{\xi}_t)$ 
in (\ref{eq-secondTerm})
as
%can be written as 
$$
\mathcal{V}_{t+1}(x_t, \tilde{\xi}_t) =
\min_{\bm{\eta}_t} 
% \left\{
\sum_{l\in[I_{t+1}]} \sum_{k\in[K_{t+1}]} q_{t+1,l}\beta_{t+1,l,k}
\left(\eta_{t, k}+ \frac{1}{1-\alpha_{t+1, k}} \bbe_{t+1|\tilde{\xi}_t} [(V_{t+1}(x_{t},\xi_{t+1})-\eta_t)_+] \right).
% \sum_{j\in[K]} p_j (\xi_j-\eta)_+)
% \CVaR_{\alpha_{t+1,k}} \left(V_{t+1}(x_{t},\xi_{t+1})| \xi_t=\tilde{\xi}_t \right) 
% \right\}.
$$
%{\color{red} 
Substituting $\mathcal{V}_{t+1}(x_t, \tilde{\xi}_t)$ into (\ref{eq-costToGoFunction-ARSRM}), we obtain 
% By combining the term $\sum_{l\in[I_{t+1}]} \sum_{k\in[K_{t+1}]} q_{t+1,l}\beta_{t+1,l,k} \eta_{t, k}$ with $c_t x_t$ and redefining the rest terms of $\mathcal{V}_{t+1}(x_t, \xi_{t+1})$ as $\mathcal{V}_{t+1}(x_t, \eta_t, \xi_{t+1})$, we obtain another equivalent reformulation of (\ref{eq-costToGoFunction-ARSRM}) as:
% {\color{red}HX: why aren't we interested in its lower approximation?}
% {\color{blue} Q: If we do this, it is just taking $\eta_t$ out and considering $\eta_t$ and $x_t$ as decision variable, and to approximate an expectation. The formulation will be more complicated and the computational efficiency is not very good. }
\begin{equation*}
% \label{eq-costToGoFunction-ARSRM}
    \begin{split}
        V_t (x_{t-1}, \tilde{\xi}_t ) := & \min_{x_t\in\mathcal{X}_t(x_{t-1},\tilde{\xi}_t), \eta_t} \; \Bigg\{c_t x_t + \sum_{l\in[I_{t+1}]} \sum_{k\in[K_{t+1}]} q_{t+1,l}\beta_{t+1,l,k} \eta_{t, k} \\
    & \qquad\qquad  + \sum_{l\in[I_{t+1}]} \sum_{k\in[K_{t+1}]}  
    \frac{q_{t+1,l}\beta_{t+1,l,k}}{1-\alpha_{t+1, k}} \bbe_{t+1|\tilde{\xi}_t} [(V_{t+1}(x_{t},\xi_{t+1})-\eta_t)_+] \Bigg\}.
    \end{split}
\end{equation*}
It is possible to construct lower approximations  of $V_t (x_{t-1}, \tilde{\xi}_t)$ for $t=1,\cdots,T$ based on the above dynamic recursive formulae. However,
we do not take this avenue because 
existing research shows that the CVaR formulations described as such may slow down the computational process, see 
\cite{STdCJS12,ZRG16}.

% This reformulation redefines the subproblem by adding a VaR representative variable $\eta_t$ to turn it into an risk neutral problem. 
% Here we use the change-of-measure-based decomposition to avoid the necessity to maintain a Value-at-Risk state variable $\eta_t$ at each stage, with the expense of computing $\lambda_t$ 
% {\color{red}HX: not defined in the above equation?}
% but no increase in the number of variables.
% Numerical simulations in \cite{STdCJS12,ZRG16} have shown that the corresponding subproblems are solved faster than maintaining the variable $\eta_t$ if the minimization formula from \cite{Ury00,RoU02} for the $\CVaR$ or the combination of expectation and $\CVaR$ is used.
% %We provide the lower bound algorithm in the following Algorithm 1.
% {\color{red}HX: The next algorithm computes the lower bounds of ??? BTW the algorithm computes both lower bounds and upper bounds.
% }
\end{remark}

Based on the discussions above, we propose 
an algorithm for computing the lower bounds of %the optimal value of
$V_t (x_{t-1}, \tilde{\xi}_t)$ in
(\ref{eq-costToGoFunction-ARSRM}). The algorithm consists of  a forward process and a backward process. The forward process computes optimal solutions of the lower approximation problems from stage $1$ to stage $T$. 
{\color{black} The backward process %iteratively 
constructs in each iteration 
a
%the 
hyperplane cut of the cost-to-go function. By convention, we set the initial cuts at each stage %by available information or as that 
with a sufficiently small 
$g_{t}^0$ %being smaller than 
below the true optimal value and $G_t^0$ being a zero vector.

\begin{algorithm}[H]
\caption{Construction of Lower Approximations and Computation of Lower Bounds }
\label{alg-lowerbound}
\SetAlgoLined

\KwIn { the maximum iteration $i^{\max}$, the size of sampling $N$,
% $g_{t,n}^0=0$, $G_{t, n}^0$ is the zero vector for $t=2,\ldots,T$ and $n=1,\ldots, N$, 
the initial cuts with $g_{2, n}^0$ and $G_{2,n}^0$ for $t=2,\ldots,T$ and $n=1,\ldots, N$.}
% {\color{red}HX: add a comment on how the initial cuts are set. }

\KwOut {the optimal solutions $\hat{x}^{i}_{t,n}$ for $t\in[T]$, $n\in[N]$, $i\in[i^{\max}]$ and lower bound $\underline{V}_{1, i^{\max}}$.}
%{\color{red} 
%if available.
%HX: what do you mean?} 
\While{$i=1,\ldots,i^{\max}$} 
{
$\rhd$ \textbf{Forward pass:} obtain decision $\hat{x}^{i}_{t,n}$, \text{for}\; $t=1,\ldots,T$, $n=1,\ldots,N$.

Solve the approximate stage $1$ problem 
$\underline{V}_{1, i} (x_0, \xi_1)$, obtain the optimal solution 
$\hat{x}^{i}_{1,n} =x_1^*$ for
$n\in[N]$ and optimal value $\underline{V}_{1, i}$ as a lower bound.

Draw $N$ samples of $\xi_t$, for $t=2,\ldots,T$. 

\For{stage $t=2,\ldots,T$ and sample paths $n=1,\ldots,N$}
{
Solve the approximate stage $t$ subproblem (\ref{eq-sub-re2}) with $\tilde{\xi}_t= \xi_{t, n}$,
% $\underline{V}_{t, i} (\hat{x}_{t-1,n}, \xi_{t, n})$ , 
obtain optimal solution $\hat{x}_{t, n}^i$. 
}
\vspace{0.3cm}
$\rhd$ \textbf{Backward pass:} construct new cuts and  update lower approximations at each stage. 

\For{stage $t=T,\ldots,2$, sample path $n=1,\ldots,N$}
{
    \For{scenario $j=1,\ldots,K_t$}
    {
    Solve the approximate stage $t$ subproblem  (\ref{eq-sub-re2}) with $\tilde{\xi}_t= \xi_{t, j}$.\
    % $\underline{V}_{t, i} (\hat{x}_{t-1,n}, \xi_{t, j})$.
    }
    Construct a
    %hyperplane 
    cut for stage $t-1$
    %: compute 
    with $G^i_{t, n}$ and
    $g^i_{t, n}$ defined as in (\ref{eq-cut-re2-a}) and (\ref{eq-cut-re2-b}).
}
}
\end{algorithm}
% \end{breakablealgorithm}

%-----------------------------------------------------------------------------------
\subsubsection{\texorpdfstring{Upper Bound of $V_t (x_{t-1}, \tilde{\xi}_t)$}{}}
%approximation of MARSRM}
\label{sec-innerApproximation-MARSRM}

{\color{black} 
We now turn 
%our attention 
to discuss construction of %building the
upper bound of $V_t (x_{t-1}, \tilde{\xi}_t)$.
%approximation of cost-to-go function defined in (\ref{eq-costToGoFunction-ARSRM}), which provides an upper bound of the optimal value of the subproblem (\ref{eq-costToGoFunction-ARSRM}).
% , which differs from the risk-neutral case.
}
% Observe that the tractable formulation (\ref{eq-ARSRM-sDiscretized}) of ARSRM $\rho_t$ depends on the scenario number of $X$, 
% which means the statistical upper bound can not be obtained. 
In the standard SDDP approach, 
the algorithm terminates if the gap between the lower bound 
%obtained from the lower approximation 
and the stochastic upper bound 
%becomes
falls 
%within 
%below
within a 
%fixed 
prescribed tolerance level, see \cite{Alex11}.
%for details.
Although a small gap may indicate that the current bound provides a good approximation of the true problem, it may cause early termination of the algorithm and subsequently
fail to provide a guaranteed  
%deterministic
quality of the current approximation (which would be of interest to the decision maker).
Moreover, for the risk averse case such as CVaR, 
%different path {\color{red}HX: do you mean different sample paths based on MC?
%}
%{\color{blue} Q: Yes.}
%by using Monte Carlo sampling scheme 
%would 
sample average approximation 
may incur more significant errors 
than risk neutral case.
%cause larger sampling error.
There are some 
%methods 
techniques to improve the statistic upper bounds
%{\color{blue} 
or stop 
%the 
iterations 
by using importance sampling or hypothesis testing respectively, see e.g. \cite{kozmik2015evaluating, homem2011sampling, garcia2023multicut}.
However, these techniques 
might be invalid or hard to implement in our problem because there could be multiple %$\CVaRs$
CVaRs
in ARSRM and the optimal values 
%from 
at the sampled paths 
%cannot
may provide enough evidence
to reject or accept the hypothesis.
%{\color{red}HX: explain why}
% {\color{blue}Q: if there is only one CVaR, we can do important sampling, but here is a combination of CVaR. }
This motivates us 
%we avoid such a statistic upper bound, instead we utilize the accurate 
to adopt a different approach \cite{park2023data} which calculates 
deterministic upper bound at each iteration
%upper bound proposed in 
%\cite{park2023data}.
%{\color{red} HX: What do you mean by accurate bound?}
%{\color{blue} Q: Statistical upper bound sometimes is not even an upper bound.}
Here we include a sketch of the method 
for self-containedness.

% Since the estimator of stage $t$ is calculated with one scenario, the corresponding ARSRM $\rho_t$ is also approximated.
% A more effective procedure available to compute an upper bound for the risk-averse case is proposed by Philpott et al. \cite{PdF13}.
% The authors develop an upper approximation scheme that provides a candidate policy and a deterministic upper bound on the policy value, using a convex combination of feasible policies and without sampling error.
% We should start from the last stage and recurse back to the first stage to obtain an estimator of the risk measure.
Suppose that 
at stage $t$ we have obtained 
upper bounds $\bar{V}_{t+1, 1},\ldots,\bar{V}_{t+1, N_t}$ 
%on 
of the values of $\rho_{Q_{t+1}}(V_{t+1}(x_t,\xi_{t+1}))$ at  $\hat{x}_{t,1},\ldots,\hat{x}_{t, N_t}$.
% {\color{red}HX: 
% is $\hat{x}_{t, N_t}$ obtained
% in the forward pass at the $N_t$-th iteration?
% }
% {\color{blue} Q: $\hat{x}_{t, N_t}$ is obtained from the forward pass with the $N_t$-th scenario $\xi_{N_t}$.}
% {\color{red}HX: In the algorithm, you were saying it provides upper bounds}
We use the points $\{(\hat{x}_{t, 1},\bar{V}_{t+1, 1}),\ldots,(\hat{x}_{t, N_t},\bar{V}_{t+1, N_t})\}$ 
in the epigraph of 
$\rho_{Q_{t+1}}(V_{t+1}(x_t,\xi_{t+1}))$
to construct 
% defines 
a convex polyhedral set $\mathcal{H}$. 
% that is 
% %a subset of 
% contained in 
% %the convex 
% the epigraph of $\rho_{Q_{t+1}}(V_{t+1}(x_t,\xi_{t+1}))$. 
% The
% %and the 
% optimal value of the following problem %is the 
% provides a lower bound of $\mathscr{H}$, which 
% %is 
% also provides  an upper bound %approximation to 
% for the value
% $\rho_{Q_{t+1}}(V_{t+1}(x_t,\xi_{t+1}))$:
By solving program
\begin{equation}
\label{eq-upper1}
    \min_{\bm{\theta}, x_t\geq 0}\;\; \sum_{n\in[N_t]}\theta_n\bar{V}_{t+1, n} , \quad 
    \st\;\; \sum_{n\in[N_t]}\theta_n \hat{x}_{t, n} =x_t, \sum_{n\in[N_t]}\theta_n=1, \theta_n \geq 0, 
    %m
\;  \text{for}\;   n\in[N_t],
\end{equation}
we obtain 
%and the 
%optimal value of the following problem %is the 
%provides 
the lowest point in set ${\cal H}$.
The optimal value of problem (\ref{eq-upper1})
%and 
%consequently 
provides an upper bound 
%the smallest value 
% a lower bound of $\mathscr{H}$, 
% which 
% %is 
% also provides 
% an upper bound %approximation to 
for the minimum of 
%value
$\rho_{Q_{t+1}}(V_{t+1}(x_t,\xi_{t+1}))$.
% where $\Theta_t:=\left\{\bm{\theta}\in\R^{N_t}_+ \;|\; \sum_{n\in[N_t]}\theta_n=1\right\}$.
%Then 
%Consequently 
Based on the discussions above, 
we consider the following upper bound on $c_tx_t+\rho_{Q_{t+1}}(V_{t+1}(x_t,\xi_{t+1}))$ under each realization of $\xi_t$ by solving 
\begin{subequations}
\label{eq-upper-bound-t}
\begin{align}
    \bar{V}_{t, n}(\tilde{\xi}_t)= \min_{\theta \in\Theta_t, x_t\geq 0}\; & \tilde{c}_t x_t+\sum_{n\in[N_t]}\theta_n\bar{V}_{t+1, n} + M_t\|y_t\|_1 \label{eq-upper-bound-t-obj} \\
    \st\quad\; & \tilde{A}_t x_t= \tilde{b}_t - \tilde{E}_t x_{t-1}, \label{eq-upper-bound-t-a} \\
    & \sum_{n\in[N_t]}\theta_n \hat{x}_{t, n} +y_t =x_t, \label{eq-upper-bound-t-b}
\end{align}
\end{subequations}
where $\widebar{V}_{T+1,n}=0$ for $n\in[N_{T+1}]$ and
% and the constant $M_t\geq 0$ 
 the constant $M_t$
is an upper bound of the Lipschitz modulus of $\bar{V}_{t+1,n}$ with respect to $x_t$.
%is not defined.}
At the first stage, we solve the following problem 
\begin{subequations}
\label{eq-upper-bound-1}
\begin{align}
    \bar{V}_1(\xi_1)= \min_{\theta \in\Theta_1, x_1\geq 0}\;\; & c_1 x_t+\sum_{n\in[N_1]}\theta_n\bar{V}_{2, n} +M_1\|y_1\|_1 \\
    \st \quad\;\;\; & A_1 x_1= b_1, \\
    & \sum_{n\in[N_t]}\theta_n x_{1, n} +y_1=x_1. 
\end{align}
\end{subequations}
%Note that t
The auxiliary variables $y_t$, $t\in[T]$ are introduced in (\ref{eq-upper-bound-t}) and (\ref{eq-upper-bound-1}) to ensure the feasibility of the upper approximation problem because the relatively complete recourse assumption does not guarantee that
a feasible solution $x_t\in\mathcal{X}_t(x_{t-1},\xi_t)$ also belongs to the convex combination \cite{park2023data}.
%then 
$y_t$ becomes nonzero if there does not exist $\theta$ that satisfies the constraints.
Therefore, a penalty proportional to a sufficiently large scalar $M_t$ is imposed in the case where $y_t$ is nonzero.
As the upper approximation 
% {\color{red}HX: I think we should replace inner/lower approximation by 
% upper/lower approximation throughout.}
$\bar{V}_{t,n}$ for all $t\in[T]\setminus\{1\}$ are calculated by the candidate solution from the backward process, the value $\bar{V}_{1}(x_1, \xi_1)$ provides tighter upper bound on the value of optimal policies along the iterations.
Algorithm \ref{alg-upperbound} states the upper bound computation scheme, which can be used in Algorithm \ref{alg-lowerbound} as the termination test.

%========================================================================================
% \newpage
% \vspace{0.5cm}
% \begin{breakablealgorithm}
% \caption{Upper bound computation}
% \label{alg-upperbound}
% \noindent
% \begin{algorithmic}[1]
% \For{$i = 1,\ldots,i^{\max}$, stage $t=T,\ldots,2$, scenario $j=1,\ldots,K_t$, $n=1,\ldots,N$}
% \State Solve (\ref{eq-upper-bound-t}) with $x_{t-1}=\hat{x}^{i}_{t-1,n}$ and $\xi_{t,j}$.
% \State Compute the optimal $\bm{\hat{\lambda}}_t$ that maximizes $\sum_{j\in[K_t]} p_{t, j} \hat{\lambda}_{t, k,j} \bar{V}_{t, n}(x_{t-1}, \xi_{t, j})$.
% \EndFor
% \State Solve (\ref{eq-upper-bound-1}) and save the optimal value as an upper bound.
% \end{algorithmic}
% \end{breakablealgorithm}

\begin{algorithm}[H]
\caption{Computation of Upper Bounds }
\label{alg-upperbound}
\SetAlgoLined
\KwIn {the optimal solutions $\hat{x}^{i}_{t,n}$ for $t\in[T]$, $n\in[N]$, $i\in[i^{\max}]$ obtained from Algorithm \ref{alg-lowerbound}.}

\KwOut {upper bound $\bar{V}_{1, i^{\max}}$.}

\For{
$i = 1,\ldots,i^{\max}$, stage $t=T,\ldots,2$, scenario $j=1,\ldots,K_t$, $n=1,\ldots,N$}
{
Solve (\ref{eq-upper-bound-t}) and obtain the optimal values $\bar{V}_{t, i}(x_{t-1}, \xi_1)$ with $x_{t-1}=\hat{x}^{i}_{t-1,n}$ and $\xi_{t,j}$, compute the optimal $\bm{\hat{\lambda}}_t$ that maximizes $\sum_{j\in[K_t]} p_{t, j} \hat{\lambda}_{t, k,j} \bar{V}_{t, n}(x_{t-1}, \xi_{t, j})$.
}

Solve (\ref{eq-upper-bound-1}) and obtain the optimal value $\bar{V}_{1, i}(x_0, \xi_1)$. 

% \KwOut {the optimal solutions $\hat{x}^{i}_{t,n}$ for $t\in[T]$, $n\in[N]$, $i\in[i^{\max}]$ and lower bound $\underline{V}_{1, i^{\max}}$.
% \end{algorithm}
\end{algorithm}

%=================================================================================

% \newpage
\section{Distributionally Preference Robust MARSRM}

% {\color{red}
% HX: We need to clarify at the very beginning whether
% the true probability distribution $Q^*$ is stage-dependent or not.
% }

%{\color{black}
In the previous section, we discuss the multistage decision making model 
with the DM's risk preferences at each stage 
%the risk measure of the 
%DM 
being characterized by ARSRM.
%and t
The probability distribution $Q_t$ of the 
random preference parameter $s$ is known and fixed at each stage, which 
means that
%implies 
the possible states of DM's risk preference 
can be observed completely. 
In practice, however, a DM's risk preference and 
their uncertain features may vary 
%with 
from stage to stage, in which case it would 
be difficult to predict $Q_t$ due to lack of complete information of risk preference.
In this section, we 
discuss the case that
%discuss the situation where 
%the DM's risk preference
%is ambiguous, that is, 
there is incomplete information to identify the true stage-dependent probability distribution $Q^*_t$ of random preference parameter $s$.
This motivates us to consider 
%the DRO 
a distributionally robust version of MARSRM where
%while 
the ambiguity arises from incomplete information on DM’s risk preferences, rather than 
exogenous uncertainty in 
%common DRO 
in the traditional distributional robust optimization models. 

\subsection{Distributional Robust MARSRM Model}
Specifically, 
we consider the distributionally preference robust version of multistage decision making (\ref{eq-MARSRM}) (DR-MARSRM)
which is defined as follows:
\begin{equation}
\label{eq-DR-MARSRM}
\begin{split}
    \min_{x_1\in \mathcal{X}_1(x_0, \xi_1)}  c_1x_1 + 
    & \sup_{Q_2\in\mathcal{Q}_2} 
    \rho_{Q_2|\mathcal{F}_1}  \Bigg[ \inf_{x_2\in\mathcal{X}_2(x_1, \xi_2)} c_2x_2 +\cdots \\
    &  + \sup_{Q_{T-1}\in\mathcal{Q}_{T-1}} \rho_{Q_{T-1}|\mathcal{F}_{T-2}} 
    \bigg[ \inf_{x_{T-1}\in\mathcal{X}_{T-1}(x_{T-2}, \xi_{T-1})} c_{T-1}x_{T-1} \\
    & \qquad\qquad\qquad + \sup_{Q_T\in\mathcal{Q}_T} 
      \rho_{Q_{T}|\mathcal{F}_{T-1}}
    \Big[ \inf_{x_T\in\mathcal{X}_T(x_{T-1}, \xi_T)} c_Tx_T \Big]\bigg]\Bigg],
\end{split}
\end{equation}
%{\color{red}HX: Is the problem time consistent?}
where 
$\mathcal{Q}_t$ is the ambiguity set of random preference parameter $s_t$ with support set $\mathcal{S}_t$,  $\rho_t$ and $\mathcal{X}_t$ are defined as in (\ref{eq-ARSRM}) and 
(\ref{eq-mathcalX_t}) for $t\in[T]\setminus\{1\}$.
We make similar assumptions for this multistage decision-making problem 
%as
to Assumption \ref{assumption-multistageproblem} to 
%guarantee the 
ensure well-definedness of DR-MARSRM.
%{\color{blue}Q: Note that the optimal policy of (\ref{eq-DR-MARSRM}) is time-consistent since the multistage risk measure is in a nested form, see e.g. \cite{Rus10,delage2015robust,shapiro2016rectangular}.
% } 
% {\color{red}HX: I thought time consistency implies nested form but not conversely. Can we write (\ref{eq-DR-MARSRM}) as
% \begin{equation}
% \label{eq-DR-MARSRM-HX-Q}
% \begin{split}
%     \min_{x_1\in \mathcal{X}_1(x_0, \xi_1), x_2,\cdots, x_{T-1}}  c_1x_1 + 
%     & \sup_{Q_2\in\mathcal{Q}_2} 
%     % \rho_{Q_2|\mathcal{F}_1}  \Bigg[ \inf_{x_2\in\mathcal{X}_2(x_1, \xi_2)} 
%     c_2x_2 +\cdots \\
%     &  + \sup_{Q_{T-1}\in\mathcal{Q}_{T-1}} \rho_{Q_{T-1}|\mathcal{F}_{T-2}} 
%     \bigg[ 
%     %\inf_{x_{T-1}\in\mathcal{X}_{T-1}(x_{T-2}, \xi_{T-1})}
%     c_{T-1}x_{T-1} \\
%     & \qquad\qquad\qquad + \sup_{Q_T\in\mathcal{Q}_T} 
%       \rho_{Q_{T}|\mathcal{F}_{T-1}}
%      \Big[ %inf_{x_T\in\mathcal{X}_T(x_{T%-1}, \xi_T)} 
%     c_Tx_T \Big]\bigg]\Bigg],
% \end{split}
% \end{equation}
% }
In the first stage, after observing the realization $\xi_1$,
% {\color{red}HX: do you mean $\xi_1$, I cannot find $\xi_0$ in the formula}
the nature chooses the worst-case probability distribution $Q_2$,
% {\color{red} HX: $Q_2$? The notation is not consistent with the previous equation.}
under which the decision $x_1$ is made. After the uncertain parameter $\xi_1$ is observed, the worst-case probability distribution $Q_2$ is chosen, and then the corresponding decision $x_2$ is made. The rest of the  process continues until reaching stage $T$. 
The dynamic decision-making process 
%is
can be described in the following flowchart:
%as follows:
\begin{equation*}
\begin{split}
& \underbrace{ \text { observation } \xi_0 
% \rightarrow \text{worst-case } Q_{1}  
\rightarrow \text { decision } \boldsymbol{x}_1 }_{\text {Stage 0 (current stage, all deterministic)} } %\rightarrow  
\\
& \qquad
\rightarrow \underbrace{ \text { observation }\xi_1 
\rightarrow \text{worst-case } Q_{2}  \rightarrow \operatorname{decision} x_2 }_{\text {Stage } 1 (\mathcal{F}_{1}-\text{measurable})} \rightarrow \cdots \\
& \qquad
\rightarrow \underbrace{ \text { observation }\xi_{t-1} 
\rightarrow \text{worst-case } Q_{t} \rightarrow \operatorname{decision} x_{t} }_{\text {Stage } t-1\ (\mathcal{F}_{t-1}-\text{measurable}) } \rightarrow \cdots \\
& \qquad
\rightarrow \underbrace{ \text { observation } \xi_{T-1} 
\rightarrow \text{worst-case } Q_{T} 
\rightarrow \operatorname{decision} x_T }_{\text {Stage } T-1 (\mathcal{F}_{T-1}-\text{measurable})} 
\rightarrow \underbrace{ \operatorname{observation} \xi_T}_{\text {Stage } T (\mathcal{F}_{T}-\text{measurable})}
\end{split}
\end{equation*}
{\color{black}
At this point, it might be helpful to relate (\ref{eq-DR-MARSRM}) to
the multistage utility preference robust optimization (MUPRO) model in \cite{liu2021multistage}.
%Part from obvious difference where \cite{liu2021multistage}
First, in \cite{liu2021multistage}, 
the DM's risk preference at each stage satisfies 
the four axiomatic properties of von Neumann-Morgenstern's expected utility theory which means that the DM's preference can be characterized by a unique deterministic utility function up to positive linear transformation. 
%whereas
In contrast, 
here (in (\ref{eq-DR-MARSRM})), the DM's preferences at each stage do not satisfy the four  axiomatic properties of the dual expected utility theory 
which means the DM's preferences cannot  be described by a unique SRM. 
Second, in \cite{liu2021multistage},
%each stage 
the DM's optimal decision at stage $t$ is made before realization of $\xi_t$  based on the worst-case expected utility of the random reward function,  
whereas here the optimal decision is made after observation of uncertainty $\xi_t$ and  based on the worst-case distribution $Q_{t+1}$, which means that the optimal decision might be less conservative.
%moreover, 
Third, by adopting ARSRM,  the worst-case distribution $Q_t$
is state/historical path independent whereas the worst-case utility function in \cite{liu2021multistage} is state/historical path dependent. 
This will significantly simplify numerical procedures for solving
problem (\ref{eq-DR-MARSRM}).
Fourth, the nested form of model (\ref{eq-DR-MARSRM}) is set-up without considering its relation to the 
holistic minimax problem where
the optimal decisions at all stages are made 
at stage one
based on 
the sequence of the worst-case probability distributions of $Q_1,\cdots,Q_{T-1}$, 
%are considered 
%at the beginning of first stage, 
this means that we are not concerned about rectagularity of the sequence of 
ambiguity sets 
and time consistency of the sequence of optimal solutions as in \cite{liu2021multistage}. Indeed, if we adopt 
the specific SRM (\ref{eq:SRM-comb-Exp-CVar})
at each stage, 
the resulting DR-MARSRM model is time inconsistent because 
%the model
CVaR does not satisfy the tower property.

Throughout this section,
%paper,
we assume that the realizations of
$s_t$ (representing the state of the DM's risk preferences) are observable, 
%from 
%at least 
e.g., based on the 
%past information
empirical data 
but the corresponding 
probabilities that these states may 
occur again is unknown. In other words,
%that is, the realizations 
the support set of the random variable 
$s_t$ 
%are 
is known, 
% {\color{red}HX: the support set of $s_t$ is unknown? 
% This is not the set-up of distributional ambiguity.
% }
but the probability $P(s_t=s_{t, l})$, 
%unknown for 
$l\in[|\mathcal{S}_t|]$, is ambiguous.
Since the support set of 
%$S$ 
$s_t$
is finite, 
we can formulate the robust cost-to-go function as
\begin{equation}
\label{eq-costToGoFunction-PDRO}
\resizebox{0.92\linewidth}{!}{$
    V_t^R(x_{t-1},  \tilde{\xi_{t}})
    := \displaystyle{\min_{x_t\in\mathcal{X}_t(x_{t-1}, \tilde{\xi}_t)}} \left\{ c_t x_t  
    + \sup_{\bm{q}_{t+1}\in \mathcal{Q}_{t+1}} 
    \sum_{k\in[K_{t+1}]} \sum_{l\in[|\mathcal{S}_{t+1}|]}  q_{t+1,l}\beta_{t+1,l,k} \CVaR_{\alpha_k}(V_{t+1}^R(x_{t},\xi_{t+1})|\xi_{t}=\tilde{\xi_{t}}
    )\right\}, 
    $}
\end{equation}
where $\mathcal{Q}_{t+1}$ is 
% {\color{blue}
an ambiguity set of probability distributions 
%with 
over the support set
$\mathcal{S}_{t+1}$. 
%which is a subset of the simplex with proper dimension.
Moreover, we can assume that 
%observe 
the 
%empirical
sample mean $\mu_t$ and sample 
covariance matrix  $\Sigma_t$ of $s_t$ 
%from previous decision
are obtainable, and this 
%which allows 
enables 
us to define the stagewise 
ambiguity sets 
%by matching 
% the estimates of the
with the estimated 
%true 
sample mean value  and 
sample covariance.
These samples may be obtained 
from empirical data or
computer simulation.
%matrix. 
%{\color{red}HX: What do you mean?}
Specifically, we consider the moment-based ambiguity set with finite support for characterizing $\mathcal{Q}_t$:
\begin{equation}
\label{eq-mathcalQ_t}
\resizebox{0.92\linewidth}{!}{$
\displaystyle{
    \mathcal{Q}_t := \left\{q_t\in\R_+^{|\mathcal{S}_t|} \left|\, \sum_{l\in[|\mathcal{S}_t|]} q_{t,l}=1, \sum_{l\in[|\mathcal{S}_t|]}q_{t,l} s_{t, l}=\mu_t, \sum_{l\in[|\mathcal{S}_t|]} q_{t, l}(s_{t,l}-\mu_t)(s_{t,l}-\mu_t)^{\top}
    =\Sigma_t \right.\right\},
    }
$}
\end{equation}
where
$|\mathcal{S}_t|$ is the sample size,
% number of scenarios of 
% $s_{t,l}$, 
% $\mu_t$ and $\Sigma_t$ 
% %is 
% are 
% %the
% %empirical 
% subjective 
% mean value  and covariance 
% of $s_t$ at stage $t$, 
see similar 
%ideas
formulations
in the definition of 
good-deal measure 
\cite{cochrane2000beyond,druenne2011good}.
% This can be understood as 
% estimated from empirical data including information on $\mu_{t-1}$ and $\Sigma_{t-1}$. 
% {\color{red}HX: how can we obtain these information when we stand 
% at the beginning of stage $t$? 
% A more reasonable way would be consider these from the previous stage.
% }
%and 
% $|\mathcal{S}_t|$ is the number of scenarios of 
% $s_{t,l}$. 
% {\color{blue}In other words, the probability distribution of $s_t$ is inferred from the mean and variance of $s_{t-1}$.}

\subsection{SDDP for Solving DR-MARSRM}

% {\color{blue}
In this section,
we consider DR-MARSRM
(\ref{eq-DR-MARSRM}) with ambiguity set $\mathcal{Q}_t$ being defined in 
%the moment-based ambiguity set of probability distribution of $s$ as in 
(\ref{eq-mathcalQ_t}).
As in the non-robust case discussed in Section 3,
we propose to solve the problem 
by SDDP. To this end, we 
%and
derive 
%corresponding outer and upper 
upper and lower 
approximations of the objective functions at each stage in an iterative process.
% for solving distributionally preference robust MARSRM (\ref{eq-DR-MARSRM}).
To tackle the additional challenges 
arising from the worst-case ARSRM at each stage, we %begin by 
need to derive the Lagrange dual of the
inner maximization problem (with $Q_t$ at stage $t$) and merge it into the outer minimization problem (with respect to $x_t$). 
% }

% {\color{red}HX: Literature review of SDDP for DRO
%At this point, we 
It might be helpful to note at this point that 
various DRO models for multistage decision making
have been developed where the probability distributions 
of the exogenous uncertainties ($\xi_t$ in this context) are ambiguous.
% Several approaches to solve MSP in distributionally robust setting by using SDDP scheme have been studied in the literature. 
% The ways of constructing t
The ambiguity sets are constructed using
%can 
%base on 
empirical moments 
%(see e.g. 
(\cite{yu2022multistage}) or statistical distance such as $\infty$-norm (\cite{huang2017study}), 
modified $\chi^2$ (\cite{philpott2018distributionally}), Wasserstein distance (\cite{pflug2014multistage, duque2020distributionally}), nested distance (\cite{gao2024data}), total variance (\cite{silva2021data}).
The main difference between  our DR-MARSRM model and 
the multistage DRO models listed above is that here we focus 
on ambiguity of endogenous uncertainty which is related to the decision maker's risk preference.  This difference poses new challenges 
for us to handle when 
we develop lower and upper approximations 
%results in subsequent differences 
%when we apply 
in the application of the SDDP method to solve problem (\ref{eq-DR-MARSRM}).

% {\color{red} More fundamental differences??

% To handle risk-aversion in the problem, we need to construct the approximation of a distributionally robust ARSRM 
% }

% Different from these models, 
% our model considers that the decision maker is described by an ARSRM, rather than expectation which is risk-neutral.
% Moreover, the uncertainty in our model is on the endogenous random parameter, that is, on the preference parameter $s_t$. 
% Park et al. \cite{park2023data} propose a data-driven scheme to solve distributionally robust MSLP where the ambiguity set is constructed as a polyhedral approximation of \textit{modified} $\chi^2$ by using using the Nadaraya-Watson regression.

\subsubsection{Moment-ased Ambiguity Set and Tractable Formulation of (\ref{eq-costToGoFunction-PDRO})}
%and 
%lower/upper approximations of $V_t^R(x_{t-1}, \tilde{\xi_{t}})$}

% Given the moment-based ambiguity set (\ref{eq-mathcalQ_t}), we first reformulate the inner robust maximization problem of
Recall that $V_t^R(x_{t-1}, \tilde{\xi_{t}})$ 
is the optimal value of the minimization problem
defined as in (\ref{eq-costToGoFunction-PDRO}).
The objective function of the problem comprises
a linear term $c_tx_t$ and 
a nonlinear term
% The optimal value is based
% The worst-case value of 
% $V_t^R(x_{t-1}, \tilde{\xi_{t}})$ with respect to
% $\bm{q}_{t+1}\in \mathcal{Q}_{t+1}$
% is defined as 
\begin{equation}
\label{eq-innerMaximization-PDRO}
    % \mathcal{V}_t^R (x_t,\xi_{t+1}) = 
    \sup_{\bm{q}_{t+1}\in \mathcal{Q}_{t+1}} \;
    \sum_{l\in[|\mathcal{S}_{t+1}|]} \sum_{k\in[K_{t+1}]} q_{t+1,l}\beta_{t+1,l,k} \CVaR_{\alpha_{t+1,k}}(V_{t+1}^R (x_{t},\xi_{t+1})|\xi_{t}=\tilde{\xi_{t}}),
\end{equation}
which is a convex function of $x_t$.
%Since this term in the cost-to-go function (\ref{eq-costToGoFunction-PDRO}) is an optimization problem rather than a risk measure as in %non-PDRO version decision making problem MARSRM 
%problem (\ref{eq-MARSRM}), w
We can use the convex minimization formulation of $\CVaR$ 
%as defined 
in (\ref{eq-CVaR-optimizationFormulation}) to %integrate the solving process to 
reformulate the term as a single linear program.
%and subsequently
%avoid the ordering step 
%at each stage 
%discussed earlier.
The next proposition states this.

\begin{proposition}
\label{prop-reformulationSecondTerm}
Let $\mathcal{Q}_{t+1}$ be 
%the ambiguity set as 
defined as in (\ref{eq-mathcalQ_t}). If 
 $\mathcal{Q}_{t+1}$
%the ambiguity set defined in (\ref{eq-mathcalQ_t}) 
is non-empty, then (\ref{eq-innerMaximization-PDRO})
is equivalent to the following linear program:
\begin{equation}
\label{eq-innermax-dro-one}
\begin{split}
    \inf_{\substack{\bm{\zeta}_{t+1}, \bm{\eta}_{t+1}, \\ \bm{\Delta}_{t+1}}} \; & \zeta^1_{t+1} + \mu_{t+1}\zeta^2_{t+1}  + \zeta^3_{t+1}\cdot \Sigma_{t+1} \\
    \st \quad\; & \sum_{k\in[K_{t+1}]} \beta_{t+1, l, k}
    \Bigg(\eta_{t+1, k}+\frac{1}{1-\alpha_{t+1, k}} \sum_{j\in[K_{t+1}]} p_{t+1,j} \Delta_{t+1, k, j} \Bigg) -\zeta^1_{t+1} - s_{t+1,l}\zeta^2_{t+1} \\
    & \hspace{6em} - \zeta^3_{t+1} \cdot (s_{t+1,l}-\mu_{t+1})(s_{t+1,l}-\mu_{t+1})^{\top} \leq 0, \quad %\forall
    \text{for}\;
    l\in[|\mathcal{S}_{t+1}|], \\
    & V_{t+1}^R (x_{t},\xi_{t+1, j}) - \eta_{t+1, k} \leq \Delta_{t+1, k,j}, \quad 
    %\forall
    \text{for}\;
    k, j\in[K_{t+1}],  \\
    & \Delta_{t+1, k, j} \geq 0, \quad  
    %\forall
    \text{for}\; k, j\in[K_{t+1}], 
\end{split}
\end{equation}
where
$\zeta^3_{t+1}\cdot \Sigma_{t+1}$ denotes the Frobenius inner product of two matrices. 
%trace.
\end{proposition}

\proof
By substituting 
(\ref{eq-CVaR-optimizationFormulation}) 
into (\ref{eq-innerMaximization-PDRO}), we obtain
\begin{align}
    & \max_{\bm{q}_{t+1}\in\mathcal{Q}_{t+1}} \min_{{\bm \eta}_{t+1}} \;\; v(\bm{q}_{t+1}, \bm{\eta}_{t+1}) \label{eq-inner-maximin} \\
    & := \sum_{l\in[|\mathcal{S}_{t+1}|]} \sum_{k\in[K_{t+1}]} q_{t+1, l} \beta_{t+1, l, k}
    \Bigg( \eta_{t+1, k}+\frac{1}{1-\alpha_{t+1, k}} \sum_{j\in[K_{t+1}]} p_{t+1,j} \left(V_{t+1}^R (x_{t},\xi_{t+1, j})-\eta_{t+1, k}\right)_+ \Bigg). \nonumber
\end{align}
%{\color{red} 
%By {\color{blue} 
By Assumption \ref{assumption-multistageproblem} (a),
%HX: which assumption?} 
$V_{t+1}^R (x_{t},\xi_{t+1, j})$
is bounded for all $j\in[K_{t+1}]$, hence
the optimal $\bm{\eta}^*_{t+1}$ is bounded. 
Since $v(\bm{q}_{t+1}, \bm{\eta}_{t+1})$ is convex in $\bm{\eta}_{t+1}$ for all $\bm{q}_{t+1}$ and linear in $\bm{q}_{t+1}$ for all 
$$\eta_{t+1, k}\in \left[\essinf V_{t+1}^R (x_{t},\xi_{t+1}),  \esssup V_{t+1}^R (x_{t},\xi_{t+1})\right]\subset\R,$$ 
and $\mathcal{Q}_{t+1}$ is convex, it follows 
%from 
by Sion's minimax theorem (see e.g. \cite{sion1958general}),
(\ref{eq-inner-maximin}) is equivalent to 
\begin{equation*}
\label{eq-inner-minimax}
\resizebox{\hsize}{!}{
$
    \displaystyle{ \min_{{\bm \eta}_{t+1}} \max_{\bm{q}_{t+1}\in\mathcal{Q}_{t+1}} } \;\; \sum_{i\in[|\mathcal{S}_{t+1}|]} \sum_{k\in[K_{t+1}]} q_{t+1, i} \beta_{t+1, i, k}
    \Bigg(\eta_{t+1, k}+\frac{1}{1-\alpha_{t+1, k}} \sum_{j\in[K_{t+1}]} p_{t+1,j} \left(V_{t+1}^R (x_{t},\xi_{t+1, j})-\eta_{t+1, k} \right)_+ \Bigg).
$
}
\end{equation*}
Furthermore, we have an equivalent reformulation of the above problem given by
\begin{align}
    \min_{{\bm \eta}_{t+1}, \bm{\Delta}_{t+1}} & \max_{\bm{q}_{t+1}\in\mathcal{Q}_{t+1}} \; \sum_{l\in[|\mathcal{S}_{t+1}|]} \sum_{k\in[K_{t+1}]} q_{t+1, l} \beta_{t+1, l, k}
    \Bigg(\eta_{t+1, k}+\frac{1}{1-\alpha_{t+1, k}} \sum_{j\in[K_{t+1}]} p_{t+1,j} \Delta_{t+1, k, j} \Bigg) \nonumber \\
    \st \quad & V_{t+1}^R (x_{t},\xi_{t+1, j}) - \eta_{t+1, k} \leq \Delta_{t+1, k,j}, \quad \forall k, j\in[K_{t+1}], \label{eq-prop-PDRO-1} \\
     & \Delta_{t+1, k, j} \geq 0, \quad  \text{for}\; k, j\in[K_{t+1}]. \nonumber 
\end{align}
By the definition of ambiguity set (\ref{eq-mathcalQ_t}), the inner maximization problem of (\ref{eq-prop-PDRO-1}) can be formulated as the following linear program:
\begin{subequations}
\label{eq-prop-PDRO-3}
\begin{align}
        \sup_{\bm{q}_{t+1}\geq0} \; & \sum_{l\in[|\mathcal{S}_{t+1}|]} \sum_{k\in[K_{t+1}]} q_{t+1, l} \beta_{t+1, l, k}
        \Bigg(\eta_{t+1, k}+\frac{1}{1-\alpha_{t+1, k}} \sum_{j\in[K_{t+1}]} p_{t+1,j} \Delta_{t+1, k, j} \Bigg) \\ 
        \st \;\;\; & \sum_{l\in[|\mathcal{S}_{t+1}|]} q_{t+1, l} =1, \\
        & \sum_{l\in[|\mathcal{S}_{t+1}|]} q_{t+1, l} s_{t+1, l} =\mu_{t+1}, \\
        & \sum_{l\in[|\mathcal{S}_{t+1}|]} q_{t+1, l} (s_{t+1, l}-\mu_{t+1})(s_{t+1, l}-\mu_{t+1})^{\top} =\Sigma_{t+1}.
\end{align}
\end{subequations}
We associate the dual variables $\zeta^1_{t+1}\in\R$, $\zeta^2_{t+1}\in\R^{d_s}, \zeta^3_{t+1}\in\R^{d_s \times d_s}$ with constraints 
%of 
(\ref{eq-prop-PDRO-3}b)- 
(\ref{eq-prop-PDRO-3}d)
% % {\color{red}HX: This is a bit vague. 
% % We had better label all three constraints
% % }
and 
write down the Lagrange dual of the above problem: 
\begin{equation}
\label{eq-prop-PDRO-2}
    \begin{split}
        \min_{\bm{\zeta}_{t+1}} \; & \zeta^1_{t+1} + \mu_{t+1} \zeta^2_{t+1} + \zeta^3_{t+1}\cdot \Sigma_{t+1} \\
        \st \; & \sum_{k\in[K_{t+1}]} \beta_{t+1, l, k}
        \Bigg(\eta_{t+1, k}+\frac{1}{1-\alpha_{t+1, k}} \sum_{j\in[K_{t+1}]} p_{t+1,j} \Delta_{t+1, k, j} \Bigg) -\zeta^1_{t+1} - s_{t+1, l} \zeta^2_{t+1}  \\
        & - \zeta^3_{t+1} \cdot (s_{t+1, l}-\mu_{t+1})(s_{t+1, l}-\mu_{t+1})^{\top} \leq 0, \quad \text{for}\; l\in[|\mathcal{S}_{t+1}|].
    \end{split}
\end{equation}
Combining (\ref{eq-prop-PDRO-1}) and (\ref{eq-prop-PDRO-2}), we obtain (\ref{eq-innermax-dro-one}).
\hfill\Box

By substituting (\ref{eq-innermax-dro-one}) into the 
%cost-to-go function 
the right hand side minimization problem in
(\ref{eq-costToGoFunction-PDRO}), we can obtain %its 
tractable formulation of the minimization problem 
as follows:
\begin{subequations}
\label{eq-costToGoFunction-PDRO-tractable}
    \begin{align}
        V_t^R (x_{t-1}, \tilde{\xi}_{t}) & \nonumber \\
        = \inf_{\substack{x_t, \bm{\zeta}_{t+1}, \\ \bm{\eta}_{t+1}, \bm{\Delta}_{t+1}}} \; & c_t x_t +\zeta_{t+1}^1 + \mu_{t+1} \zeta^2_{t+1} + \zeta_{t+1}^3\cdot \Sigma_{t+1} \\
        \st \quad\; & \sum_{k\in[K]} \beta_{t+1, l, k}
        \Bigg(\eta_{t+1, k}+\frac{1}{1-\alpha_{t+1, k}} \sum_{j\in[K_{t+1}]} p_{t+1,j} \Delta_{t+1, k, j} \Bigg) -\zeta_{t+1}^1 -s_{t+1,l}\zeta_{t+1}^2 \nonumber \\
        & \hspace{5em} - \zeta_{t+1}^3 \cdot (s_{t+1,l}-\mu_{t+1})(s_{t+1,l}-\mu_{t+1})^{\top} \leq 0, \quad \text{for}\; l\in[|\mathcal{S}_{t+1}|], \\
        & V_{t+1}^R (x_t, \xi_{t+1, j}) - \eta_{t+1, k} \leq \Delta_{t+1, k, j}, \quad \text{for}\; k, j \in[K_{t+1}], \\
        & \tilde{A}_t x_t=\tilde{b}_t-\tilde{E}_t x_{t-1}, \label{constraint-robust-staget} \\
        & x_t \geq 0, \Delta_{t+1, k, j} \geq 0, \quad  \text{for}\; k, j\in[K_{t+1}]. 
    \end{align}
\end{subequations}
%Now

\subsubsection{
%Moment-based ambiguity set and 
\texorpdfstring{Lower Bound of $V_t^R(x_{t-1}, \tilde{\xi_{t}})$}{}
}

Next, we 
%can 
propose an iterative approach 
which generates a sequence of 
lower bounds
of $ V_t^R (x_{t-1}, \tilde{\xi}_{t})$
%by developing 
based on lower linear 
envelop of $V_{t+1}^R (x_t, \xi_{t+1, j})$:
%build the piecewise linear outer and upper approximations to obtain the lower and upper bound of $V_t^R (x_{t-1}, \tilde{\xi}_{t})$ with (\ref{eq-costToGoFunction-PDRO-tractable}) which are iteratively refined in the algorithm.
% {\color{red} The lower approximation of the $i$th iteration to approximate $V_t^R (x_{t-1}, \tilde{\xi}_{t})$ can be defined as follows:
% HX: Not very clear, please explain why it is a lower approximation. 
% }
\begin{subequations}
    \label{eq-costToGoFunction-PDRO-tractable-lower}
    \begin{align}
        \underline{V}_{t, i}^R(x_{t-1}, \tilde{\xi}_{t}) & \nonumber \\
        = \inf_{\substack{x_t, \bm{\zeta}_{t+1}, \\ \bm{\eta}_{t+1}, \bm{\Delta}_{t+1}}} \; & c_t x_t +\zeta_{t+1}^1 + \mu_{t+1} \zeta_{t+1}^2 + \zeta_{t+1}^3\cdot \Sigma_{t+1} \\
        \st \quad\; & \sum_{k\in[K]} \beta_{t+1, l, k}
        \Bigg(\eta_{t+1, k}+\frac{1}{1-\alpha_{t+1, k}} \sum_{j\in[K_{t+1}]} p_{t+1,j} \Delta_{t+1, k, j} \Bigg) -\zeta_{t+1}^1- s_{t+1, l}\zeta_{t+1}^2 \nonumber \\
        & \hspace{2em} - \zeta_{t+1}^3 \cdot (s_{t+1, l}-\mu_{t+1})(s_{t+1, l}-\mu_{t+1})^{\top} \leq 0, \quad \text{for}\; l\in[|\mathcal{S}_{t+1}|], \\
        & g_{t+1, j, n} +G_{t+1, j, n} x_t - \eta_{t+1, k} \leq \Delta_{t+1, k, j}, \quad \text{for}\; n\in\mathcal{N}_{t+1, i}, k, j \in[K_{t+1}], \label{eq-eq-costToGoFunction-PDRO-tractable-lower-c} \\
        & \tilde{A}_t x_t=\tilde{b}_t-\tilde{E}_t x_{t-1}, \label{constraint-robust-staget-lin} \\
        & x_t \geq 0, \Delta_{t+1, k, j} \geq 0, \quad  \text{for}\; k, j\in[K_{t+1}], 
    \end{align}
\end{subequations}
where $\mathcal{N}_{t+1, i}$ is the index set of previous cuts in the first $i$ iterations, $G_{t+1, j, n}$ and $g_{t+1, j, n}$ are the subgradient and intercept of $V_{t+1}^R (x_t, \xi_{t+1, j})$.
% {\color{blue}
Since 
$
V_{t+1}^R (x_t, \xi_{t+1, j}) \geq 
g_{t+1, j, n} +G_{t+1, j, n} x_t$, then 
inequality 
(\ref{eq-costToGoFunction-PDRO-tractable}c)
implies (\ref{eq-costToGoFunction-PDRO-tractable-lower}c)
which means the feasible set of problem
(\ref{eq-costToGoFunction-PDRO-tractable}) is smaller 
that of problem (\ref{eq-costToGoFunction-PDRO-tractable-lower}) and hence $\underline{V}_{t, i}^R(x_{t-1}, \tilde{\xi}_{t})$ provides a lower bound of 
$       V_t^R (x_{t-1}, \tilde{\xi}_{t})$.
% }

Let $\varrho_{t, j}$ be the 
Lagrange multiplier of  constraint (\ref{constraint-robust-staget-lin}) with $\tilde{A}_t = A_{t, j}, \tilde{b}_t = b_{t, j}$ and $\tilde{E}_t = E_{t, j}$ for $j\in[K_t]$. Then a new cut can be constructed based on 
a subgradient of $V_{t+1}^R (\cdot, \xi_{t+1, j})$ at $\hat{x}_{t-1}$
and the intercept 
%is 
%are defined as
\begin{equation}
\label{eq-cut-DRO}
% \begin{split}
    G_{t, j} = -\varrho_{t, j} E_{t, j}, \;
    g_{t, j} = \underline{V}^R_{t, i} (\hat{x}_{t-1}, \xi_{t, j}) -G_{t, j}\hat{x}_{t-1}.
    % \varrho_{t, j} E_{t, j} \hat{x}_{t-1}.
% \end{split}
\end{equation}
% {\color{red}HX: How is the intercept figured out?}
% {\color{blue} Q: By writing the Lagrange dual problem, The optimal value $\underline{V}^R_{t, i} (\hat{x}_{t-1}, \xi_{t, j})$ is equal to the objective of (4.11) plus other non-zero term, other terms with variables all equal to zero.}

\subsubsection{
%approximations 
\texorpdfstring{Upper Bound 
of $V_t^R(x_{t-1}, \tilde{\xi_{t}})$}{}}

Analogous to (\ref{eq-upper-bound-t}), we can also construct a program which gives rise to
an upper bound of 
% By using similar scheme for constructing the upper approximation of MARSRM as discussed in Section \ref{sec-innerApproximation-MARSRM}, 
% the corresponding upper approximation of the $i$th iteration to obtain the upper bound of
$V_t^R (x_{t-1}, \tilde{\xi}_{t})$ 
%can be defined as:
\begin{subequations}
    \label{eq-costToGoFunction-PDRO-tractable-upper}
    \begin{align}
        \bar{V}_{t, i}^R(x_{t-1}, \tilde{\xi}_{t}) & \nonumber \\
        = \inf_{\substack{x_t, \bm{\zeta}_{t+1}, \\ \bm{\eta}_{t+1}, \bm{\Delta}_{t+1}}} \; & c_t x_t +\zeta_{t+1}^1 + \mu_{t+1} \zeta_{t+1}^2 + \zeta_{t+1}^3\cdot \Sigma_{t+1} \\
        \st \quad\; & \sum_{k\in[K]} \beta_{t+1, l, k}
        \Bigg(\eta_{t+1, k}+\frac{1}{1-\alpha_{t+1, k}} \sum_{j\in[K_{t+1}]} p_{t+1,j} \Delta_{t+1, k, j} \Bigg) -\zeta_{t+1}^1- s_{t+1, l} \zeta_{t+1}^2 \nonumber \\
        & \hspace{2em} - \zeta_{t+1}^3 \cdot (s_{t+1, l}-\mu_{t+1})(s_{t+1, l}-\mu_{t+1})^{\top} \leq 0, \quad \text{for}\; l\in[|\mathcal{S}_{t+1}|], \\
        & \bar{V}_{t+1, i}^R(x_t, \xi_{t+1, j}) - \eta_{t+1, k} \leq \Delta_{t+1, k, j}, \quad \text{for}\; k, j \in[K_{t+1}], \\
        & \sum_{n\in \mathcal{N}_{t, i}} \theta_{n, j} \bar{V}_{t+1, i}^R(x_{t, n}, \xi_{t+1, j}) + M_{t+1} \|y_{t, j}\|_1 = \bar{V}_{t+1, i}^R(x_t, \xi_{t+1, j}), \nonumber \\
        & \hspace{20em} \quad \text{for}\; j \in [K_{t+1}],  \\
        & \sum_{n\in \mathcal{N}_{t, i+1}} \theta_{j, n} x_{t, n} + y_{t, j} = x_t, \quad \text{for}\; j \in [K_{t+1}],  \\ 
        & \sum_{n\in \mathcal{N}_{t, i+1}} \theta_{j, n} = 1, \quad \text{for}\; j \in [K_{t+1}], \label{}\\
        & \tilde{A}_t x_t=\tilde{b}_t-\tilde{E}_t x_{t-1},  \\
        & x_t \geq 0, \Delta_{t+1, k, j} \geq 0, \quad \text{for}\; k, j\in[K_{t+1}].
    \end{align}
\end{subequations}
%Different
Here constraints (\ref{eq-costToGoFunction-PDRO-tractable-upper}d)-(\ref{eq-costToGoFunction-PDRO-tractable-upper}f) 
correspond to (\ref{eq-upper-bound-t-obj})-(\ref{eq-upper-bound-t-b}).
It is important to note that differing from the non-robust case in (\ref{eq:change-of-measure-based subproblem}), 
here we construct multiple scenario-dependent cuts 
(\ref{eq-eq-costToGoFunction-PDRO-tractable-lower-c})
%in the above approximation problems, we use a multicut version, that is, for each realization of $\xi_t$, 
for the cost-to-go function $V_t^R(\cdot, \tilde{\xi}_{t})$ as opposed to a single/average of the scenario-based cuts in (\ref{eq:change-of-measure-based subproblem}).
This is because here we have exploited 
the optimization reformulation of CVaR (\ref{eq-CVaR-optimizationFormulation}) instead of dual formulation of CVaR (\ref{eq-CRM-dual}).
Likewise, the upper bound $\bar{V}_{t, i}^R(\cdot, \tilde{\xi}_{t})$ is calculated in 
%, which provide 
%lower approximated by $\underline{V}_{t, i}^R(\cdot, \tilde{\xi}_{t})$ 
a different way from (\ref{eq-upper-bound-t})
with scenario-based constraints (\ref{eq-costToGoFunction-PDRO-tractable-upper}d)-(\ref{eq-costToGoFunction-PDRO-tractable-upper}f).
The optimization problems at each iteration are tractable since the outer and upper approximations are linear programs, which can be solved efficiently %solvable usin
by off-the-shelf solvers. 
% In what follows we present the algorithm to solve problem (\ref{eq-DR-MARSRM}).

\subsubsection{Algorithms for Computing Lower and Upper Bounds}
{\color{black}
Summarizing the discussions above, we are able to propose an 
algorithm  which generates a sequence of 
lower bounds and another algorithm which generates a sequence of upper bounds. A combination of the two algorithms provides a sequence of intervals containing the true optimal value.
%We can will show that 
These intervals are guaranteed to 
 shrink to 
 % a singleton .  
%{\color{red}
the true optimal value with probability one.
%}

% for solving problem (\ref{eq-DR-MARSRM}).

\begin{algorithm}[H]
\caption{Construction of Lower Approximations of DR-MARSRM}
\label{alg-lowerbound-dro}
\SetAlgoLined

\KwIn {The maximum iteration $i^{\max}$, 
the initial cuts with $g_{t, j}^0$ and $G_{t, j}^0$ for $t=2,\ldots,T$ and $j=1,\ldots, K_t$.}

\KwOut {The optimal solutions $\hat{x}^{i}_{t,n}$ for $t\in[T]$, $n\in[N]$, $i\in[i^{\max}]$ and lower bound $\underline{V}^R_{1, i^{\max}}$.}
%{\color{red} 
%if available.
%HX: what do you mean?} 
\While{$i=1,\ldots,i^{\max}$} 
{
$\rhd$ \textbf{Forward pass:} obtain decision $\hat{x}^{i}_{t}$, $t=1,\ldots,T$.

Solve the approximate stage $1$ problem 
$\underline{V}^R_{1, i} (x_0, \xi_1)$, obtain the optimal solution 
$\hat{x}^{i}_{1} =x_1^*$ and optimal value $\underline{V}^R_{1, i}$ as a lower bound.

Draw a sample of $\xi_t$, for $t=2,\ldots,T$. 

\For{stage $t=2,\ldots,T$}
{
Solve the approximate stage $t$ subproblem (\ref{eq-costToGoFunction-PDRO-tractable-lower}) with $\tilde{\xi}_t= \xi_{t}$,
% $\underline{V}_{t, i} (\hat{x}_{t-1,n}, \xi_{t, n})$ , 
obtain optimal solution $\hat{x}_{t}^i$. 
}
\vspace{0.3cm}
$\rhd$ \textbf{Backward pass:} construct new cuts and  update lower approximations at each stage. 

\For{stage $t=T,\ldots,2$}
{
    \For{scenario $j=1,\ldots,K_t$}
    {
    Solve the approximate stage $t$ subproblem  (\ref{eq-costToGoFunction-PDRO-tractable-lower}) with $\tilde{\xi}_t= \xi_{t, j}$, 
    construct a
    %hyperplane 
    cut for stage $t-1$
    %: compute 
    with $G^i_{t, j}$ and
    $g^i_{t, j}$ defined as in (\ref{eq-cut-DRO}).
    }
}
}
\end{algorithm}

\begin{algorithm}[H]
\caption{Computation of Upper Bounds of DR-MARSRM}
\label{alg-upperbound-dro}
\SetAlgoLined
\For{
$i = 1,\ldots,i^{\max}$, stage $t=T,\ldots,2$, scenario $j=1,\ldots,K_t$, $n=1,\ldots,N$}
{
Solve (\ref{eq-costToGoFunction-PDRO-tractable-upper}) and obtain the optimal value $\bar{V}^R_{t,i}(x_{t-1},\xi_t)$ with $x_{t-1}=\hat{x}^{i}_{t-1,n}$ and $\xi_t=\xi_{t,j}$.
}

Solve (\ref{eq-costToGoFunction-PDRO-tractable-upper}) for $t=1$ and obtain the optimal value as an upper bound $\bar{V}^R_{1,i}(x_0,\xi_1)$. 

% \KwOut {the optimal solutions $\hat{x}^{i}_{t,n}$ for $t\in[T]$, $n\in[N]$, $i\in[i^{\max}]$ and lower bound $\underline{V}_{1, i^{\max}}$.
% \end{algorithm}
\end{algorithm}
}

% \vspace{1em}
% \begin{breakablealgorithm}
% \caption{lower approximation}
% \label{alg:DRO}
% \noindent
% \textbf{Initialize:} Set the maximum iteration $it^{\max}$, 
% % $g_{t, k}^0$, $G_{t, k}^0$ is the zero vector for $t\in[T]/\{1\}$ and $k\in[K_t]$, 
% set the initial cut $g_{t, j}^0$ and $G_{t, j}^0$, $t\in[T]/\{1\}$ and $j\in[K_t]$ if available.
% \begin{algorithmic}[1]
%     \While{$i=1,\ldots,it^{\max}$} 
%     \Statex \hspace{1.1em} {$\rhd$ \textbf{Forward pass:} obtain decision $\hat{x}^{i}_{t}$, $t=1,\ldots,T$}
%     \State \multiline{Solve the approximate stage $1$ problem 
%     $\underline{V}_{1, i} (x_0, \xi_1)$, obtain the optimal solution 
%     $\hat{x}_1^*$ and optimal value $\underline{V}_{i}$ as the lower bound.}
%     \State sample $\tilde{\xi}_t$, for $t=2,\ldots,T$.
%     \For{stage $t=2,\ldots,T$}
%         \State solve the approximate stage $t$ subproblem $\underline{V}_{t, i} (\hat{x}_{t-1}, \tilde{\xi}_t)$.
%     \EndFor 
%     \Statex \hspace{1.1em} {$\rhd$ \textbf{Backward pass:} obtain cuts}
%     \For{stage $t=T,\ldots,2$, scenario $j=1,\ldots,K_t$}
%     \State \multiline{solve the approximate stage $t$ subproblem 
%         $\underline{V}_{t, i} (\hat{x}_{t-1}, \xi_{t, j})$, create hyperplane cut for stage $t-1$: compute $G_{t, j}$ and
%         $\psi_{t, j}$ as in (\ref{eq-cut-DRO}) }
%     \EndFor
% \EndWhile \\
% \Return the optimal solutions $\hat{x}^{i}_{t}$ for $t\in[T]$, $it\in[it^{\max}]$.
% \end{algorithmic}
% \end{breakablealgorithm}

%=======================================================================
\subsection{Convergence of Algorithm~\ref{alg-lowerbound-dro}-\ref{alg-upperbound-dro}}

% {\color{red}HX: I think we should include an algorithm, otherwise it is not sensible to talk about convergence }
%{\color{red}
In this section, we discuss convergence of 
%the SDDP
Algorithm~\ref{alg-lowerbound-dro}-\ref{alg-upperbound-dro}.
%}
% To derive the convergence of
% the sequence of lower/upper bounds 
% obtained from Algorithms \ref{alg-lowerbound-dro} and \ref{alg-upperbound-dro}, 
% %the DPRO-ARSRM problem under the computational scheme presented in the previous sections, 
% we make the following assumption.
% % to derive the convergence.
% \begin{assumption}
% \label{asmp-basicsolution}
%     The solution $\hat{x}_{t, i}$ obtained in forward process
%     and the dual solutions $\varrho_{t, k, i}$ obtained in the backward process 
%     % {\color{red} 
%     of Algorithm \ref{alg-lowerbound-dro}
%     % }
%     % {\color{red} 
%     are basic feasible solutions for $t\in[T]$  in  each iteration $i\in[i^{\max}]$.
%     % HX: What do you really mean and how to justify the assumption?
%     % }.
% \end{assumption}
% %{\color{red}
% The assumption is 
% %used 
% needed to ensure 
% %the finiteness 
% well-definedness 
% of 
% %different 
% cutting planes 
% in the backward process.
% This assumption 
% %which is mild since linear
% %programs 
% can be easily satisfied 
% %by using the 
% when simplex method is used to solve linear programs.
% {\color{red} This raises further questions as to which linear programs 
% are concerned and whether they are solved by simplex methods. I suggest we remove the assumption is it is not essential.
% }
%}
Note that %we have not presented 
the convergence analysis to be presented in this section does not automatically 
cover 
%for 
the convergence of Algorithm~\ref{alg-lowerbound}-\ref{alg-upperbound} because the ambiguity set ${\cal Q}_t$
does not reduce to a singleton (the true distribution of $s_t$) even when 
the mean value $\mu_t$ and the 
covariance 
$\Sigma_t$ are true. 
However, we can demonstrate 
%the
convergence of Algorithm~\ref{alg-lowerbound}-\ref{alg-upperbound}
in a similar way.
% Though the DRO problem (\ref{eq-costToGoFunction-PDRO}) cannot reduce to the discretized problem (\ref{eq-costToGoFunction-ARSRM}) if we set the stagewise mean and covariance as the true counterparts {\color{red} HX: this sentence is incomplete}. 
% The convergence proof presented in this section can be similarly derived for MARSRM problem {\color{red} HX: why do we have to say this?}.
We first show that the cost-to-go function $\underline{V}_t^R (\cdot, \xi_t)$ is piecewise linear convex 
%in
with respect to the state variable $x_t$.

\begin{lemma}
\label{lem-piecewiselinear}
    %{\color{red}
Let  $V_t^R(x_{t-1},\tilde{\xi_{t}})$
be defined as in (\ref{eq-costToGoFunction-PDRO}), and assume that 
%      Assume that (\ref{eq-costToGoFunction-PDRO}) has finite optimal values 
is finite valued
for all realizations of $\xi_t$ and
    %for 
    $t\in [T]$.
    %}
    % At each stage $t\in[T]$, 
    Then
%    the cost-to-go function
$\underline{V}_t^R (\cdot, \xi_t)$ is a piecewise linear convex function on $\mathcal{X}_{t-1}$ with a finite number of pieces for any $\xi_t\in \Xi_t$ and feasible $x_{t-1}$.
\end{lemma}

\proof
We 
%start
prove by induction backwards 
%for 
from stage $t:=T$ 
%and 
with $\xi_{T, j}$ for $j\in[K_T]$. By definition
\begin{align*}
    \underline{V}_T^R (x_{T-1}, \xi_{T, j}) & = \min_{x_T} \left\{c_{T, j} x_T: A_{T, j} x_T = b_{T, j}- E_{T, j} x_{T-1}; x_T\geq 0\right\} \\
    & = \max_{\tau_{T, j}} \left\{\varrho_T(-E_{T, j} x_{T-1}+ b_{T, j}) : \varrho_T A_{T, j} \leq c_{T, j}\right\} \\
    & = \max_{\kappa\in[|\mathcal{K}_T|]}
    \left\{\varrho_{T, j, \kappa}(-E_{T, j} x_{T-1}+ b_{T, j})\right\},
\end{align*}
where $\varrho_{T,\kappa}$, $\kappa\in \mathcal{K}_{T, j}$ denotes 
the 
% {\color{red} 
vertices of polyhedron 
$%\varrho_T(-E_{T, j} x_{T-1}+ b_{T, j}) :
\{\varrho_T: \varrho_T A_{T, j} \leq c_{T, j}\}$
% }
%extreme points of the dual problem of stage $T$ 
and $\mathcal{K}_{T, j}$ is a finite set. 
The second equality follows from strong duality and the assumption of relatively complete recourse.
This shows that 
$ \underline{V}_T^R (x_{T-1}, \xi_{T, j})$ is piecewise linear and convex with a finite number of pieces, and hence the conclusion
%that the result holds 
for $t=T$. 
% {\color{red}
% The finite number of pieces is
% %guaranteed by
% due to the fact that
% %if we compute basic feasible solutions of
% the number of basic solutions of the 
%  dual problem
%  %, the number of such basic optimal solutions 
%  is finite, and hence the number of cutting planes is finite, see e.g.~\cite{ruszczynski2003decomposition,Alex11}.
% HX: I feel this explanation is irrelevant in this context. I will delete it. PLease let me know if I misunderstand. 
% }
Next, suppose that $\underline{V}_{t+1}^R (\cdot, \tilde{\xi}_{t+1})$ is convex and piecewise linear with finite pieces, i.e.,
%Then we have 
$\underline{V}_{t+1}^R (x_t, \xi_{t+1, j}) =\max_{\kappa\in \mathcal{K}_{t+1, j}} \{g_{t+1, j, \kappa} + G_{t+1, j, \kappa} x_t\}$ for some 
index set $\mathcal{K}_{t+1, j}$. Then
%and problem (\ref{eq-costToGoFunction-PDRO-tractable}) can be reformulated as
\begin{subequations}
    \label{eq-costToGoFunction-PDRO-tractable-analysis}
    \begin{align}
        \underline{V}_t^R(x_{t-1}, \tilde{\xi}_{t}) & \nonumber \\
        = \inf_{\substack{x_t, \bm{\zeta}_{t+1}, \\ \bm{\eta}_{t+1}, \bm{\Delta}_{t+1}}} \; & c_t x_t +\zeta_{t+1}^1 + \mu_{t+1} \zeta_{t+1}^2 + \zeta_{t+1}^3\cdot \Sigma_{t+1} \\
        \st \quad\; & \sum_{k\in[K_{t+1}]} \beta_{t+1, l, k}
        \left(\eta_{t+1, k}+\frac{1}{1-\alpha_{t+1, k}} \sum_{j\in[K_{t+1}]} p_{t+1, j} \Delta_{t+1, k, j} \right) -\zeta_{t+1}^1- s_{t+1, l}\zeta_{t+1}^2 \nonumber \\
        & \hspace{2em} - \zeta_{t+1}^3 \cdot (s_{t+1, l}-\mu_{t+1})(s_{t+1, l}-\mu_{t+1})^{\top} \leq 0, \quad \text{for}\; l\in[|\mathcal{S}_{t+1}|], \\
        & g_{t+1, j, \kappa} +G_{t+1, j, \kappa} x_t - \eta_{t+1, k} \leq \Delta_{t+1, k, j}, \quad \text{for}\; \kappa\in \mathcal{K}_{t+1,j}, k, j \in[K_{t+1}], \label{constraint-plus1-ms} \\
        & \tilde{A}_t x_t=\tilde{b}_t-\tilde{E}_t x_{t-1}, \\
        & x_t \geq 0, \Delta_{t+1, k, j} \geq 0, \quad  \text{for}\; k, j\in[K_{t+1}]. 
    \end{align}
\end{subequations}
%whose 
The Lagrange dual of the problem above  can be written as 
%derived as 
\begin{align*}
    \max_{\substack{\bm{\varrho}_t, \bm{p}_{t+1} \\ \bm{\tau}_{t+1} }} \;\; & \varrho_t (\tilde{b}_t- \tilde{E}_t x_{t-1})+ \sum_{k, j\in[K_{t+1}]} \sum_{\kappa\in\mathcal{K}_{t+1, j}} \tau_{t+1, k, j, \kappa} g_{t+1, j, \kappa} \\
    \st \quad \; & \sum_{l\in[|\mathcal{S}|]} p_{t+1, l} =1, \\
    & \sum_{l\in[|\mathcal{S}_{t+1}|]} p_{t+1, l} s_{t+1, l} =\mu_{t+1}, \\
    & \sum_{l\in[|\mathcal{S}_{t+1}|]} p_{t+1, l} (s_{t+1, l}-\mu_{t+1})(s_{t+1, l}-\mu_{t+1})^{\top}=\Sigma_{t+1},  \\
    & \sum_{l\in[|\mathcal{S}_{t+1}|]} p_{t+1, l} \beta_{t+1, l, k} - \sum_{j\in[K_{t+1}]} \sum_{\kappa\in\mathcal{K}_{t+1, j}} \tau_{t+1, k, j, \kappa} =0, \quad \text{for}\; k\in[K_{t+1}],  \\
    & \sum_{l\in[|\mathcal{S}_{t+1}|]} p_{t+1, l} \frac{\beta_{t+1, l, k}}{1-\alpha_{t+1,k}} - \sum_{\kappa\in\mathcal{K}_{t+1, j}} \tau_{t+1, k, j, \kappa} \geq 0, \quad \text{for}\; k,j \in[K_{t+1}], \\ 
    & c_t+\varrho_t\tilde{A}_t +\sum_{k, j\in[K_{t+1}]} \sum_{\kappa\in\mathcal{K}_{t+1, j}} \tau_{t+1, k,j,\kappa}G_{t+1, j,\kappa} \geq 0.
\end{align*}
The assumption on the finite-valuedness of $\underline{V}_{t+1}^R (x_t, \xi_{t+1, j})$
ensures the strong duality holds. 
The maximum of the above dual problem 
%can be achieved 
can be attained at finite set of 
extreme points,
%i.e.
denoted by $\mathcal{K}_t$. 
Thus
%then 
\begin{equation}
\label{eq-obj-dual}
    \underline{V}_t^R (x_{t-1},\tilde{\xi}_{t}) = \max_{l\in\mathcal{K}_t} \; \varrho_{t, l}(\tilde{b}_t-\tilde{E}_t x_{t-1} )
    + \sum_{k, j\in[K_{t+1}]} \sum_{\kappa\in\mathcal{K}_{t+1, k}} \tau_{t+1, k, j, \kappa, l} g_{t+1, j, \kappa},
\end{equation}
which implies that $\underline{V}_t^R (x_{t-1}, \tilde{\xi}_{t})$ is 
%a 
piecewise linear 
%function 
and convex with respect to $x_{t-1}$ with finite number of pieces. 
\hfill \Box

The piecewise linear convexity of the cost-to-go function 
%implies that
$\underline{V}_t^R (x_{t-1}, \xi_t)$ for all $\xi$ implies that it 
can be restored by a set of support
%ing 
%hyperplanes
functions, 
although the number of
%hyperplanes 
such planes could be prohibitively large. 
Next,
we establish that the 
%outer
lower and upper approximations of the function
may provide
valid lower and upper bounds for
the optimal value of problem (\ref{eq-costToGoFunction-PDRO}).
%{eq-mathcalQ_t}).
% of {\color{red} ???? please add}
%the multistage problem respectively.

\begin{theorem}
\label{lem:convg-low-up-bnds}
Let  $L_{t}^R, t\in [T]$ be 
the Lipschitz modulus of $V^R_{t}(x_{t-1}, \xi_{t})$ under the 1-norm.
Suppose that there exist positive constants  $\{M_t\}_{t\in [T]}$ 
%be a constant
%positive constants.
such that 
    % For the upper approximation (\ref{eq-costToGoFunction-PDRO-tractable-upper}), suppose 
    $L_{t}^R \leq M_{t}$ for all $t\in[T]$.
    %, where $L_{t}^R$ is the Lipschitz constant for $V^R_{t}(x_{t-1}, \xi_{t})$ under the 1-norm. 
     % For the upper approximation (\ref{eq-costToGoFunction-PDRO-tractable-upper}),
    % For any iteration $i$, t
    Then the upper and lower approximations 
    derived in (\ref{eq-costToGoFunction-PDRO-tractable-lower}) and (\ref{eq-costToGoFunction-PDRO-tractable-upper})
    provide valid bounds for the optimal value of the subproblem (\ref{eq-costToGoFunction-PDRO})
    %, that is, 
    in the sense that
    \begin{equation}
    \label{eq-boundproof}
        \underline{V}^R_{t+1, i}(\cdot, \xi_{t+1}) \leq V^R_{t+1}(x_t, \xi_{t+1}) \leq \bar{V}^R_{t+1, i}(x_t, \xi_{t+1}), \quad %\text{for}\; t\in[T-1], \; 
        \forall \xi_{t+1}\in\Xi_{t+1}
    \end{equation}
    for $t\in[T-1]$ and 
     iterations $i=1,2,3,\cdots$.
\end{theorem}

\proof
We prove (\ref{eq-boundproof}) by induction. 
We begin 
%the proof by 
by considering the backward process. Recall the lower bound problem at the $\tilde{i}$-th iteration:
\begin{subequations}
    \label{eq-costToGoFunction-PDRO-tractable-lower-n}
    \begin{align}
        \underline{V}_{t, \tilde{i}}^R(x_{t-1}, \tilde{\xi}_{t}) & \nonumber \\
        = \inf_{\substack{x_t, \bm{\zeta}_{t+1}, \\ \bm{\eta}_{t+1}, \bm{\Delta}_{t+1}}} \; & c_t x_t +\zeta_{t+1}^1 + \mu_{t+1} \zeta_{t+1}^2 + \zeta_{t+1}^3\cdot \Sigma_{t+1} \\
        \st \quad\; & \sum_{k\in[K_{t+1}]} \beta_{t+1, l, k}
        \left(\eta_{t+1, k}+\frac{1}{1-\alpha_{t+1, k}} \sum_{j\in[K_{t+1}]} p_{t+1, j} \Delta_{t+1, k, j} \right) -\zeta_{t+1}^1- s_{t+1,l}\zeta_{t+1}^2 \nonumber \\
        & \hspace{2em} - \zeta_{t+1}^3 \cdot (s_{t+1,l}-\mu_{t+1})(s_{t+1,l}-\mu_{t+1})^{\top} \leq 0, \quad \text{for}\; l\in[|\mathcal{S}_{t+1}|], \\
        & g_{t+1, j, n} +G_{t+1, j, n} x_t - \eta_{t+1, k} \leq \Delta_{t+1, k, j}, \quad \text{for}\; n\in \mathcal{N}_{t+1, \tilde{i}}, k, j \in[K_{t+1}], \label{constr-cut} \\
        & \tilde{A}_t x_t=\tilde{b}_t-\tilde{E}_t x_{t-1}, \label{constr-sub} \\
        & x_t \geq 0, \Delta_{t+1, k, j} \geq 0, \quad  \text{for}\; k, j\in[K_{t+1}], 
    \end{align}
\end{subequations}
where $\mathcal{N}_{t+1, \tilde{i}}$ is the index set representing the total previous cuts in the first $\tilde{i}$ for $V_t^R(x_{t-1}, \tilde{\xi}_{t})$. 
\underline{At stage $T$}, for any feasible solution $x_{T-1}$ and $\xi_T\in\Xi_t$, let $\varrho_T$ be the dual variable associated with constraint (\ref{constr-sub}). Then 
\begin{equation*}
    \begin{split}
        g_{T, j} +G_{T, j} x_{T-1}
        & = \varrho_T (\tilde{b}_T-\tilde{E}_T x_{T-1}) \\
        & \leq \max\{\varrho_T (\tilde{b}_T-\tilde{E}_T x_{T-1}), \tilde{A}_T\varrho_T \leq \tilde{c}_T, \varrho_T\in\R^{d_T} \} \\
        & = \min \{\tilde{c}_T x_T: \tilde{A}_Tx_T+ \tilde{E}_T x_{T-1} = \tilde{b}_T, x_T\in\R^{d_T}_+ \} \\
        & = V_T^R (x_{T-1}, \tilde{\xi}_T),
    \end{split}
\end{equation*}
which means constraint 
(\ref{eq-costToGoFunction-PDRO-tractable}c) is tighter than constraint (\ref{eq-costToGoFunction-PDRO-tractable-lower-n}c)
for the case that $t=T-1$. Since all other constraints are the same, then we
assert that the optimal values of the two programs satisfy
$\underline{V}_{T-1, i}^R(x_{T-2}, \tilde{\xi}_{T-1}) \leq V_{T-1}^R(x_{T-2}, \tilde{\xi}_{T-1})$.
For the upper bound, 
%first 
we have $\bar{V}_{T, i}^R(x_{T-1}, \tilde{\xi}_{T}) = V_{T}^R(x_{T-1}, \tilde{\xi}_{T})$ at stage $T$ because %{\color{red}Since there is no cost-to-go function at stage $T$??
$V_{T+1}^R(x_{T}, \tilde{\xi}_{T+1})\equiv 0$.
%}
By the convexity of $V_{T}^R(\cdot, \tilde{\xi}_{T})$, we have 
\begin{equation}
\label{eq-convexity-V}
\resizebox{0.9\linewidth}{!}{$
    V_{T}^R(x_{T-1}, \xi_{T, j}) - \varpi_T y_{T-1, j} 
    \leq V_{T}^R \left(\sum_{n\in\mathcal{N}_{t+1, i}} \theta_{n, j} \hat{x}_{T-1, n}, \xi_{T, j}\right) \leq \sum_{n\in\mathcal{N}_{t+1, i}} \theta_{n, j} V_{T}^R( \hat{x}_{T-1, n}, \xi_{T, j}), 
$}
\end{equation}
where $\varpi_T$ is a subgradient of $V_{T}^R(\cdot, \xi_{T, j})$ at $x_{T-1}$. 
Therefore, 
\begin{equation*}
\begin{split}
     \sum_{n\in\mathcal{N}_{t+1, i}} \theta_{n, j} \bar{V}_{T, i}^R( \hat{x}_{T-1, n}, \xi_{T, j}) + M_T \|y_{T-1, j}\|_1 & = \sum_{n\in\mathcal{N}_{t+1, i}} \theta_{n, j} V_{T}^R( \hat{x}_{T-1, n}, \xi_{T, j}) + M_T \|y_{T-1, j}\|_1 \\
    & \geq \sum_{n\in\mathcal{N}_{t+1, i}} \theta_{n, j} V_{T}^R( \hat{x}_{T-1, n}, \xi_{T,j}) + \|\varpi_T \|_1 \|y_{T-1, j}\|_1 \\
    & \geq \sum_{n\in\mathcal{N}_{t+1, i}} \theta_{n, j} V_{T}^R( \hat{x}_{T-1, n}, \xi_{T,j}) + \|\varpi_T \|_1 \|y_{T-1, j}\|_2 \\
    & \geq V_{T}^R \left( \sum_{n\in\mathcal{N}_{t+1, i}} \theta_{n, j} \hat{x}_{T-1, n}, \xi_{T, j} \right) + \varpi_T y_{T-1, j} \\
    & \geq V_{T}^R( \hat{x}_{T-1}, \xi_{T, j}),
\end{split}
\end{equation*}
where $\hat{x}_{T-1}=\sum_{n\in\mathcal{N}_{t+1, i}} \theta_{n, j} \hat{x}_{T-1, n}$.
The first inequality holds because $M_T \geq L_T^R$ and the subgradient of $V_{T}^R( \hat{x}_{T-1}, \xi_{T, j})$ under $1$-norm is bounded by $M_T$. The convexity and Cauchy-Schwartz inequality imply the third inequality, where the last inequality follows from (\ref{eq-convexity-V}).
The inequalities above 
%This shows 
imply that constraint (\ref{eq-costToGoFunction-PDRO-tractable-upper}) is tighter than (\ref{eq-costToGoFunction-PDRO-tractable}d), and hence the optimal values of the two programs satisfy
$V_{T-1, i}^R(x_{T-2}, \tilde{\xi}_{T-1}) \leq \bar{V}_{T-1}^R(x_{T-2}, \tilde{\xi}_{T-1})$, this in turn show that 
$$
\underline{V}_{T-1, i}^R(x_{T-2}, \tilde{\xi}_{T-1}) \leq V_{T-1}^R(x_{T-2}, \tilde{\xi}_{T-1}) \leq \bar{V}_{T-1, i}^R(x_{T-2}, \tilde{\xi}_{T-1}), 
$$
that is, the claim holds for stage $T-1$.

Suppose now that (\ref{eq-boundproof}) holds for fixed $t':=t+1\leq T-1$. We show that the inequalities hold for $t':=t$.
%, i.e., 
% $$
% \underline{V}_{t+1, i}^R(x_{t}, \tilde{\xi}_{t+1}) \leq V_{t+1}^R(x_{t}, \tilde{\xi}_{t+1}) \leq \bar{V}_{t+1, i}^R(x_{t}, \tilde{\xi}_{t+1}). 
% $$
The first inequality in (\ref{eq-boundproof}) 
implies that the cuts generated for the lower approximations are valid for $t':=t+1$, i.e.,
\begin{equation*}
    \max\big\{g_{t, j, n} +G_{t, j, n} x_t, n\in \mathcal{N}_{t+1, i}\big\} 
    \leq V_{t+1}^R(x_{t}, \xi_{t+1, j}).
\end{equation*}
Let $\mathcal{D}_t$ and $\mathcal{D}_{t, i}$ denote respectively the feasible regions of problems (\ref{eq-costToGoFunction-PDRO}) and (\ref{eq-costToGoFunction-PDRO-tractable-upper}) defining $V_{t+1}^R(x_{t}, \xi_{t+1, j})$ and $\underline{V}_{t+1, i}^R(x_{t}, \xi_{t+1, j})$, 
let
$(\varrho_t, \tau_{t+1})$ be 
the Lagrange multipliers of
%a dual feasible extreme point of 
problem (\ref{eq-costToGoFunction-PDRO-tractable-lower-n}) associated with constraints (\ref{constr-sub}) and (\ref{constr-cut}). 
Using similar result to that in (\ref{eq-obj-dual}), we obtain
\begin{equation*}
\begin{split}
    g_{t-1} +G_{t-1} x_t 
    & = \varrho_{t}(\tilde{b}_t-\tilde{E}_t x_{t-1} )
    + \sum_{k, j\in[K_{t+1}]} \sum_{\kappa\in \mathcal{N}_{t+1, i}} \tau_{t+1, k, j, \kappa} g_{t+1, j, \kappa} \\
    & \leq \max \left\{\varrho_{t}(\tilde{b}_t-\tilde{E}_t x_{t-1} )
    + \sum_{k, j\in[K_{t+1}]} \sum_{\kappa\in \mathcal{N}_{t+1, i}} \tau_{t+1, k, j, \kappa} g_{t+1, j, \kappa}: (\varrho_t, \tau_{t+1}) \in \mathcal{D}_{t, i} \right\} \\
    & \leq \max \left\{\varrho_{t}(\tilde{b}_t-\tilde{E}_t x_{t-1} )
    + \sum_{k, j\in[K_{t+1}]} \sum_{\kappa\in \mathcal{N}_{t+1, i}} \tau_{t+1, k, j, \kappa} g_{t+1, j, \kappa}: (\varrho_t, \tau_{t+1}) \in \mathcal{D}_{t} \right\} \\
    & \leq V_t^R(x_{t-1}, \tilde{\xi}_{t}).
\end{split}
\end{equation*}
The first inequality holds since $(\varrho_t, \tau_{t+1})\in\mathcal{D}_{t, i}$, the second inequality holds since $\mathcal{D}_{t, i}\subset \mathcal{D}_{t}$, and the last inequality is due to strong duality. 
Thus $\underline{V}_{t, i}^R(x_{t-1}, \tilde{\xi}_{t})\leq V_t^R(x_{t-1}, \tilde{\xi}_{t})$.

For the upper approximation, by hypothesis we have $V_{t}^R(x_{t-1}, \tilde{\xi}_{t})\leq \bar{V}_{t, i}^R(x_{t-1}, \tilde{\xi}_{t})$.
Replace $\bar{V}_{t, i}^R(x_{t-1}, \tilde{\xi}_{t})$ with $V_{t}^R(x_{t-1}, \tilde{\xi}_{t})$ in the constraints of problem (\ref{eq-costToGoFunction-PDRO-tractable-upper}), we obtain a new problem $\bar{V}_{t, i}^{\rm R_L}(x_{t-1}, \tilde{\xi}_{t})$, we have $\bar{V}_{t, i}^{\rm R_L}(x_{t-1}, \tilde{\xi}_{t}) \leq \bar{V}_{t, i}^R(x_{t-1}, \tilde{\xi}_{t})$.
Next we prove $\bar{V}_{t, i}^{\rm R_L}(x_{t-1}, \tilde{\xi}_{t})$ is an upper bound for $V_{t}^R(x_{t-1}, \tilde{\xi}_{t})$.
Let $\varpi_t$ be a subgradient of $\bar{V}_{t+1}^R(\cdot, \xi_{t+1, j})$ at $x_t$. We have
\begin{equation*}
\begin{split}
    V_{t+1}^R(x_{t}, \xi_{t+1, j}) & \leq 
    \sum_{n\in\mathcal{N}_{t+1, i}} \theta_{n, j} V_{t+1}^R( \hat{x}_{t, n}, \xi_{t+1, j}) + \varpi_t y_{t,j} \\
    & \leq 
    \sum_{n\in\mathcal{N}_{t+1, i}} \theta_{n, j} V_{t+1}^R( \hat{x}_{t, n}, \xi_{t+1, j}) + \|\varpi_t\|_2 \|y_{t,j}\|_2 \\
    & \leq 
    \sum_{n\in\mathcal{N}_{t+1, i}} \theta_{n, j} V_{t+1}^R( \hat{x}_{t, n}, \xi_{t+1, j}) + \|\varpi_t\|_1 \|y_{t,j}\|_1 \\
    & \leq  
    \sum_{n\in\mathcal{N}_{t+1, i}} \theta_{n, j} V_{t+1}^R( \hat{x}_{t, n}, \xi_{t+1, j}) + M_t \|y_{t,j}\|_1 \\
    & = \bar{V}_{t+1}^R(x_{t}, \tilde{\xi}_{t+1}).
\end{split}
\end{equation*}
Consequently we have 
$$
V_{t}^R(x_{t-1}, \tilde{\xi}_{t}) \leq \bar{V}_{t, i}^{\rm R_L}(x_{t-1}, \tilde{\xi}_{t}) \leq \bar{V}_{t, i}^R(x_{t-1}, \tilde{\xi}_{t})
$$
and 
$$
\underline{V}_{t, i}^R(x_{t-1}, \tilde{\xi}_{t}) \leq V_{t, i}^R(x_{t-1}, \tilde{\xi}_{t}) \leq \bar{V}_{t, i}^R(x_{t-1}, \tilde{\xi}_{t}).
$$
The proof is complete.
\hfill \Box

%In the following 
The theorem
%, we give the finite time convergence of 
states that Algorithm~\ref{alg-lowerbound-dro}-\ref{alg-upperbound-dro} terminate in a finite number of iterations for the prescribed tolerance.
%the SDDP algorithm for solving DR-MARSRM problem follows the proof of \cite{park2023data}, which is used in the distributionally robust risk neutral case.

\begin{theorem}

    The sequences of lower bounds and upper bounds generated by Algorithm~\ref{alg-lowerbound-dro}-\ref{alg-upperbound-dro} 
    %for solving 
    converge to the optimal value of problem (\ref{eq-DR-MARSRM}) %converges to an %approximate 
    %optimal solution 
    in a finite number of iterations with probability one. %{\color{red}HX: do we have convergence of optimal solutions?}
  %  {\color{blue} Q: We do not have.}
\end{theorem}

% {\color{red}Where is Lemma~\ref{lem:convg-low-up-bnds} used in the proof?}

\proof
We first discuss the convergence of the lower bounds.
%approximation. 
By Lemma \ref{lem-piecewiselinear}, $V_t(\cdot, \tilde{\xi}_t)$ is a piecewise linear convex function with a finite number of pieces. 
Moreover, the number of scenarios at each stage is finite and there is a nonzero probability for any possible scenario to occur in the forward step since scenarios are generated by Monte Carlo sampling. 
% {\color{red} Thus, any possible scenario occurs infinitely many times in the forward step unless the algorithm terminates.
% HX: rephrase the sentence}
% {\color{blue}
If the algorithm does not terminate, then any possible scenario may occur infinitely many times in the forward step.
%}

As shown in \cite{Alex11}, the lower bound convergence holds if an optimal policy $\{\hat{x}_t (\tilde{\xi}_t)\}_{t=1}^T$ 
%for 
of 
% {\color{red}the current HX: which iteration?} 
% {\color{blue} This is a general optimality condition, one policy satisfies the dynamic optimality.}
the lower approximation problem satisfies the following dynamic 
%programming
optimality condition
\begin{equation}
\label{eq-dunamic-opt}
\begin{split}
    \left\{\hat{x}_t(\tilde{\xi}_t)\right\}_{t=1}^T \in & \arg\min_{x_t \in \mathcal{X}_t (\hat{x}_{t-1} (\tilde{\xi}_{t-1}), \tilde{\xi}_t)}  
    \Bigg\{c_t x_t \\
    & + \max_{\bm{q}_{t+1}\in \mathcal{Q}_{t+1}} \sum_{i\in[|\mathcal{S}_{t+1}|]} \sum_{k\in[K_{t+1}]} q_{t+1,i}\beta_{t+1,i,k} \CVaR_{\alpha_k}(V_{t+1}^R(x_{t},\xi_{t+1})|\xi_{t}=\tilde{\xi_{t}}
    ) \Bigg\}. 
\end{split}
\end{equation}
Let $\{\hat{x}_{t, i} (\tilde{\xi}_t)\}_{t=1}^T$ be a policy obtained by lower approximation $\underline{V}_{t, i}^R(\cdot, \tilde{\xi}_t)$ for all $t\in [T]\setminus\{1\}$ 
%{\color{blue} 
in the $i$th iteration.
%{\color{red}HX: is the definition a repeat?}
Suppose that (\ref{eq-dunamic-opt}) does  not hold for some stage $t\in [T]\setminus\{1\}$ and possible scenario, i.e., a policy $\{\hat{x}_{t, i} (\tilde{\xi}_t)\}_{t=1}^T$ for the current lower approximation is not optimal. 
Let stage $\tilde{t}$ be the largest stage $t$ such that $\hat{x}_{\tilde{t}, i} (\tilde{\xi}_{\tilde{t}})$ does not satisfy (\ref{eq-dunamic-opt}), i.e., a candidate solution $\hat{x}_{\tilde{t}, i}=\tilde{x}_{\tilde{t}, i} (\xi_{\tilde{t}})$ is suboptimal for the current scenario. 
Then, at some iteration $\tilde{i}>i$, we add a new cut corresponding to $\hat{x}_{\tilde{t}, i}$, updating the current approximation $\underline{V}^R_{\tilde{t}, i}$. 
Similarly, a new cut is added until (\ref{eq-dunamic-opt}) holds for stage $\tilde{t}$. 
Such cut generations continue for some stage $t<\tilde{t}$ when (\ref{eq-dunamic-opt}) does not hold. 
After a sufficiently large number of iterations, (\ref{eq-dunamic-opt}) holds for all stages and possible
scenarios. 
This completes the proof for the lower bound convergence.
% {\color{red}HX: I would expect a contradiction to the assumption ``Suppose that (\ref{eq-dunamic-opt}) does  not hold for some stage $t\in [T]\setminus\{1\}$''}

Next, we discuss convergence of 
the upper bounds.
%convergence.
Let $\bar{i}$ be the iteration such that the lower bound convergence holds at any iteration $i\geq \bar{i}$.
% {\color{red} 
%By the above 
Based on the discussions above, we conclude that 
the sequence of the lower bounds converges
% convergence under the Assumption \ref{asmp-basicsolution}
% implies that, 
after iteration $\bar{i}$, 
%then 
and hence
there 
%only 
exists only 
a finite number of optimal policies $\{\hat{x}_t(\xi_t)\}_{t=1}^T$ after 
%the convergence. 
the $\bar{i}$th iteration.
In other words, for each scenario, we have only a finite number of sequences of candidate solutions $\{\hat{x}_t\}_{t=1}^T = \{\hat{x}_t(\xi_t)\}_{t=1}^T$ at which the upper bound approximation is evaluated (recall that candidate solutions are obtained by solving the lower bound problem in the algorithm). 
Hence, we will show that the gap between the lower and upper approximation is         zero at those finite number of candidate solutions.

Without loss of generality, let us assume that there is only one optimal policy. Let $\{\hat{x}_{t, j_t}\}_{t=1}^T$ be the sequence of the optimal
%sequence of 
solutions for scenario $j_t\in[K_t]$ after the lower bound
convergence. Here, we show 
\begin{equation}
\label{eq-lower-upper}
    \bar{V}^R_t(\hat{x}_{t-1, j_{t-1}}, \xi_{t, j}) = \underline{V}^R_t(\hat{x}_{t-1, j_{t-1}}, \xi_{t, j}) \quad \text{for}\; j\in[K_{t}], \text{ for some iteration } i\geq \bar{i}
\end{equation}
for all stages.

Let $\bar{i}_t$ be the iteration for which (\ref{eq-lower-upper}) holds at stage $t$. 
At the stage $T$, (\ref{eq-lower-upper}) holds after the
upper approximation at stage $T$ is evaluated at all possible candidate solutions $\bar{x}_{T, j}, j\in[K_T]$, since $\bar{V}^R_{T+1, i} (\cdot, \cdot) = \underline{V}^R_{T+1, i} (\cdot, \cdot)=0$.
Proceeding by induction, for $t=T,\ldots,2$, there exists iteration $i\geq \bar{i}_{t+1}$ such that (\ref{eq-lower-upper}) holds at stage $t$ because $\bar{V}^R_{t+1, i} (\hat{x}_{t, j_t}, \xi_{t+1, j}) = \underline{V}^R_{t+1, i} (\hat{x}_{t, j_t}, \xi_{t+1, j})$ for all $j, j_t\in[K_{t+1}]$. 
This completes the convergence of the upper bound.
\hfill \Box

% {\color{red}
%========================================================================================

\subsection{Step-like Approximation of Risk Spectra and Its Impact on Optimal Values}
%Error bounds on the optimal value of MARSRM}

The computational schemes outlined in the preceding
sections rely heavily on step-like risk spectrum $\sigma$. In practice, the risk spectrum of a DM may not be step-like. Consequently, there is a need to 
estimate the effect of the modelling error incurred by the step-like approximation.
In this section, we give a theoretical guarantee that such error is containable.
% In the previous sections, we discuss the multistage problem with step-like risk spectrum $\sigma$.
% In this section, we discuss the errors arising from the step-like approximation and its propagation to the optimal value in the multistage problem.
Let $\tilde{\sigma}$ and $\tilde{\rho}_Q$ denote the general risk spectrum and ARSRM with step-like random risk spectrum as defined in (\ref{eq-ARSRM-step}) respectively, $\tilde{V}_t$ and $V_t$ denote the optimal value of the corresponding subproblems without and with the step-like approximation, that is, $V_t(x_{t-1}, \xi_t)$ be defined as in (\ref{eq-costToGoFunction-ARSRM}) and 
$$
\tilde{V}_t(x_{t-1}, \xi_t):=\min_{x_t\in\mathcal{X}_t(x_{t-1},\xi_t)} \; c_tx_t+
\tilde{\rho}_{Q_{t+1}}( \tilde{V}_{t+1}(x_{t},\xi_{t+1})). 
$$
In order to estimate the error arise from the step-like approximation, we first define the distance in the space of risk spectrum.
Let $\mathscr{G}$ denote the set of all risk spectra $\sigma\in\mathcal{L}_1[0,1]$ such that $\sigma$ is nonnegative, nondecreasing and normalized, and $\mathscr{G}_J$ denote the set of all nonnegative, nondecreasing and normalized step-like risk spectra $\sigma_J$ over $[0,1]$ with breakpoints $0=z_0<\cdots<z_{J+1}=1$.
By definition, $\mathscr{G}_J\subset \mathscr{G}$.
We call $\sigma_J\in\mathscr{G}_J$ a projection of $\sigma$ onto the space $\mathscr{G}_J$ if $\sigma_J(z):=\sigma_j$ for $z\in[z_j, z_{j+1})$ such that $\sigma\in[\sigma(z_j), \sigma(z_{j+1})]$, $j=0,\ldots,J$ and $\int_0^1 \sigma_J(z) d z=\sum_{j=0}^J \sigma_j(z_{j+1}, z_j)=1$.
Let $\Phi:[0,1]\to\R_+$ be a positive function with $\int_0^1 \Phi(z) d z<+\infty$ and $\mathcal{G}$ be a set of measurable functions $\varphi: [0,1] \to \R$ such that $|\varphi(z)|\leq\Phi(z)$ for $t\in [0,1]$.
For any $\sigma_1, \sigma_2\in\mathscr{G}$, let
\begin{equation}
\label{eq:psi-distance}
    \dd_{\Phi} (\sigma_1, \sigma_2):= \sup_{\varphi\in\mathcal{G}} |\langle \varphi, \sigma_1 \rangle - \langle \varphi, \sigma_2 \rangle|,
\end{equation}
where $\langle \varphi, \sigma \rangle=\int_0^1 \varphi(z) \sigma(z) d z$.
$\dd_{\psi} (u,v)$ is a pseudo-metric in the sense that $\dd_{\psi} (u,v)=0$ implies
$\langle g, u\rangle - \langle g, v\rangle=0$ for all $g\in \mathscr{G}$ but not necessarily $u=v$ unless
the set $\mathscr{G}$ is sufficiently large. 
From practical point of view,
$\mathscr{G}$ may be regarded as a set of test functions such as the quantile functions of some prospects
and 
{\color{black}$\dd_{\psi} (u,v)=0$} 
%the equality 
means the investor does not tell the difference between $u$ and $v$ for this particular set of investment prospects.
The pseudo-metric provides a tighter measurement for discrepancy of two risk spectra than $\|\cdot\|_1$-norm, see \cite[Remark 3.1]{WaX20}.

% Since $u,v\in \mathscr{L}^1[0,1]$, (\ref{eq:psi-distance}) is satisfied
% when $\psi\in \mathscr{L}^\infty[0,1]$. 
%The function $\psi$ is not
%In that case, $\dd_{\psi} (u,v)$ is analogous to the total variation metric in probability theory \cite{gibbs2002choosing}.
% Our interest here is in the case that $\psi\in \mathscr{L}^p[0,1]$ for some $p> 1$, and we may choose $u,v\in \mathscr{L}^q[0,1]$ with $\frac{1}{p}+\frac{1}{q} =1$ to guarantee (\ref{eq:psi-distance}). 

% {\color{red} 
% The definition of the space of risk spectrum, projection, and Lipschitz.
% }

\begin{assumption}
\label{as-Phi} 
For $t\in[T]$, let $\mathscr{G}_t:=\left\{F^{\leftarrow}_{\tilde{V}_{t}(x_{t-1}, \xi_{t})}(\cdot): x_{t-1}\in\mathcal{X}_{t-1}\right\}$. 
    There exist 
    % a constant $p_t\in[1, +\infty)$ and 
    a positive function $\Phi_t$ with $\int_0^1 \Phi_t(z) d z<+\infty$ such that 
    %the following growth condition in the space of quantile functions holds:
    \begin{equation*}
        \left| F^{\leftarrow}_{\tilde{V}_{t}(x_{t-1}, \xi_{t})}(z) \right| \leq \Phi_t(z), \quad \text{for}\; z\in[0,1].
    \end{equation*}  
\end{assumption}
% {\color{red}HX: where is  $\Phi_t\in\mathcal{L}_{p_t}[0,1]$ used?
% }
%{\color{red}HX: Discuss how strong this condition is.}
%{\color{blue} 
To see how 
Assumption \ref{as-Phi} may be satisfied, we note that
%satisfied under some circumstances. For instance, if 
under Assumption \ref{assumption-multistageproblem} (c), we can show recursively that there exist 
$a_t(\xi_t)$ and $b_t(\xi_t)$ such that 
$$
a_t(\xi_t) \leq \tilde{V}_{t}(x_{t-1}, \xi_{t}) \leq b_t(\xi_t), \forall x_{t-1}\in\mathcal{X}_{t-1}, \xi_t\in\Xi_t, 
$$ 
for $t\in[T]$. 
% {\color{red} Since
% $\mathcal{X}_{t-1}(x_{t-2}, \xi_{t-1})$ 
% is a nonempty polyhedral set, 
% we may set: Clarify why}
%and letting 
Consequently we may set 
$\Phi_t(z):= \max\left\{F^{\leftarrow}_{a_t(\xi_t)}(z), F^{\leftarrow}_{b_t(\xi_t)}(z)\right\}$ for all $t\in[0.1]$, where $F^{\leftarrow}_{a_t(\xi_t)}(z)$
and 
$F^{\leftarrow}_{b_t(\xi_t)}(z)$ denote the 
quantile functions of $a_t(\xi_t)$ and $b_t(\xi_t)$ respectively. Assumption \ref{as-Phi} is satisfied  when 
the two quantile functions are integrable. 
Let $\tilde{\rho}_Q(\cdot)$ and $\rho_Q(\cdot)$ denote the ARSRM without and with step-like random risk spectrum defined as in (\ref{eq-costToGoFunction-PDRO}) and (\ref{eq-ARSRM-step}), and let
$$
\tilde{V}_t(x_{t-1}, \xi_t):=\min_{x_t\in\mathcal{X}_t(x_{t-1},\xi_t)} \; \left\{ c_tx_t+
\max_{Q_{t+1}\in\mathcal{Q}_{t+1}} \tilde{\rho}_{Q_{t+1}}( \tilde{V}_{t+1}(x_{t},\xi_{t+1}))\right\}.
$$
% where $\tilde{\rho}_{Q_{t+1}}$ is defined as in (\ref{eq-ARSRM}).
The next theorem quantifies 
%the difference between 
the error in terms of the optimal values when we use 
%the optimal value of 
problem (\ref{eq-costToGoFunction-PDRO}) as an approximation to 
%based on the two ARSRM measures.
the one with general risk spectrum.

% (\ref{eq-ARSRM-step}) respectively, $\tilde{V}_t$ and $V_t$ denote the optimal value of the corresponding subproblems, that is, $V_t(x_{t-1}, \xi_t)$ be defined as in (\ref{eq-costToGoFunction-ARSRM}) and 
% $$
% \tilde{V}_t(x_{t-1}, \xi_t):=\min_{x_t\in\mathcal{X}_t(x_{t-1},\xi_t)} \; c_tx_t+
% \tilde{\rho}_{Q_{t+1}}( \tilde{V}_{t+1}(x_{t},\xi_{t+1})). 
% $$
\begin{theorem}
\label{thm-error-step-DR}
% Let $V_t(x_{t-1}, \xi_t)$ be defined as in (\ref{eq-costToGoFunction-PDRO}) and 
% $$
% \tilde{V}_t(x_{t-1}, \xi_t):=\min_{x_t\in\mathcal{X}_t(x_{t-1},\xi_t)} \; c_tx_t+
% \max_{Q_{t+1}\in\mathcal{Q}_{t+1}} \tilde{\rho}_{Q_{t+1}}( \tilde{V}_{t+1}(x_{t},\xi_{t+1})), 
% $$
% where $\tilde{\rho}_{Q_{t+1}}$ is ARSRM defined as in (\ref{eq-ARSRM}) with general $\sigma$. 
%Let 
Let 
$\sigma_{s_\nu}\in\mathscr{G}$ be Lipschitz continuous with modulus being bounded by $L_\nu$ and $\tilde{\sigma}_{s_{t,\nu}}$ be a projection of $\sigma_{s_{t,\nu}}$ on $\mathscr{G}_J$ for $\nu\in|\mathcal{S}_t|$. Under Assumption \ref{as-Phi}, 
%hold for each $t\in[T]$.
%Then 
\begin{equation}
\label{eq-thm-errorbound}
    \left|V^R_t(x_{t-1}, \xi_t)-\tilde{V}^R_t(x_{t-1}, \xi_t)\right| \leq \sum_{t'\in[T]\setminus[t]} \max_{\nu\in[|\mathcal{S}_t|]} 
    L_\nu \beta_J\int_0^1 \Phi_{t'}(z) d z \quad \text{for}\; t\in[T],
\end{equation}
%{\color{blue} 
where $\beta_J:= \max_{j\in[J]} (z_{j+1}-z_j)$ and $z_j$ is the breakpoint of the approximate step function of $\sigma$ for $j\in[J]$.
%}

\end{theorem}

\proof
%{\color{blue}
We prove (\ref{eq-thm-errorbound}) by induction.
Suppose that (\ref{eq-thm-errorbound}) holds at stage $t+1$, i.e.,
%we have 
\begin{equation*}
    \left| V^R_{t+1}(x_{t}, \hat{\xi}_{t+1})- \tilde{V}^R_{t+1}(x_{t}, \hat{\xi}_{t+1}) \right| \leq \sum_{t'\in[T]\setminus[t+1]} \max_{\nu\in[|\mathcal{S}_{t+1}|]} L_\nu\beta_J\int_0^1 \Phi_{t'}(z) d z. 
\end{equation*}
For any $x_{t-1}$ and $\xi_t$, $t\in[T]$,
\begin{align}
    & |V^R_t(x_{t-1}, \xi_t)-\widetilde{V}^R_t(x_{t-1}, \xi_t)| \nonumber \\
    & 
    \leq \left| \min_{x_t\in\mathcal{X}_t(x_{t-1},\xi_t)} \left\{c_t x_t+ \max_{Q_{t+1}\in\mathcal{Q}_{t+1}} \rho_{Q_{t+1}}( V^R_{t+1}(x_{t},\xi_{t+1})) \right\}- \right. \nonumber \\
    & \hspace{15em} \left. \min_{x_t\in\mathcal{X}_t(x_{t-1},\xi_t)}  \left\{c_t x_t+ \max_{Q_{t+1}\in\mathcal{Q}_{t+1}} \tilde{\rho}_{Q_{t+1}}( \tilde{V}^R_{t+1}(x_{t},\xi_{t+1})) \right\} \right|
    \nonumber \\
    & \leq \max_{x_t\in\mathcal{X}_t(x_{t-1},\xi_t)} \max_{Q_{t+1}\in\mathcal{Q}_{t+1}} 
    |\rho_{Q_{t+1}}( V^R_{t+1}(x_{t},\xi_{t+1}))-\tilde{\rho}_{Q_{t+1}}( \tilde{V}^R_{t+1}(x_{t},\xi_{t+1}))| \nonumber \\
    & \leq \max_{x_t\in\mathcal{X}_t(x_{t-1},\xi_t)} \max_{Q_{t+1}\in\mathcal{Q}_{t+1}} |\rho_{Q_{t+1}}( V^R_{t+1}(x_{t},\xi_{t+1}))- \rho_{Q_{t+1}}( \tilde{V}^R_{t+1}(x_{t},\xi_{t+1})) | \nonumber \\
    &  \hspace{12em} +|\rho_{Q_{t+1}}( \tilde{V}^R_{t+1}(x_{t},\xi_{t+1})) - \tilde{\rho}_{Q_{t+1}}( \tilde{V}^R_{t+1}(x_{t},\xi_{t+1}))| \nonumber \\
    & \leq \max_{x_t\in\mathcal{X}_t(x_{t-1},\xi_t)} \max_{Q_{t+1}\in\mathcal{Q}_{t+1}}
    \left|\sum_{l\in[|\mathcal{S}_{t+1}|]} q_{t+1,l} \sum_{j\in[K_{t+1}]} \psi_{t+1, l,j} (V^R_{t+1}(x_{t},\xi_{t+1, j})- \tilde{V}^R_{t+1}(x_{t},\xi_{t+1, j})) 
    % {\color{red}\sum_{l\in[|\mathcal{S}_{t+1}|]} q_{t+1,l}\psi_{t+1, l,j}??}
    \right|
    \nonumber \quad\text{(by (\ref{eq-ARSRM-sDiscretized}))}\\
    & \hspace{12em} + \left|
    \sum_{l\in[|\mathcal{S}_{t+1}|]} q_{t+1, l} \int_0^1 F^{\leftarrow}_{\tilde{V}^R_{t+1}(x_{t},\xi_{t+1})}(z) (\sigma(z,s_{t,l})-\tilde{\sigma}(z,s_{t,l})) dz \right|, 
    \label{eq-error-expectation} \nonumber \\
    % & \leq \max_{x_t\in\mathcal{X}_t(x_{t-1},\xi_t)} \max_{Q_{t+1}\in\mathcal{Q}_{t+1}}  \hat{\bbe}_{Q_{t+1}}[| V^R_{t+1}(x_{t},\xi_{t+1})- \tilde{V}^R_{t+1}(x_{t},\xi_{t+1}) |] + 
    % \max_{Q_{t+1}\in\mathcal{Q}_{t+1}} \sum_{l\in[|\mathcal{S}_{t+1}|]} q_{t+1, l} \dd_{\Phi_{t+1}}(\sigma_{s_{t,l}}, \tilde{\sigma}_{s_{t,l}}) \nonumber \\
    & \leq \max_{x_t\in\mathcal{X}_t(x_{t-1},\xi_t)} \max_{k\in[K_{t+1}]} | V^R_{t+1}(x_{t},\xi_{t+1, k})- \tilde{V}^R_{t+1}(x_{t},\xi_{t+1, k}) | +
    \max_{l\in[|\mathcal{S}_{t+1}|]} \dd_{\Phi_{t+1}}(\sigma_{s_{t,l}}, \tilde{\sigma}_{s_{t,l}}) \nonumber \\
    & \leq \max_{x_t\in\mathcal{X}_t(x_{t-1},\xi_t)} \max_{k\in[K_{t+1}]} | V^R_{t+1}(x_{t},\xi_{t+1, k})- \tilde{V}^R_{t+1}(x_{t},\xi_{t+1, k}) | +
    \max_{\nu\in[|\mathcal{S}_{t+1}|]} L_\nu\beta_J\int_0^1 \Phi_{t+1}(z) d z \nonumber.
\end{align}
The fifth inequality holds because 
% where 
% $V_{t+1}=\tilde{V}_{t+1}=0$, 
$$
\psi_{t+1,l,k}:=\int_{\pi_{t+1, k}}^{\pi_{t+1, k+1}} \sigma(z,s_{t,l})dz=\sum_{j\in[J]\cup\{0\}} \sigma_j(s_{t,l}) \int_{\pi_{t+1, k}}^{\pi_{t+1, k+1}}\ind_{[z_j,z_{j+1})}(z) dz,
$$
% {\color{red}
% and $\hat{\bbe}_{Q_{t+1}}$ is the distorted expectation with 
% $$
% % \mathbb{P}(\xi_{t+1} = \xi_{t+1, k})
% \hat{\bbe}_{Q_{t+1}}[V^R_{t+1}(x_{t},\xi_{t+1, k})- \tilde{V}^R_{t+1}(x_{t},\xi_{t+1, k})] = 
% =\sum_{l\in[|\mathcal{S}_{t+1}|]} q_{t+1,l}\psi_{t+1, l, k}.
% $$
% }
% Since 
and
$\sum_{k\in[K_{t+1}]} \sum_{l\in[|\mathcal{S}_{t+1}|]} q_{t+1,l}\psi_{t+1, l,k}=1$, %and 
the last inequality follows from \cite[Proposition 4.1]{WaX20}.
% Assume for stage $t+1$ that 
Then
\begin{align*}
    |V^R_t(x_{t-1}, \xi_t)-\tilde{V}^R_t(x_{t-1}, \xi_t)| & \leq \sum_{t'\in[T]\setminus[t+1]} \max_{\nu\in[|\mathcal{S}_{t'}|]} L_\nu\beta_J\int_0^1 \Phi_{t'}(z) d z
    +\max_{\nu\in[|\mathcal{S}_{t+1}|]} L_\nu\beta_J\int_0^1 \Phi_{t+1}(z) d z  \\
    & \leq \sum_{t'\in[T]\setminus[t]} \max_{\nu\in[|\mathcal{S}_{t'}|]} L_\nu\beta_J\int_0^1 \Phi_{t'}(z) d z.
\end{align*}
The proof is complete. 
\hfill \Box

The error bound can be established analogously for problem (\ref{eq-MARSRM})
when it is used to approximate its counterpart with general risk spectrum.
The next corollary states this.
% {\color{blue}
\begin{corollary}
\label{cor-error-step}
Suppose the conditions in Theorem \ref{thm-error-step-DR} hold. Then
% Let 
% $\sigma_{s_\nu}\in\mathscr{G}$ be Lipschitz continuous with modulus being bounded by $L_\nu$ and $\tilde{\sigma}_{s_{t,\nu}}$ be a projection of $\sigma_{s_{t,\nu}}$ on $\mathscr{G}_J$ for $\nu\in|\mathcal{S}_t|$. 
% Under Assumption \ref{as-Phi},
    %hold for each $t\in[T]$ and 
    %Then 
    \begin{equation*}
        |V_t(x_{t-1}, \xi_t)-\tilde{V}_t(x_{t-1}, \xi_t)|\leq \sum_{t'\in[T]\setminus[t]} \sum_{\nu\in[|\mathcal{S}_{t'}|]} 
        q_{t', \nu} 
        L_\nu\beta_J\int_0^1 \Phi_{t'}(z) d z,\;\;  
       \text{for}\; t\in[T],
    \end{equation*}
where $q_{t,\nu} = \mathbb{P}(s_t=s_{t, \nu})$, for $\nu\in[|\mathcal{S}_t|]$ and $t\in[T]\setminus\{1\}$.
\end{corollary}

\section{Numerical Experiments}

We have carried out some 
numerical experiments 
to examine the performance of the proposed multistage linear program models with ARSRM and DR-MARSRM and the computational schemes.
%from applying inner and lower approximation algorithms, 
%and compare the results with the deterministic risk-neutral and risk-averse model on an asset allocation model with transaction costs.
In this section we report the test results.
%some numerical results to 

\subsection{Multistage Asset Allocation Problem with Transaction Costs}

A decision maker has an initial total capital normalized to $1$ at time $t=0$. The DM
plans to invest the capital to some assets 
with random returns and adjust the allocations at each time period $t$ 
over the fixed number of time periods $t=1,2,\cdots,T$.
%time horizon.
At the beginning of time period 
$t$,
%allocate a certain fund to some assets
%to maximize the profit over a time horizon.
%At stage $t$ the decision 
%The total asset is normalized to $1$. 
%Let 
$x_t$ is the vector of 
re-balancing/allocations 
of all available funding at the time to
% %the
% %vector of proportions 
% allocations 
% of all available funding to 
% %be allocated to 
assets
%denote the decision vector 
%which
%the allocations (
in units of a multiple of a base currency, and $\xi_t$ denotes 
the vector of 
gross returns,
%at the end of period $t$.
%per stage, 
%that is, 
which are the ratios of the asset 
prices at the end of 
time period $t$
%to that 
over the asset prices at the end of
%in stage 
time period $t-1$.
%which represents the only random parameters in the model.
We consider the case where transaction costs are proportional to the value of the assets sold or bought. 
Let $f_t$ denote the relative cost, the balancing equation between stage $t-1$ and stage $t$ portfolios has to be modified to include the total
cost of $f_t\bm{e}^T|x_t-x_{t-1}|$, where $\bm{e}$ is the vector with unit components and the absolute value operation
$|\cdot|$
%function applies 
is taken component-wise.
%Converting to linear formulation,
%we obtain the following model (see 
The problem is considered in \cite{KoM15}.
Here we revisit it with our proposed MARSRM
model:
\begin{equation*}
\begin{split}
    V_{t}(x_{t-1},\xi_{t}) = \min_{x_t\geq0,z_t} \; & -\bm{e}^{\top}x_t+\rho_{Q_{t+1}}(V_{t+1}(x_t,\xi_{t+1})) \\
    \st \;\;\; & \bm{e}^{\top}x_t+f_t\bm{e}^{\top}z_t =\xi_t^{\top}x_{t-1}, \\
    & |x_t-x_{t-1}|\leq z_t.
\end{split}
\end{equation*}
%We s
Suppose that all the initial capital is used for investment into stocks, meaning transaction costs are constant for the first stage and can be therefore omitted. 
The first-stage stochastic program is in the following form:
\begin{equation*}
\begin{split}
    V_1(x_0,\xi_1) = \min_{x_1\geq 0} \;\; & -\bm{e}^{\top} x_1+\rho_{Q_2}(V_2(x_1,\xi_2)) \\
    \st \;\; & \bm{e}^{\top}x_1 =1.
\end{split}
\end{equation*}

In the case that the investor is ambiguous about the true probability distribution of random preference $s$ at each stage, we propose the distributionally preference robust counterpart of the multistage asset allocation problem to mitigate the risk arising from the ambiguity with the following subproblem at stage $t$:
\begin{equation*}
\begin{split}
    V^R_{t}(x_{t-1},\xi_{t}) = \min_{x_t\geq0,z_t} \; & -\bm{e}^{\top}x_t+ 
    \max_{Q_{t+1}\in\mathcal{Q}_{t+1}} \rho_{Q_{t+1}}(V^R_{t+1}(x_t,\xi_{t+1})) \\
    \st \;\;\; & \bm{e}^{\top} x_t+f_t\bm{e}^{\top}z_t =\xi_t^{\top}x_{t-1}, \\
    & |x_t-x_{t-1}|\leq z_t
\end{split}
\end{equation*}
and first-stage problem 
\begin{equation*}
\begin{split}
    V^R_{1}(x_0,\xi_{1}) = \min_{x_1\geq 0} \;\; & -\bm{e}^{\top}x_1
    +\max_{Q_2\in\mathcal{Q}_2} \rho_{Q_2}(V^R_2(x_1,\xi_2)) \\
    \st \;\; & \bm{e}^{\top}x_1 =1,
\end{split}
\end{equation*}
where the ambiguity set $\mathcal{Q}_{t+1}$ is state-dependent moment-based ambiguity set defined in (\ref{eq-mathcalQ_t}).

In \cite{kozmik2015evaluating}, the authors formulate the decision making as a multistage problem where the DM's risk preference is characterized by the combination of expectation and $CVaR$ and stagewise-dependent risk aversion coefficient and confidence level.

\subsection{Set-up of The Tests}

%For our experimental setup, the gross return 
In this experiment, we assume that 
$\xi_t$ is 
%{\color{red} 
interstage independent
%HX: what do you mean?}
%{\color{blue} Q: I mean 
which means that $\xi_t$ and $\xi_{t+1}$ are
%is 
independent.
%{\color{blue} 
%but $\xi_{t}$ has dependence across assets.
%}
%{\color{red}, but have dependence across stock market indices. Rephrase this properly} 
We assume that $\xi_{t,r}$ has a log-normal distribution of the form $\exp(\phi_{t,r})$, where 
% $b_{t,r}\in\R_+$ and 
$\phi_{t,r}, r\in\mathcal{R}$ has a multivariate normal distribution with mean $\tilde{\mu}_{t,r}$ and variance $\tilde{\sigma}_{t,r}^2$ and correlation coefficient $\tilde{\gamma}_{t,r,r'}$.
We consider the random spectral function $\sigma(\cdot,s)=\lambda\ind_{[0,1]}(\cdot)+(1-\lambda)\frac{1}{1-\alpha}\ind_{[\alpha,1]}(\cdot)$ with $s = (\lambda,\alpha)$ and corresponding step-like approximation: 
\begin{equation}
    \sigma_{t}(z,\lambda,\alpha) =
    \begin{cases}
        \lambda, & \text{for\;}z\in[0,z_{i}),  \\
        \left(1-\lambda z_{i}+(1+z_{i+1})\frac{1-\lambda\alpha}{1-\alpha}\right)/(z_{i}-z_{i+1}), & \text{for\;}z\in[z_{i},z_{i+1}), \\
        \lambda+\frac{1-\lambda}{1-\alpha}, & \text{for\;} z\in[z_{i+1},1),
    \end{cases}
\end{equation}
where $i$ is the index such that $\alpha\in[z_{i},z_{i+1})$.
We set the number of breakpoints as $J$ and $z_0,\ldots,z_{J+1}$ as $\{0,\frac{1}{J},\ldots,\frac{J-1}{J},1\}$. 
The scenarios tree are generated in the form of re-combining tree with $K_t$ realizations per stage.
All of the tests are carried out in MATLAB R2023b installed on a PC (16GB, CPU 2.3 GHz)
with an Intel Core i7 processor. We use Gurobi under YALMIP environment \cite{lofberg2004yalmip} to solve all linear subproblems.
Table \ref{tab-parameter} summarizes the values of the parameters and the size of tree we use in our model.

\begin{table}[htb]
    \centering
    \caption{Parameters and of problem}
    \vspace{-0.2cm}
    \begin{tabular}{ccc}
    \hline
       Index set of assets $\mathcal{R} = \{1,\ldots,4\}$ 
        % {\color{red}\inmat{HX: what is this?}} {\color{blue}\inmat{Q: the index set of asset.}} 
        &  &  \\
        % \hline
        $\tilde{\mu}_{\cdot, r} = (0.6, 0.6, 0.6, 0.6)$ &  &  \\
        $\tilde{\sigma}_{\cdot, r} = (0.3, 0.3, 0.3, 0.3)$ & &  \\
        $\tilde{\gamma}_{\cdot, r, r'} = 0.5$ \; $\inmat{for}\;  r, r'\in\{1, 3\}$, $r\neq r'$ &    &  \\
        $N_t = 1000$, $\tilde{N} = 10$ &  &  \\
        $J = 10$ &  &  \\
        $K_t\in\{100, 50, 20, 10, 10\} \; \inmat{for}\;  t\in\{2, 3, 5, 10\}$ &  &  \\
        % 15 & 20 & $\approx 10^{15}$ \\
        \hline
    \end{tabular}
    \label{tab-parameter}
\end{table}

%\textbf{(i) MARSRM.}
\subsubsection{MARSRM}
% For the non-DRO multistage problem, 
% % experiments where the probability distribution of $s$ is given, 
At each stage,
we use the Voronoi tessellation to generate the stagewise probability distribution of $s$ with the same pre-selected support $\mathcal{S}_t = \mathcal{S}=\{s_1,\ldots,s_{\tilde{N}}\}$.
We generate points $s_1,\ldots,s_N$ by using the beta distribution $\rm{Be}(2,2)$ and 
let $\{s_1, \ldots, s_N\}$ be a Voronoi tessellation centred at $\mathcal{S}$, i.e., $S_l\subset\{s\in\mathcal{S}: \|s-s_l\|=\min_{l'\in[N]}\|s-s_{l'}\|\}$ for $l=1,\ldots,[\tilde{N}]$ are pairwise disjoint subsets forming a partition of $\mathcal{S}$.
Let $s_1,\ldots,s_{N_t}$ be an iid sample of $s$ where the sample size $N_t$ is usually larger than $N$ and $N_{t,l}$ denote the number of samples falling into area $S_l$.
Define the probability distribution $Q_t$ with
\begin{equation*}
    q_{t, l} := \sum_{l \in[\tilde{N}]} \frac{N_{t, l}}{N_t} \delta_{s_l}(\cdot),
\end{equation*}
where $\delta_{s_l}(\cdot)$ denotes the Dirac probability measure located at $s_l$.
% We sample $N$ paths at each iteration.
% Table \ref{tab-size} shows the sizes of the empirical scenario trees for the test problem.
We test and compare the following versions of the MARSRM model:
\begin{enumerate}[(i)]
    \item MARSRM: multistage problem with average ramdomized spectral risk measure (\ref{eq-MARSRM}).
    \item Risk neutral MS: multistage problem with expectation (by setting $\lambda = 1$ and $\alpha = 0$ in (\ref{eq:sigmaCombination})).
    \item Mild risk averse MS: multistage problem with the combination $\lambda \bbe + (1-\lambda) \CVaR$ (we set $\lambda = 0.5$ and $\alpha = 0.8$ in (\ref{eq:sigmaCombination})).
    \item Strong risk averse MS: multistage problem with the combination $\lambda \bbe + (1-\lambda) \CVaR$ (we set $\lambda = 0.8$ and $\alpha = 0.8$ in (\ref{eq:sigmaCombination})).
    \item DR-MARSRM: Distributionally preference robust multistage problem with moment-cased ambiguity set. 
\end{enumerate}

%\textbf{(ii) MS-PDRO-ARSRM.}
\subsubsection{DR-MARSRM}
We build the 
%the 
shrinking moment-based ambiguity set $\mathcal{Q}_t$ for the preference parameter, where we generate the support set $\mathcal{S}_t = \{s_{t,1}, \ldots, s_{t,10*t}\}$ with the true uniform distribution $U(0, 1)$, the corresponding 
mean and variance are defined with the sample mean and sample variance.
We test the convergence and the out-of-sample performance of DR-MARSRM.

% \begin{table}[htb]
%     \centering
%     \begin{tabular}{ccc}
%     \hline
%         Stages (T) & Descendants per node ($K_t$) & Total scenarios \\
%         \hline
%         2 & 100 & 100 \\
%         3 & 50 & 2500 \\
%         4 & 50 & 125000 \\
%         5 & 20 & 160000 \\
%         10 & 10 & $10^{9}$ \\
%         % 15 & 20 & $\approx 10^{15}$ \\
%         \hline
%     \end{tabular}
%     \caption{Sizes of re-combining trees}
%     \label{tab-size}
% \end{table}

% All problem are solved by Gurobi through YALMIP in Matlab 2023b on a PC with 2.3GHz CPU and 32GB RAM.

\subsection{Test Results}

%In this subsection, we 

% We have carried out numerical experiments %investigate 
% to examine the performance of
% SDDP algorithms for solving
% the multistage decision making discussed in Section 3.1 and 4.1. In this section, we report the experimental results.
% % Test results are carried out from the following perspectives.

% \noindent
% \textbf{
%(i)
\subsubsection{SDDP Algorithm for Solving MARSRM}
We begin by examining convergence of 
the sequence of upper and lower bounds
%evolution 
for the optimal value of MARSRM when Algorithm~\ref{alg-lowerbound}-\ref{alg-upperbound} is applied to solve the problem,
%with respect to , 
with different number of stages
%and iterations with 
$T\in\{2, 3, 5, 10\}$. 
The lower approximation algorithm is a standard variation of SDDP with 10 scenarios used in each forward pass.
We run each upper approximation algorithm 10 times, using the states visited after $10$ to $100$ cuts have been added in the lower approximation.
Figure \ref{fig:boundsEvolution} shows that the upper and lower bounds move closer as the number of iterations increases. We can see that for $T=2,3,5$,
the gap of the two bounds goes to zero after 7 iterations. In the case that $T=0$, it takes about 10 iterations.
Table \ref{tab:bound-ARSRM} 
%Moreover, 
%we present 
lists the 
lower and upper bounds, 
their gaps
%optimality gap for the SAA problem 
and respective 
%corresponding
computational time.
%for outer and upper approximations. 
%in Table \ref{tab:bound-ARSRM}.
We can see the details
of the bounds and 
%that 
the gaps 
%between the bounds reduces 
as the number of iterations (cuts) increases in the lower approximation.
%Compared to 
The computational time is a bit long because 
%(\ref{eq-ARSRM-sDiscretized}) the 
multiple CVaRs are involved 
%in the ARSRM as defined
in (\ref{eq-ARSRM-sDiscretized}).
%in the literature (although the models are different), the computational 
% numerical results of obtained with the 
% outer and upper approximation, respectively, for several numbers of cuts.
% In addition to the bounds, Table \ref{tab:bound-ARSRM} shows the optimality gap for the SAA problem and the computational time for both approximations.

\begin{figure}[!ht]
    \centering 
    \subfigure{
    \includegraphics[width=0.23\linewidth]{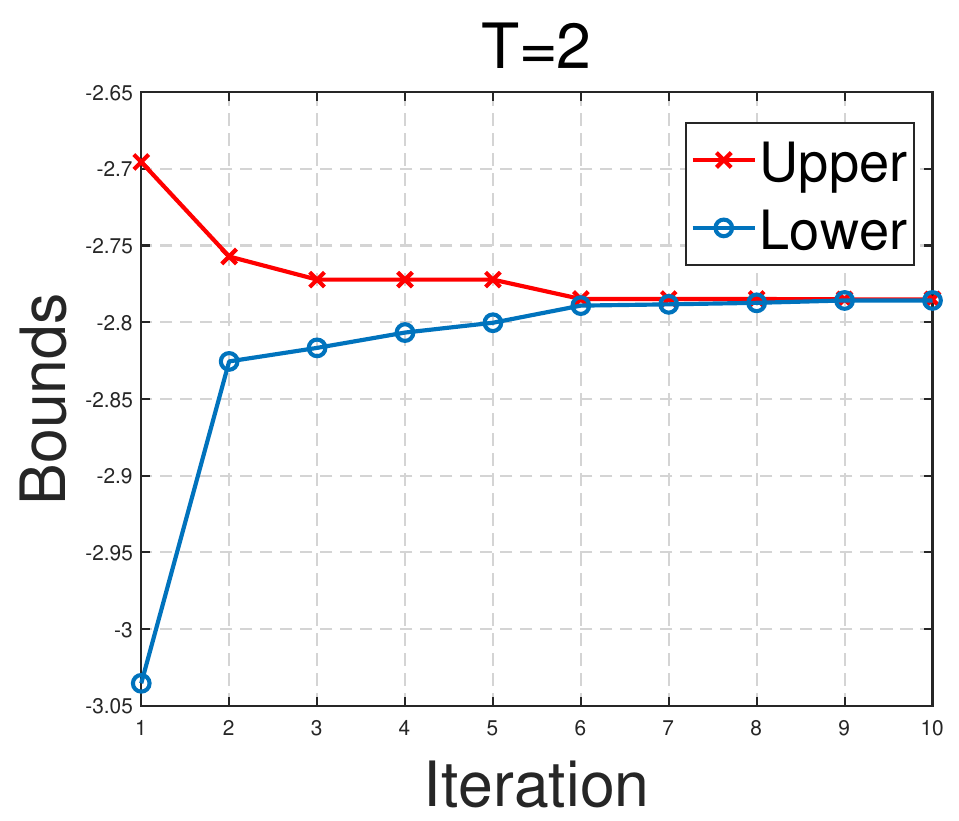}
    }
    \hspace{-1em}
    \subfigure{
    \includegraphics[width=0.23\linewidth]{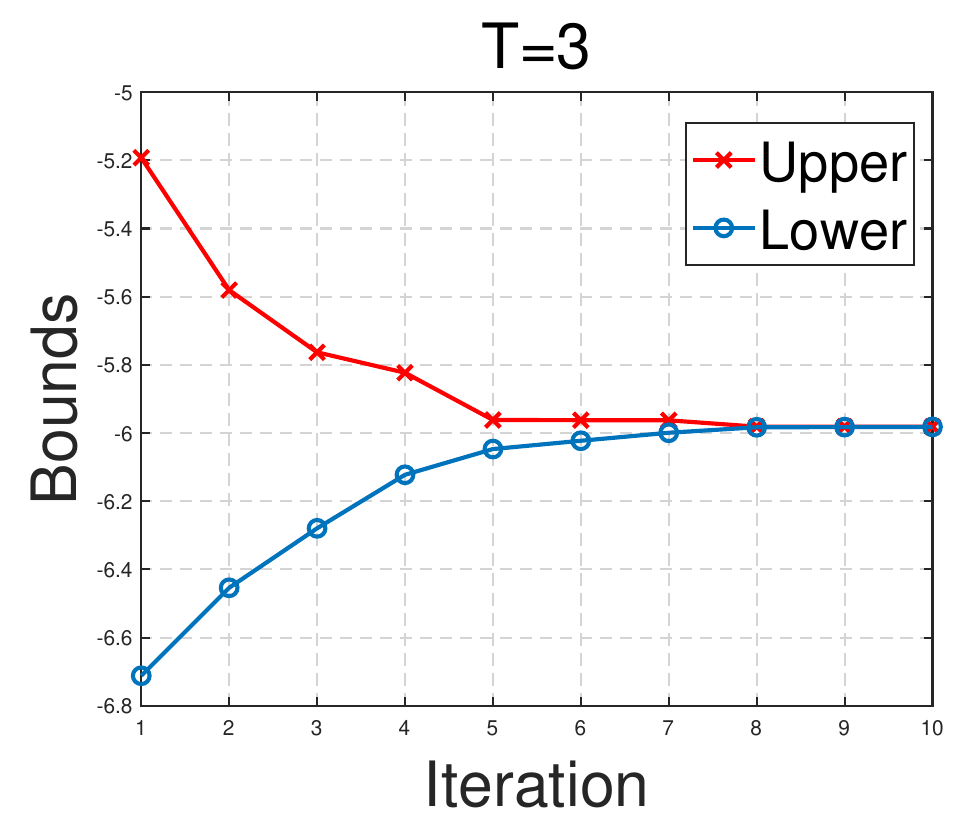}
    }
    \hspace{-1em}
    \subfigure{
    \includegraphics[width=0.23\linewidth]{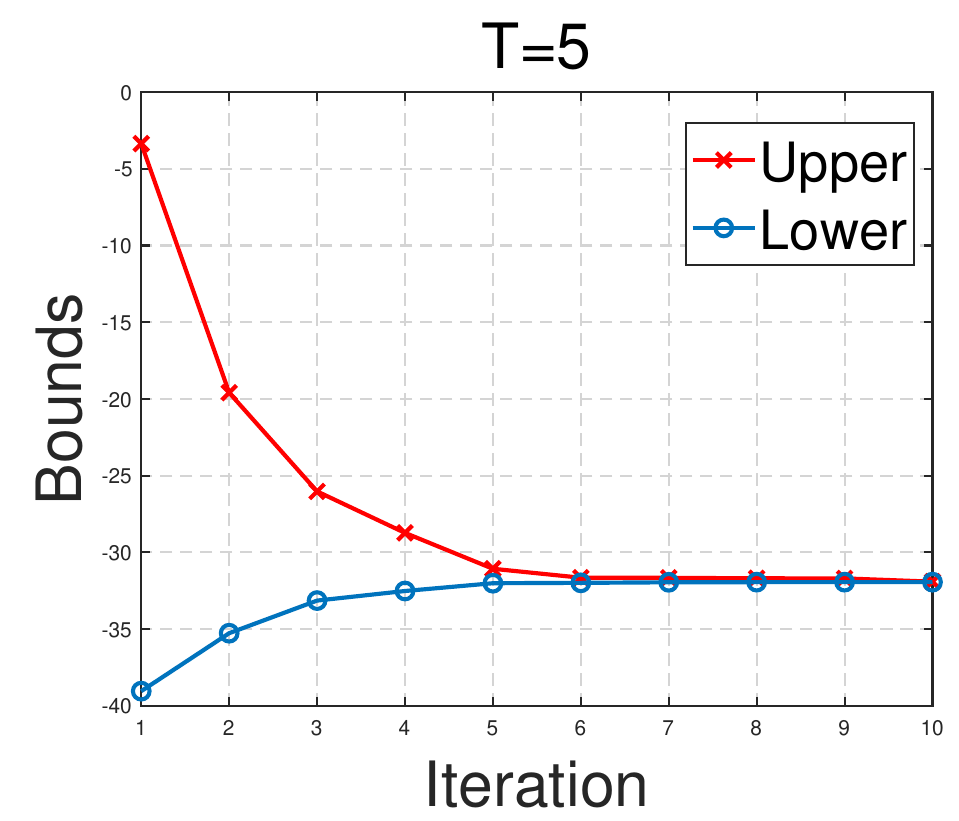}
    }
    \hspace{-1em}
    \subfigure{
    \includegraphics[width=0.23\linewidth]{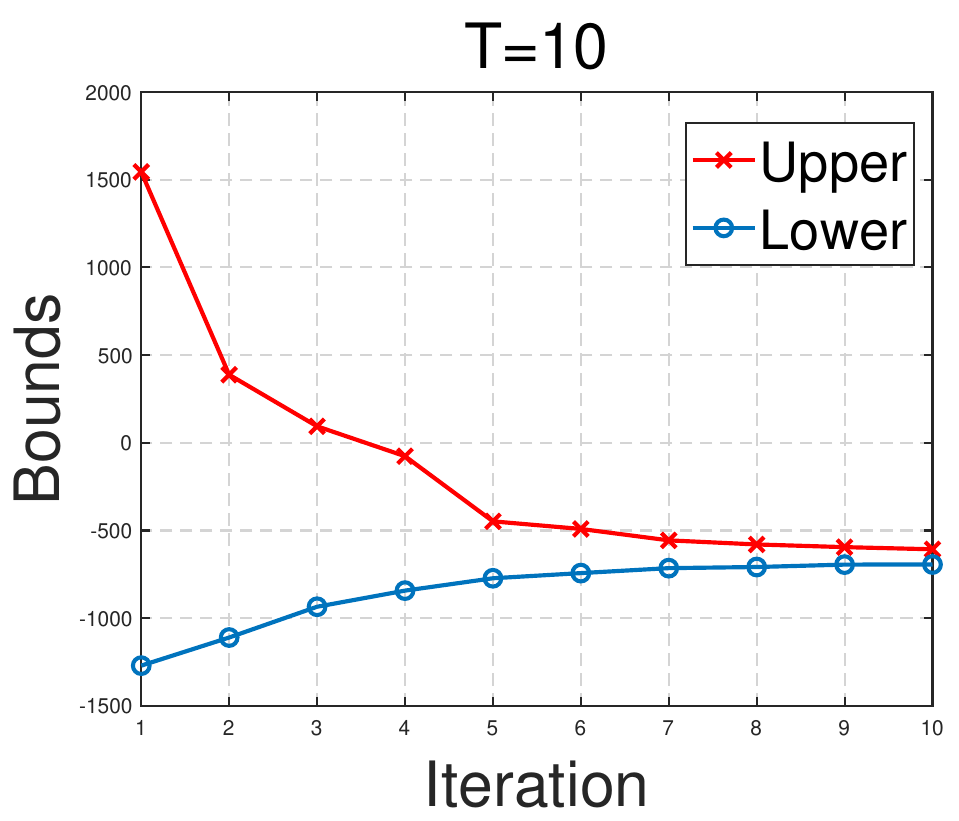}
    }
    \vspace{-0.2cm}
    \captionsetup{font=footnotesize}
    \caption{
    %U
    Convergence of upper and lower bounds %evolution of 
    of the optimal value of MARSRM}
    \label{fig:boundsEvolution} 
\end{figure}

\begin{table}[ht]
    \centering
    \footnotesize
    \captionsetup{font=footnotesize}
    \caption{Results of MARSRM, $T=10$}
    \vspace{-0.1cm}
    \begin{tabular}{cccccc}
    \hline
        Cuts & Lower bound & Upper bound & Gap & Time lower (s) & Time upper (s) \\
        \hline
        10 & -1270.7 & 1544.6 & 2815.3 & 79.0 & 34.5 \\
        20 & -1110.4 & 388.1 & 1498.5 & 35.4 & 93.8 \\
        30 & -934.9 & 94.0 & 1028.9 & 57.8 & 212.9 \\
        40 & -843.1 & -76.4 & 766.6 & 77.7 & 431.7 \\
        50 & -772.5 & -447.9 & 324.6 & 112.6 & 774.4 \\
        60 & -743.1 & -491.3 & 251.8 & 164.1 & 1284.4 \\
        70 & -714.8 & -556.4 & 158.4 & 221.0 & 1996.4 \\
        80 & -708.4 & -579.7 & 128.7 & 290.1 & 2945.0 \\
        90 & -694.6 & -595.0 & 99.6 & 364.0 & 4159.0\\
        100 & -694.1 & -607.2 & 86.8 & 451.0 & 5695.2 \\
        \hline
    \end{tabular}
    \label{tab:bound-ARSRM}
\end{table}

% {\color{red} 

% We also test computational efficiency by comparing 
%the test results 
We also compare lower bounds and their convergence 
for the optimal values of different %multistage problems
%: with
%including
version of multistage ARSRM models
including MARSRM,
risk neutral case, 
mild risk averse case and strong risk averse case
%w.r.t. iterations 
as the number of iterations increases, 
with different number of stages
%with 
$T\in\{2, 3, 5, 10\}$. 
Figures \ref{fig:upperboundsComparison} and \ref{fig:lowerboundsComparison} 
%show 
depict evolution of 
upper and lower bounds
%evolution of different 
for each of the problems, 
where the red line corresponds to ARSRM, the purple line corresponds to risk neutral multistage problem ($\lambda = 0$, $\alpha = 0$), the blue line corresponds to mild risk averse multistage problem ($\lambda = 0.5$, $\alpha = 0.8$), and the yellow line corresponds to strong risk averse multistage problem ($\lambda = 0.8$, $\alpha = 0.8$). We can see that the algorithm also displays similar convergence when it is applied to other versions of MARSRM.

% }

\begin{figure}[!ht]
    \centering 
    \subfigure{
    \includegraphics[width=0.2\linewidth]{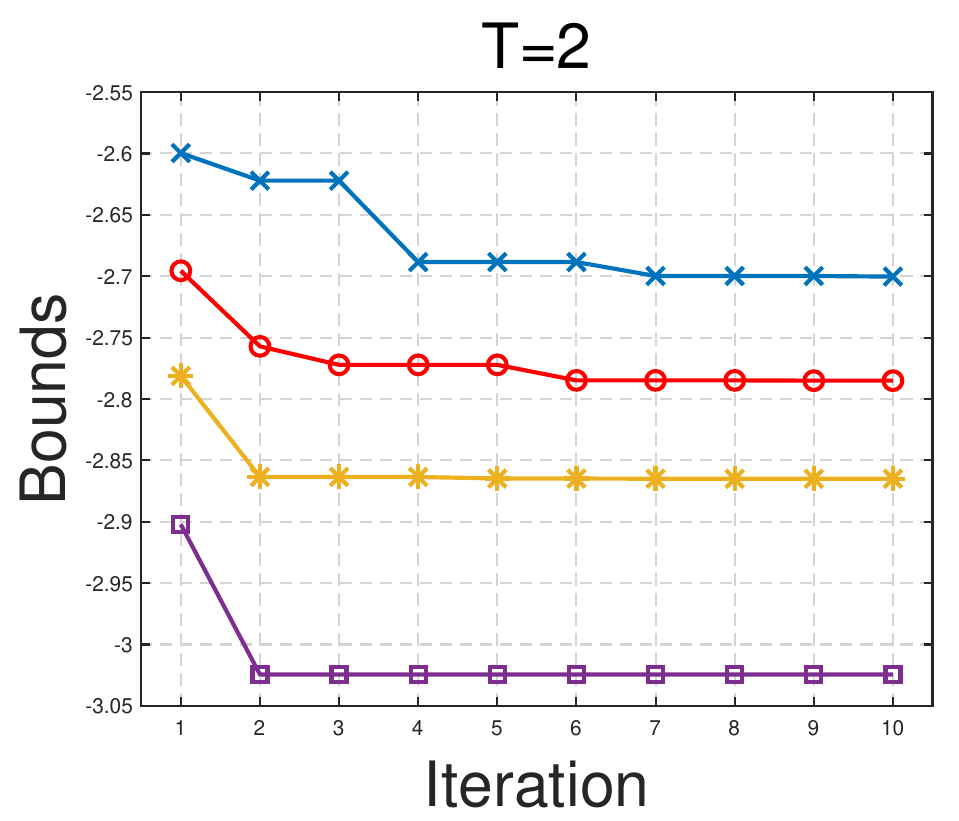}
    }
    \hspace{-1em}
    \subfigure{
    \includegraphics[width=0.2\linewidth]{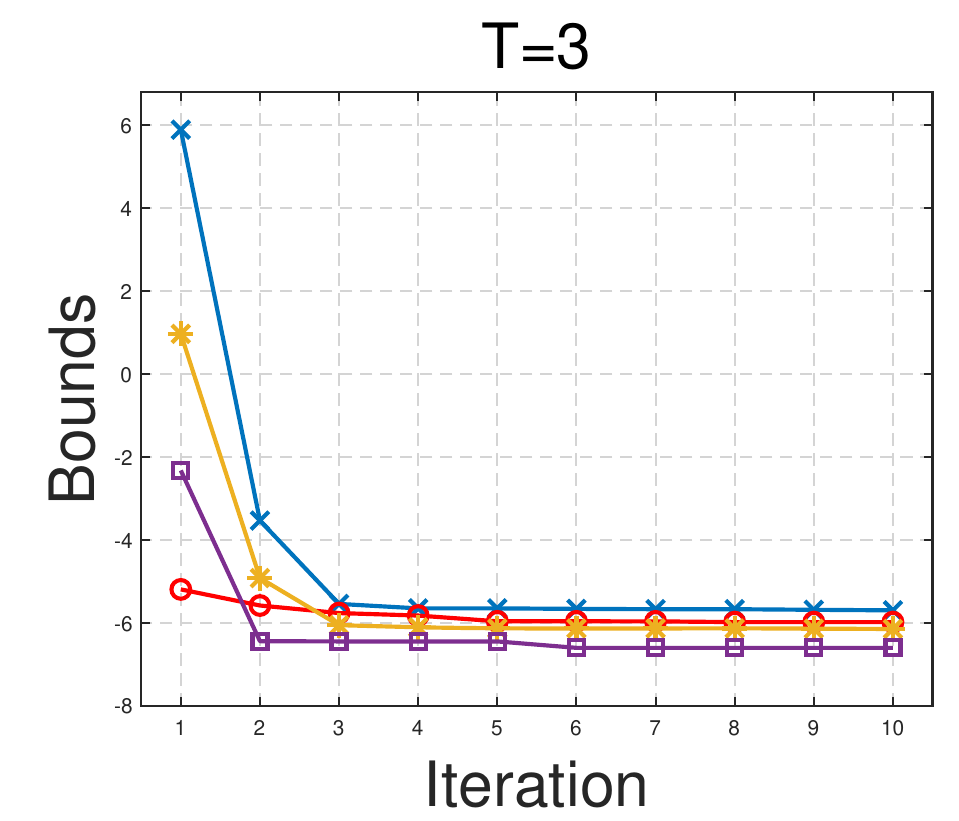}
    }
    \hspace{-1em}
    \subfigure{
    \includegraphics[width=0.2\linewidth]{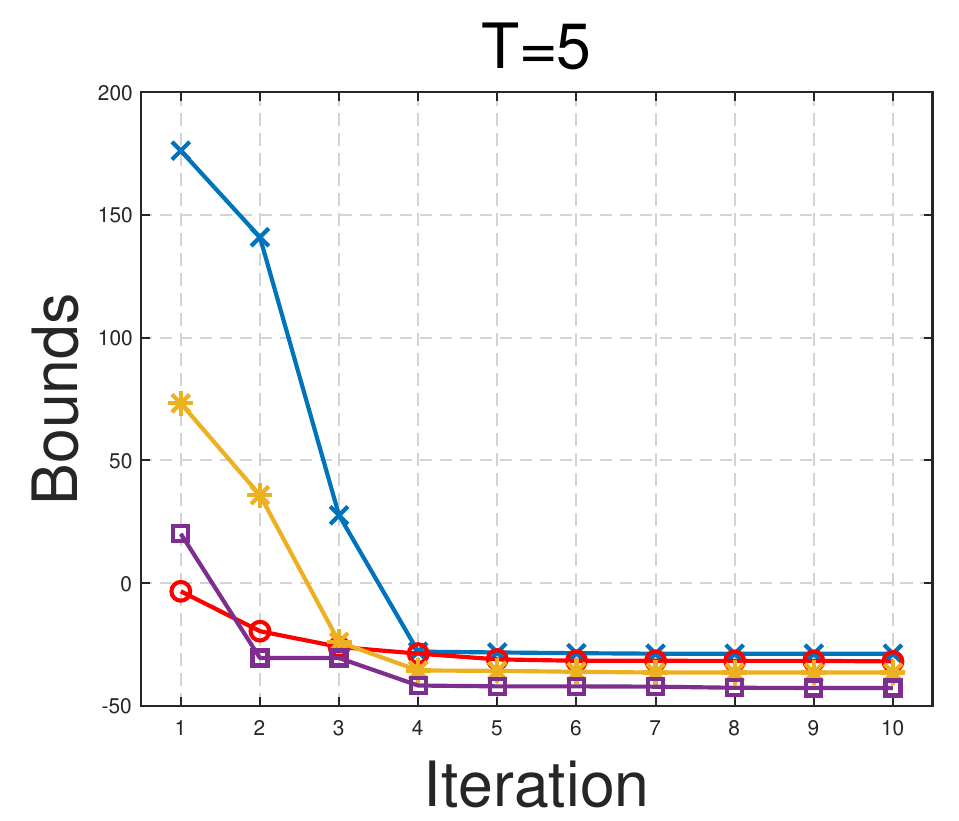}
    }
    \hspace{-1em}
    \subfigure{
    \includegraphics[width=0.2\linewidth]{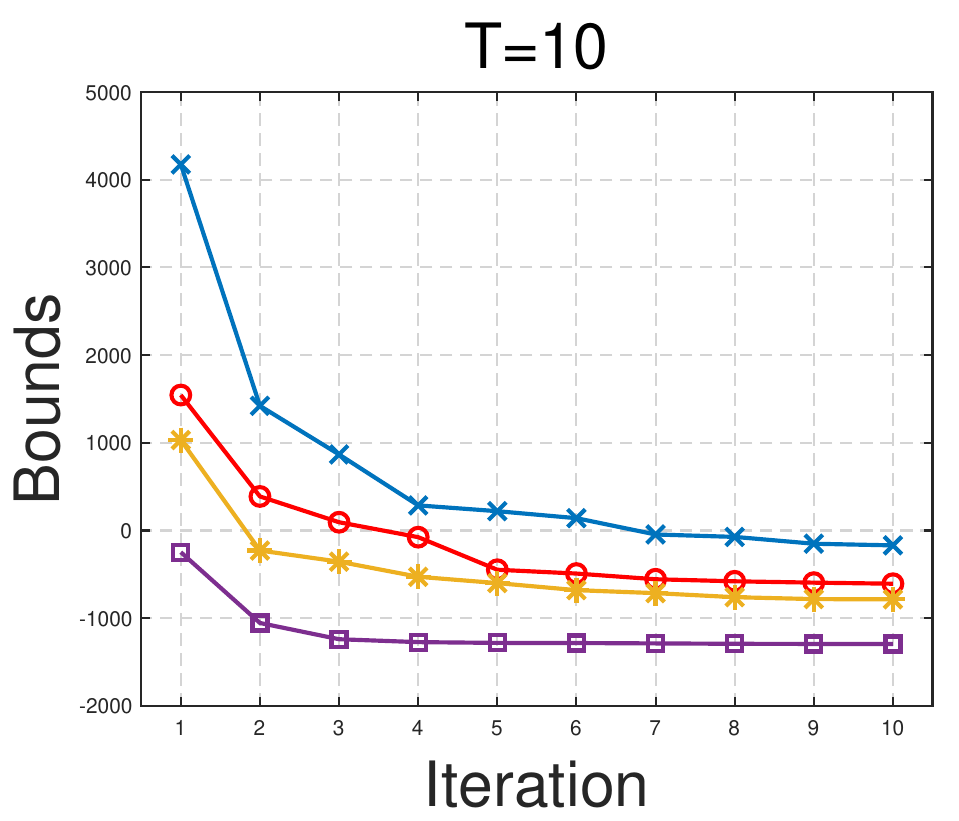}
    }
    \hspace{-1em}
    \subfigure{
    \includegraphics[width=0.15\linewidth]{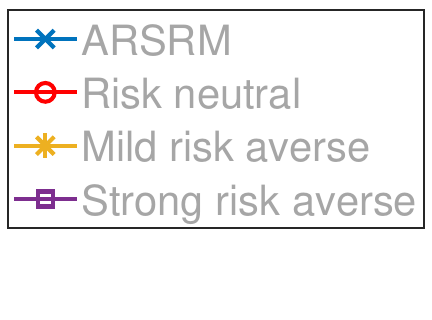}
    }
    \vspace{-0.2cm}
    \captionsetup{font=footnotesize}
    \caption{\uline{Upper bound} evolution of different multistage problems.
    }
    \label{fig:upperboundsComparison} 
\end{figure}

\begin{figure}[!ht]
    \centering 
    \subfigure{
    \includegraphics[width=0.2\linewidth]{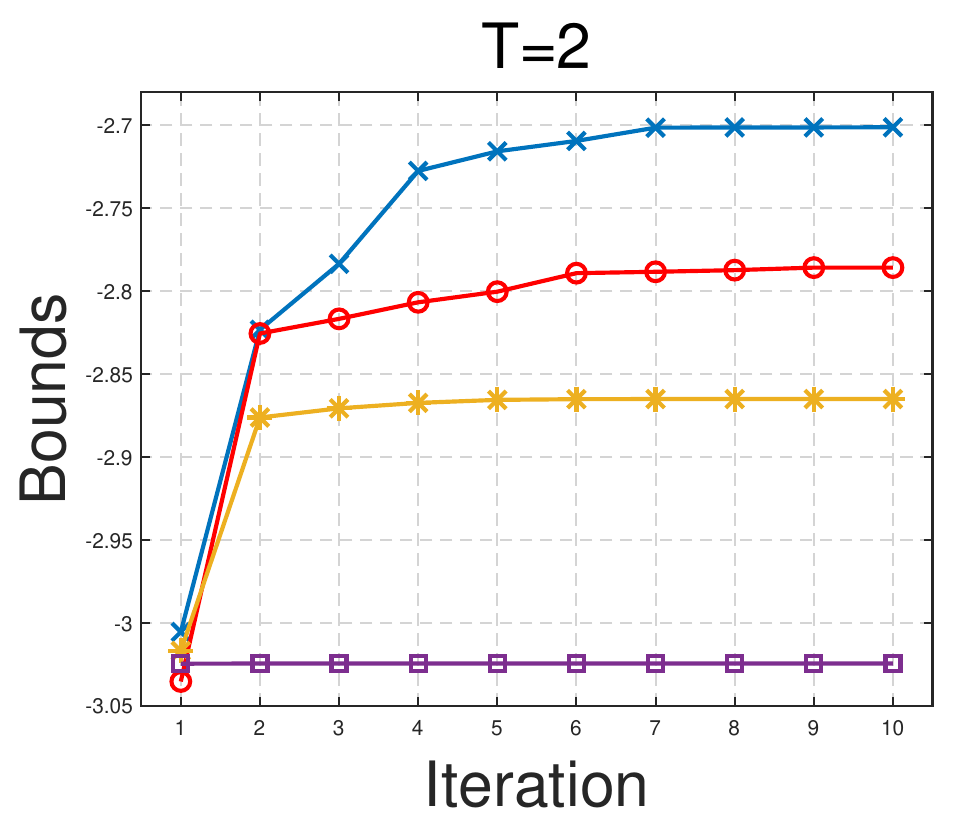}
    }
    \hspace{-1em}
    \subfigure{
    \includegraphics[width=0.2\linewidth]{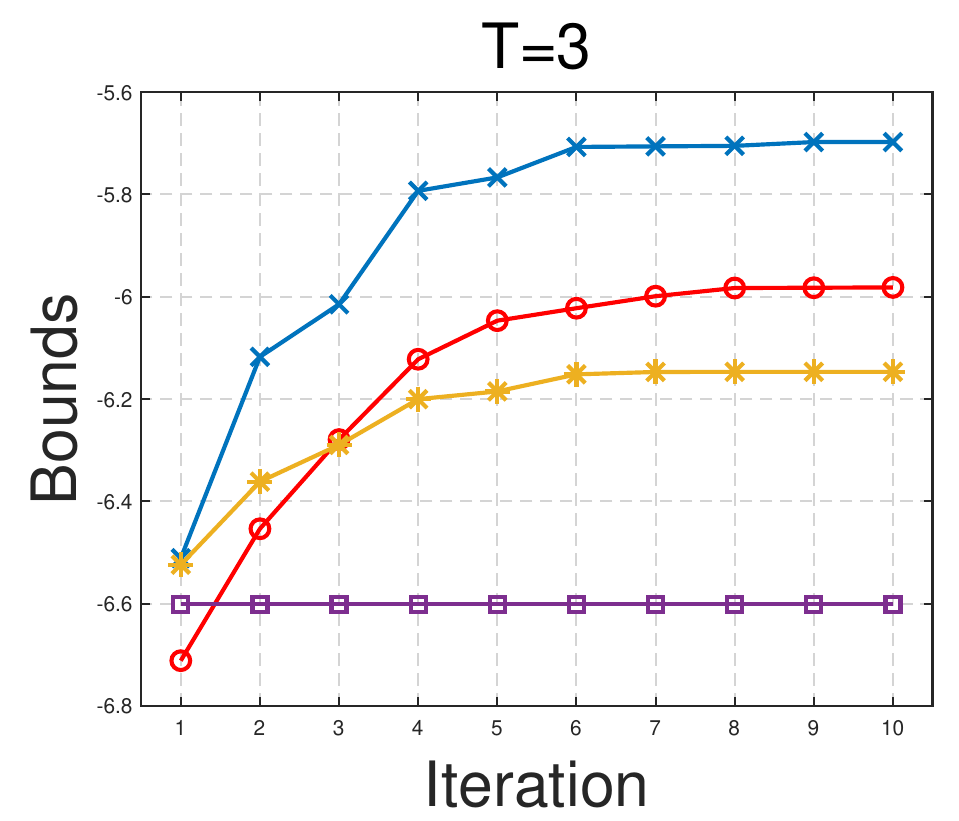}
    }
    \hspace{-1em}
    \subfigure{
    \includegraphics[width=0.2\linewidth]{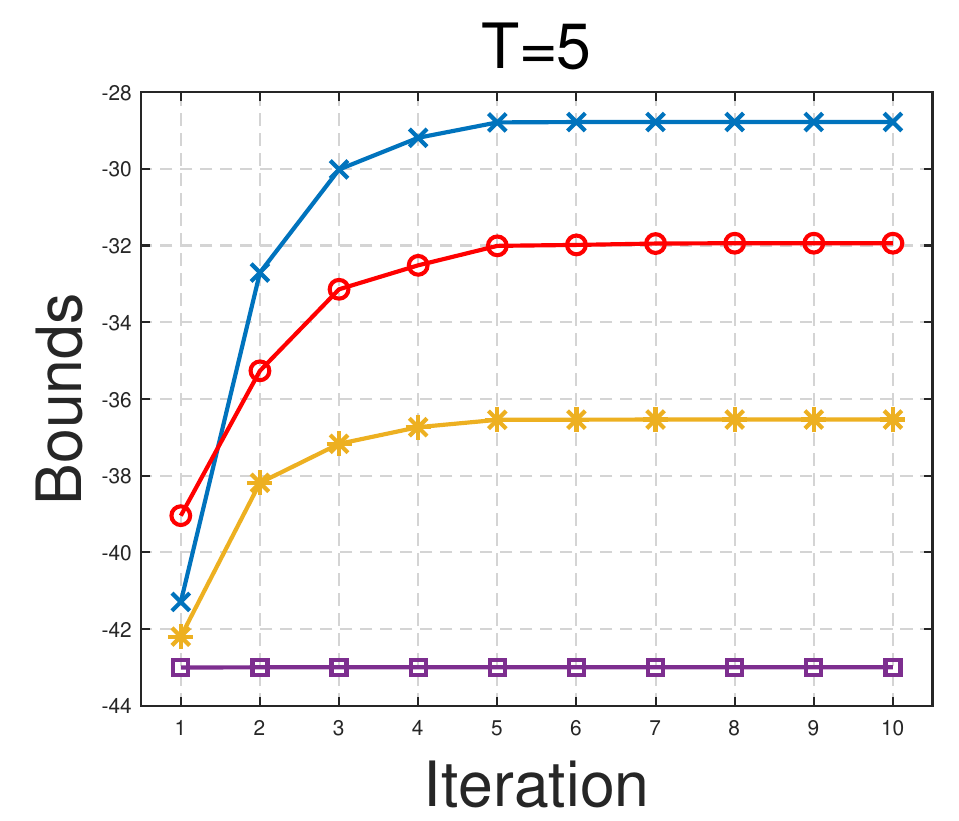}
    }
    \hspace{-1em}
    \subfigure{
    \includegraphics[width=0.2\linewidth]{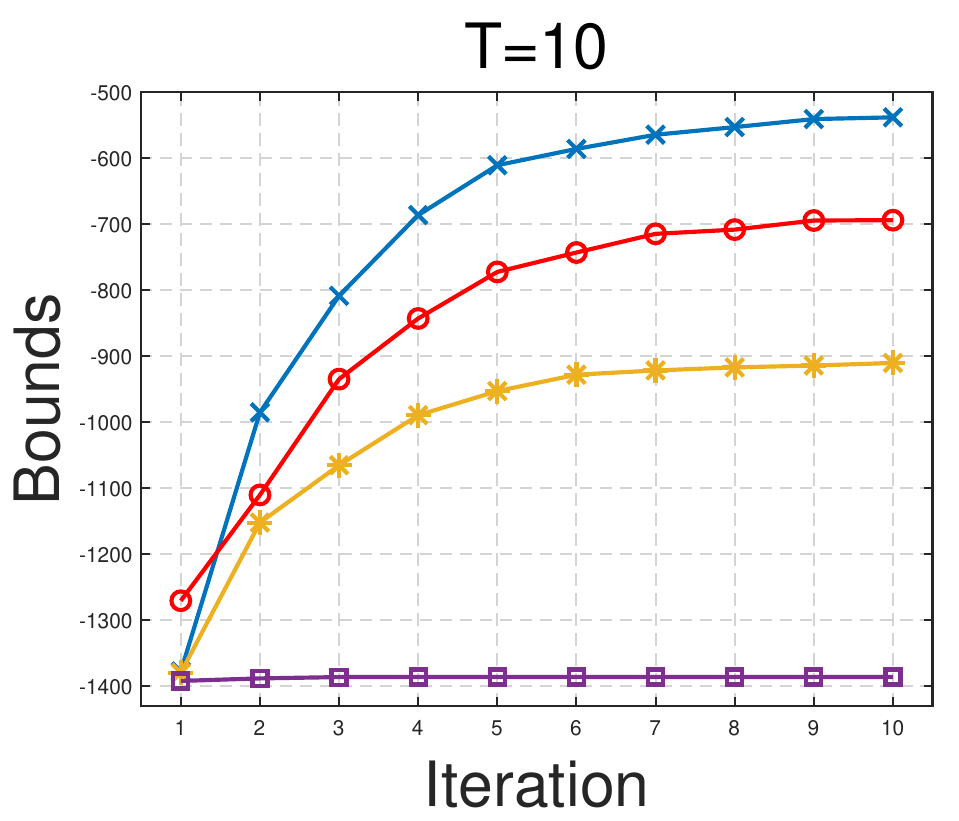}
    }
    \hspace{-1em}
    \subfigure{
    \includegraphics[width=0.15\linewidth]{comparisonlegend.pdf}
    }
    \vspace{-0.2cm}
    \captionsetup{font=footnotesize}
    \caption{
   Comparison of  \uline{lower bounds} 
    %evolution 
    and their convergence
    for the optimal values of different multistage problems.
    }
    \label{fig:lowerboundsComparison} 
\end{figure}

% \noindent
% \textbf{
%(i) 
\subsubsection{SDDP algorithm for solving DR-MARSRM.}
This part of numerical tests is 
%concerned with 
focused on the SDDP approach for solving the  DR-MARSRM model. 
We repeat the tests that we have done for the MARSRM model. Figure \ref{fig:boundsEvolutionDRO} depicts convergence of the lower and upper bounds
as the number of iterations (cuts) increases and Table \ref{tab:bound-ARSRM-DRO} reports the details of bounds, gaps and computational time. We can see that Algorithm~\ref{alg-lowerbound-dro}-\ref{alg-upperbound-dro} delivers even better overall numerical performances primarily because we use the multicut version of piecewise linear approximation.
%distributional robust MARSRM

% We 
% %we will test the
% examine convergence 
% % and 
% % out-of-sample performance 
% of the algorithm.
% %DRO 
% % model.
% We compare the upper and lower bounds for the optimal value of the DR-MARSRM with those 
% of risk averse multistage problem 
% %with respect to different numbers of stage and iterations 
% as the number of iterations increases, with
% different number of stages $T \in\{2, 3, 5, 10\}$. 
% Moreover, 
% we 
% %demonstrate 
% carry out comparative out-of-sample tests on the upper and lower bounds for the DR-MARSRM problem, risk neutral and risk averse multistage problems.

\begin{table}[ht]
    \centering
    \footnotesize
    \captionsetup{font=footnotesize}
    \caption{Results of DR-MARSRM, $T=10$}
    \vspace{-0.1cm}
    \begin{tabular}{cccccc}
    \hline
        Cuts & Lower bound & Upper bound & Gap & Time lower (s) & Time upper (s) \\
        \hline
        1 & -13.1 & 3.2 & 16.3 & 66.7 & 32.3 \\
        2 & -12.9 & -4.6 & 8.4 & 22.8 & 70.5 \\
        3 & -12.9 & -7.0 & 5.9 & 32.9 & 143.1 \\
        4 & -12.9 & -9.1 & 3.8 & 39.4 & 271.3 \\
        5 & -12.9 & -9.8 & 3.1 & 51.9 & 511.7 \\
        6 & -12.9 & -10.7 & 2.2 & 63.4 & 765.8 \\
        7 & -12.9 & -10.9 & 2.0 & 77.0 & 1151.2 \\
        8 & -12.9 & -11.1 & 1.8 & 91.5 & 1692.2 \\
        9 & -12.9 & -11.8 & 1.1 & 99.2 & 2250.6 \\
        10 & -12.9 & -12.5 & 0.4 & 114.3 & 3046.3 \\
        \hline
    \end{tabular}
    \label{tab:bound-ARSRM-DRO}
\end{table}

\begin{figure}[!ht]
    \centering 
    \subfigure{
    \includegraphics[width=0.23\linewidth]{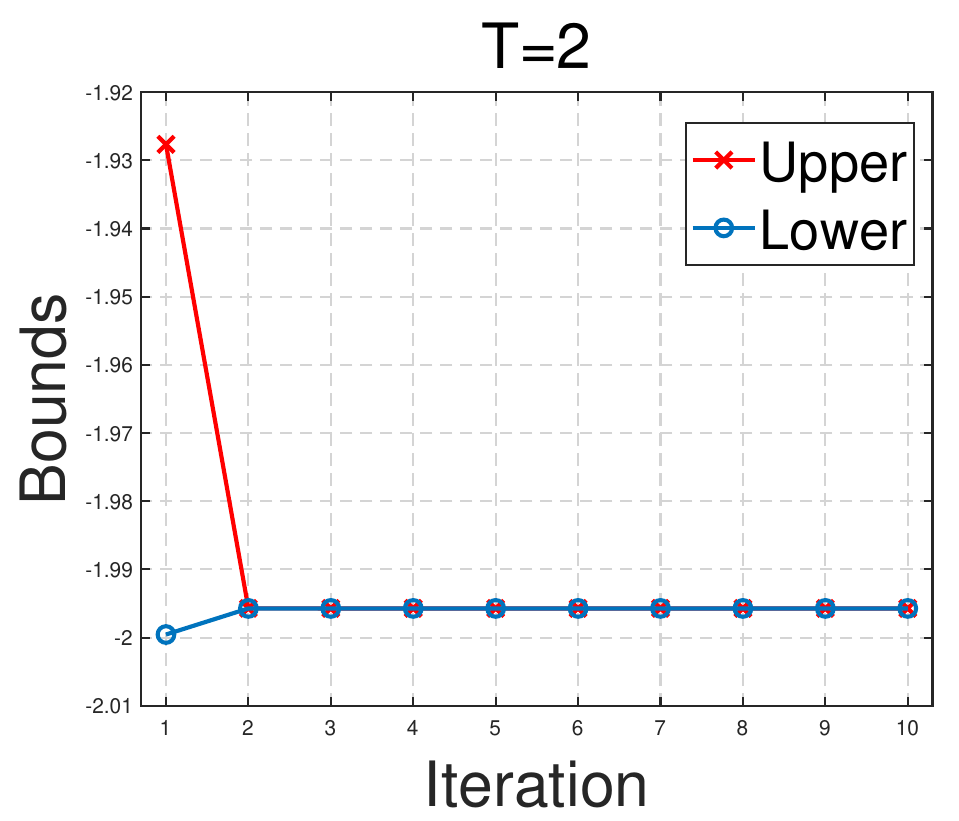}
    }
    \hspace{-1em}
    \subfigure{
    \includegraphics[width=0.23\linewidth]{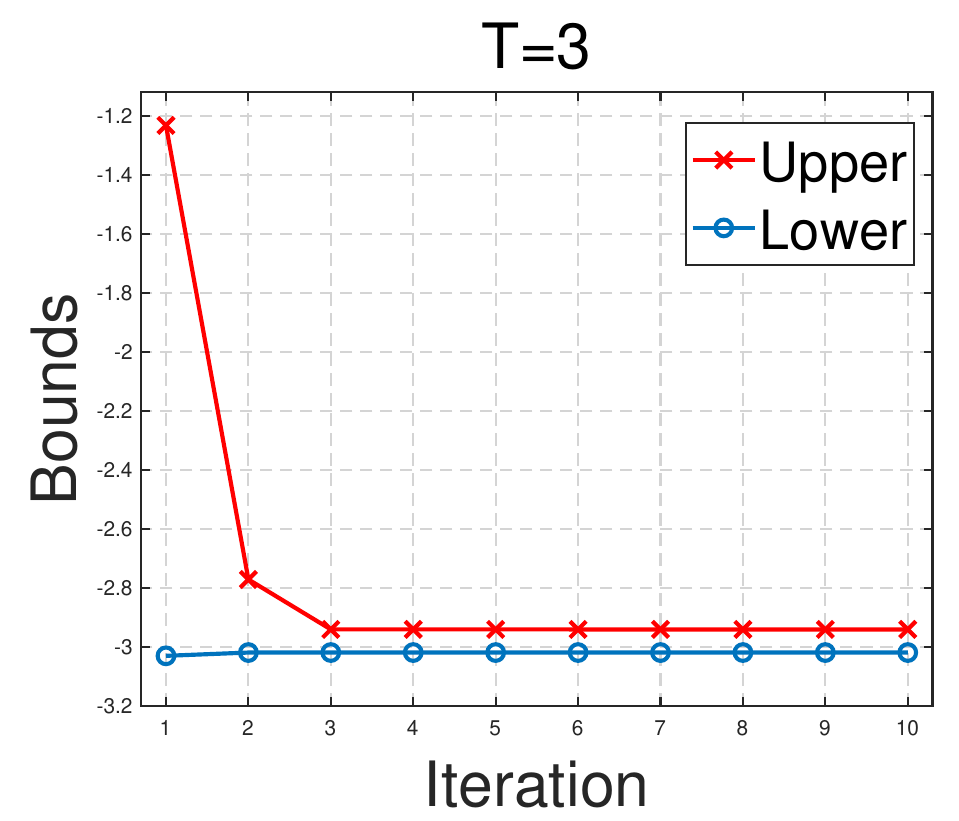}
    }
    \hspace{-1em}
    \subfigure{
    \includegraphics[width=0.23\linewidth]{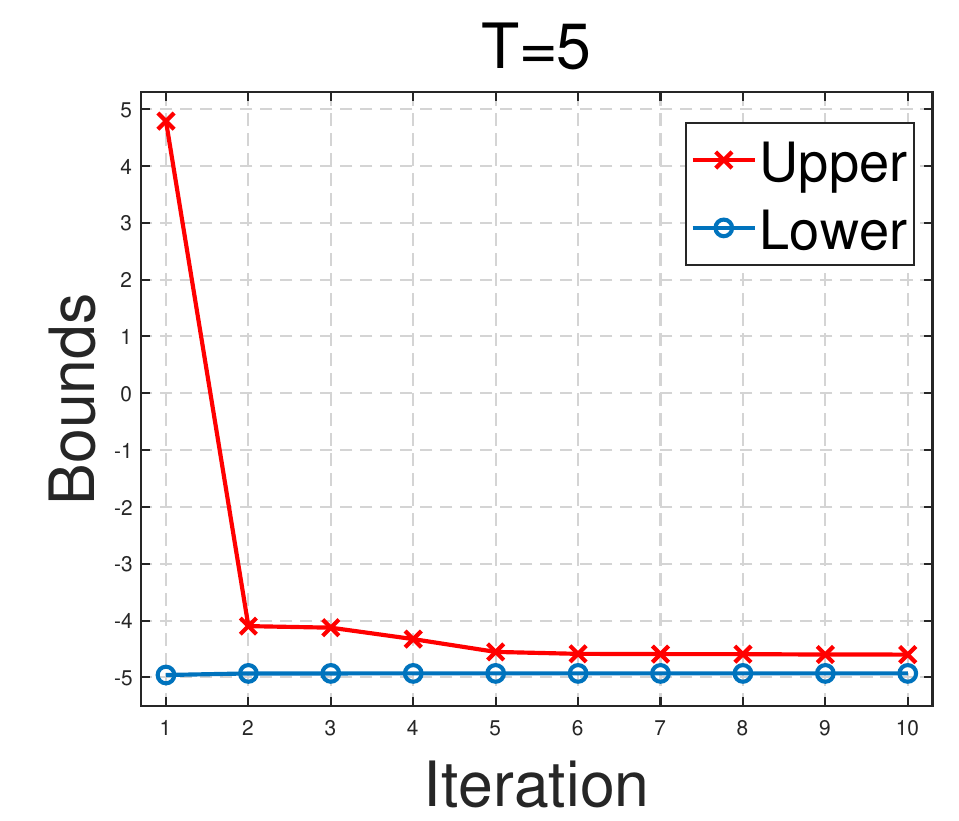}
    }
    \hspace{-1em}
    \subfigure{
    \includegraphics[width=0.23\linewidth]{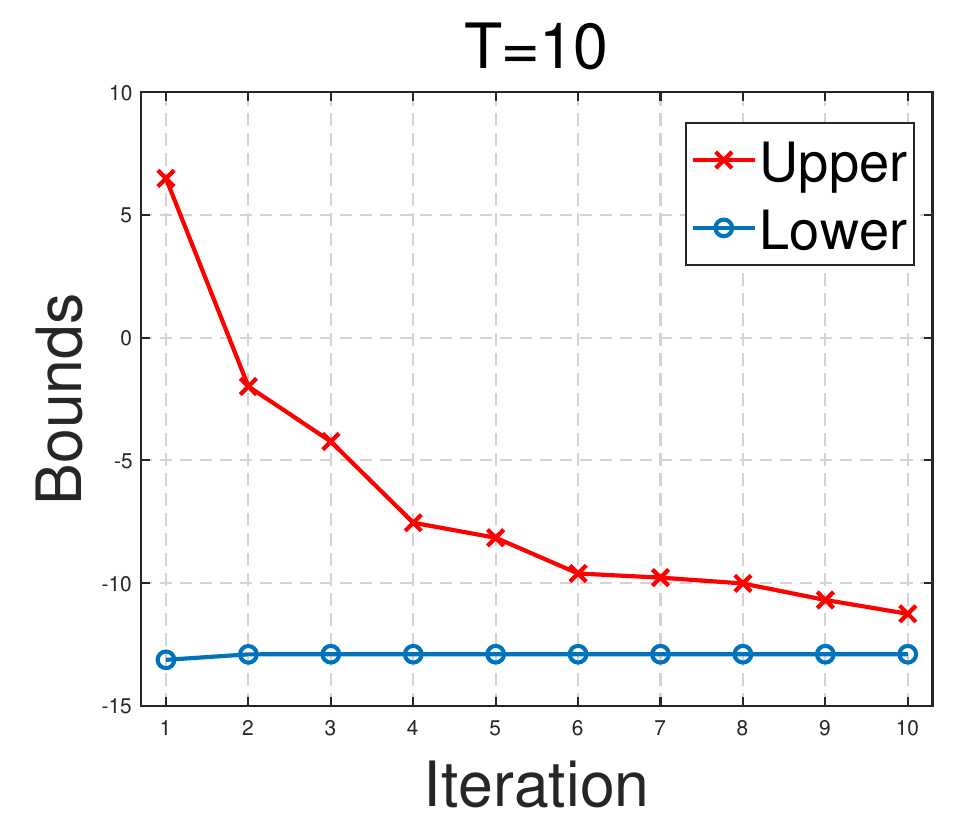}
    }
    \vspace{-0.2cm}
    \captionsetup{font=footnotesize}
    \caption{
    %U
    Convergence of upper and lower bounds %evolution
    for the optimal value of DR-MARSRM.}
    \label{fig:boundsEvolutionDRO} 
\end{figure}

% \begin{table}[ht]
%     \centering 
%     \begin{tabular}{cccc}
%     \hline
%         Stages (T) & lower bound & Upper bound (True) & Upper bound (Statistical) \\
%         \hline
%         2 & -2.0463 & -2.0463 & -1.9497 \\
%         3 & -3.2723 & -3.2723 & -3.2890 \\
%         4 & -4.4177 & - & -4.2436 \\
%         5 & -5.9437 & - & -5.8337 \\
%         10 &  -14.5440 & - & -13.0543 \\
%         15 & -25.5286 & - & -24.4513 \\
%         \hline
%     \end{tabular}
%     \caption{Risk neutral results}
%     \label{tab:bound-risk-neutral}
% \end{table}

\section{Concluding remarks}

{\color{black}
In this paper, we propose a 
multistage distributional 
preference robust risk minimization model where the DM's risk preferences at each stage 
can be described by SRMs but 
display some kind of randomness and the optimal decision is based on the worst-case ARSRM. 
%The ambiguity set of
We demonstrate how the well-known SDDP methods may be effectively used to solve this kind of problem by recursively updating lower and upper bounds of the optimal values at each stage.
A key condition for the SDDP method to be applicable
is stagewise independence of the DM's 
preferences at different stages. Since the ambiguity set at each stage  in the current model is constructed via moment-type conditions, it naturally 
raises a question as to whether 
the conditions 
%at different stages are dependent.
can be updated from the moment-conditions at the previous stage and subsequently how to solve the resulting multistage robust minimization problem. Another issue
which might be of interest is to extend the current robust model 
to the case that the DM's preferences
at each stage are described by random law invariant coherent risk measures (LICRM) to accommodate 
a wider 
%converge 
range of DM's preferences 
and this might be done
%given 
by exploiting the relationship %with 
between LICRM and SRM, see \cite{liu2024stackelberg} for one-stage case. We leave them for future research.

}

\end{document}